\documentstyle{article}
\mag=\magstep1
\def\Bbb#1{\bf#1}
\def\P{{\Bbb P}}
\def\Q{{\Bbb Q}}
\def\R{{\Bbb R}}
\def\Z{{\Bbb Z}}
\def\C{{\Bbb C}}
\def\p{\mbox{\scriptsize\Bbb P}}
\def\O{{\cal O}}
\def\E{{\cal E}}
\def\F{{\cal F}}
\def\G{{\cal G}}
\def\S{{\cal S}}

\def\op{\oplus}
\def\ol#1{\overline{#1}}
\def\proj{{\Bbb Proj}\hspace{2pt}}
\def\Pic{{\rm Pic}}
\def\qed{~ \hbox{\rule{6pt}{8pt}}}
\def\dfrac#1#2{\displaystyle\frac{#1}{#2}}
\def\Label#1{\label{#1}}
\newcounter{thm}
\def\Thm(#1)#2{
\refstepcounter{subsection}\vspace{5pt}\par\noindent(\arabic{section}.\arabic{subsection}) \setcounter{subsubsection}{0} {\bf Theorem.} \setcounter{thm}{\value{subsection}} \Label{thm:#1} {\it #2} \vspace{3pt}
}
\def\Cor(#1)#2{
\refstepcounter{subsubsection}\vspace{5pt}\par\noindent(\arabic{section}.\arabic{subsection}.\arabic{subsubsection}) {\bf Corollary.} \setcounter{thm}{\value{subsubsection}} \Label{thm:#1} {\it #2} \vspace{3pt}
}
\def\Sup(#1)#2{
\refstepcounter{subsection}\vspace{5pt}\par\noindent(\arabic{section}.\arabic{subsection}) \setcounter{subsubsection}{0} {\bf Supplement.} \setcounter{thm}{\value{subsection}} \addtocounter{subsection}{-1} \Label{thm:#1} {#2} \addtocounter{subsection}{1} \vspace{3pt}
}
\def\Titem{
\refstepcounter{subsubsection}\vspace{2pt}\\{\rm(\arabic{section}.\arabic{subsection}.\arabic{subsubsection})~}
}
\def\Eq(#1)#2{
\refstepcounter{subsubsection}\vspace{6pt}\\ (\arabic{section}.\arabic{subsection}.\arabic{subsubsection})~ \hfill $\displaystyle#2$ \qquad\hfill \setcounter{equation}{\value{subsubsection} \Label{eqn:#1} \vspace{6pt}}
}
\def\newpar{
\refstepcounter{subsection}\vspace{5pt}\noindent(\arabic{section}.\arabic{subsection})
}

\def\Longarrow#1#2{\smash{\mathop{\hbox to 8mm{\rightarrowfill}}\limits^{#1}_{#2}}}
\def\Downarrow#1{\hbox{\Huge$\downarrow$}\llap{$\vcenter{\hbox{$\scriptstyle#1\atop{}$\hspace{.8em}}}$}}
\def\Swarrow#1{\hbox{\hspace{.5em}\Huge$\swarrow$}\rlap{$\vcenter{\hbox{\hspace{-.8em}$\scriptstyle#1$}}$}}
\def\ex#1{\Titem {\it Exclusion of No}.#1.}
\def\exs#1{\Titem {\it Exclusion of Nos}.#1.}
\def\cons#1{\Titem {\it Construction from the data of No}.#1.}

\begin{document}

\title{Weak Fano threefolds with del Pezzo fibration}
\author{Kiyohiko TAKEUCHI \\
}
\date{}
\maketitle

\begin{abstract}
    This article treats smooth weak Fano 3-folds $V$ having an extremal ray of type $D$.
    We further assume that the pluri-anti-canonical morphism of $V$ contracts only a finite number of curves, i.e., the anti-canonical model of $V$ is terminal.
    The contraction morphism corresponding to the extremal ray of type $D$ is a del Pezzo fibration of degree $d$ for $1\leq d\leq6$ or $d=8$, $9$.
    Smooth weak Fano 3-folds with an extremal ray of type $D$ of degree $\ne6$ are classified into $47$ deformation types.
\end{abstract}

\section{Introduction}

    Del Pezzo surfaces are one of the most investigated subjects in the classical algebraic geometry.
    A smooth surface $S$ is called a {\it del Pezzo surface} if its anti-canonical divisor $-K_S$ is ample, and it is {\it of degree} $d$ if the self-intersection number $(-K_S)^2$ is $d$.
    The degree $d$ is an integer from 1 to 9, as is well known.
    The extension of this notion to higher dimensions gives us the notion of Fano varieties.

    A smooth 3-fold $V$ is called a {\it Fano 3-fold} if its anti-canonical divisor $-K_V$ is ample.
    Fano 3-folds have been classified by Fano, Iskovskih, Shokurov, Mori, and Mukai.
    We consider 3-folds under weaker conditions.
    Since nef \& big divisors are on the boundary of an ample cone, a weak Fano 3-fold is defined as follows in this paper: a smooth 3-fold $V$ is called a {\it weak Fano 3-fold} if its anti-canonical divisor $-K_V$ is nef \& big.
    Reid\cite{weakFano} has already introduced weak Fano 3-folds which may have canonical singularities, but here we confine ourselves weak Fano 3-folds having no singularities.

    We want to classify weak Fano 3-folds.
    In this article, we treat only those with an extremal ray of type $D$.
    To explain more precisely, we will recall some notions.

    First, we introduce the notion of del Pezzo fibration (cf. \cite{corti}, \cite{fujita}) which is the extension of the notion of del Pezzo surfaces to the relative case.
    A projective morphism $f : V\to C$ from a 3-fold $V$ onto a smooth curve $C$ is called a {\it del Pezzo fibration} if its generic fiber $V_\eta$ is a del Pezzo surface.
    It is {\it of degree} $d$ if $V_\eta$ is of degree $d$.
    Fujita\cite{fujita} investigates these fibrations from the viewpoint of the classification of polarized manifolds.

    Secondly, we recall the extremal ray theory (cf. \cite{extr}).
    For a normal variety $X$, the set of numerical equivalence classes of Cartier divisors
$$N^1(X)=\{\Pic(X)/\equiv\}\otimes\R$$
and the set of numerical equivalence classes of 1-cycles
$$N_1(X)=\{Z_1(X)/\equiv\}\otimes\R$$
are dual vector spaces of finite dimension.
    Let $\ol{NE}(X)\subset N_1(X)$ be the closure of the convex cone generated by numerical equivalence classes of effective 1-cycles of $X$.
    If the canonical divisor $K_X$ is not nef, then $\ol{NE}(X)$ has an extremal ray.
    When $X$ is a smooth 3-fold, the extremal ray falls into one of three types: type E (exceptional (divisorial) type), type C (conic bundle type), and type D (del Pezzo fiber type).
    The contraction morphism $f : X \to Y$ corresponding to the ray of type $D$ is a del Pezzo fibration with relative Picard number 1.

    Finally, the birational mapping called {\it flop} appears in the minimal model theory and plays an important role in it.
    We briefly recall the notion of flop, and 
refer \cite{KM} for more details.
    A {\it $D$-flopping contraction} is a proper birational morphism $f:X\to Y$ between normal varieties $X$ and $Y$ such that the exceptional set ${\rm Ex}(f)$ has codimension at least two in $X$, that $K_X$ is numerically $f$-trivial, and that the $\Q$-Cartier divisor $-(K_X+D)$ is $f$-ample.
    A normal variety $X^+$ together with a proper birational morphism $f^+:X^+\to Y$ is called a {\it $D$-flop} of $f$ if $K_{X^+}+D^+$ is $\Q$-Cartier $f$-ample for $D^+$ the proper transform of $D$ on $X^+$ and if the exceptional set ${\rm Ex}(f^+)$ has codimension at least two in $X^+$.
    Reid\cite{can} treats the $D$-flops of non-singular $3$-folds, and Koll\'ar\cite{flop} proves the existence of $D$-flops for the $D$-flopping contractions of $3$-folds with $\Q$-factorial canonical singularities.
    In the $3$-fold case, ${\rm Ex}(f)$ of the flopping contraction $f:X\to Y$ has dimension $1$, and a curve $C\subset{\rm Ex}(f)$ called a {\it flopping curve}.
    The flopping curve $C\subset X$ satisfies $(C\cdot K_X)=0$.
    The number of the curves $C\subset X$ with $(C\cdot K_X)=0$ is finite.
    For the sake of convenience, we say that a curve $C\subset X$ is {\it $K$-trivial} if $(C\cdot K_X)=0$.

    The aim of this article is to classify weak Fano 3-folds $V$ with an extremal ray of type $D$.
    The 3-folds $V$ have a structure of del Pezzo fibration with relative Picard number 1.
    Under the additional assumption that $V$ has only finite $K$-trivial curves, these 3-folds $V$ are classified into finite types up to deformation.

    Section \ref{sec:main} describes the main result.
    In Section \ref{sec:prel}, we recall the general result that a del Pezzo fibration has a canonical morphism to a (weighted) projective space bundle, and derive three inequalities playing an important role in this paper.
    The extremal rays of type $D$ are classified into the following three subtypes by the lengths of the rays:
$$\begin{array}{r@{~}l}
\mbox{Type $D_1$} :&\mbox{the generic fiber } X_\eta \mbox{ is a del Pezzo surface of degree } d, ~ 1\leq d\leq 6, \\
\mbox{Type $D_2$} :& f : X \to Y ~\mbox{ is a quadric bundle, \qquad and} \\
\mbox{Type $D_3$} :& f : X \to Y ~\mbox{ is a $\P^2$-bundle.}
\end{array}$$
    Thus we divide the argument into several cases according to the above types.
    In Section \ref{sec:quad}, we deal with the case of type $D_2$, and in Sections \ref{sec:df2} to \ref{sec:df1}, the case of type $D_1$ of $d\ne0$.

    The author is grateful to Prof. Shigeru Mukai for pointing out various mistakes, to Prof. Shigefumi Mori for suggesting some ideas and correcting several missing terminology, and to Prof. Takashi Maeda and Prof. Hiromichi Takagi for valuable conversation.
    This research was partially supported by the Ministry of Education, Science, Sports and Culture, Grant-in-Aid for Scientific Research (C)(2), 12640048, 2000-2001.

\section{Results} \label{sec:main}
\renewcommand{\arraystretch}{1.3}
\def\vd{\:\lower-1.ex\hbox{\rm :}\hskip-0.65ex\lower.8ex\hbox{\rm :}\:}

\newpar
    In this section, we state the theorems in this paper.
    We consider a weak Fano $3$-fold $V$ such that the number of $K$-trivial curves is finite.
    Namely, we consider a smooth projective $3$-fold $V$ such that its pluri-anti-canonical linear system $|-mK_V|$ is free and that the number of curves contracted by the corresponding morphism $\phi_{|-mK_V|}$ is finite.
    In this case, the image $\ol{V}=\phi_{|-mK_V|}(V)$ of this morphism, the anti-canonical model of $V$, is a projective $3$-fold with only terminal singularities.
    We additionally assume that $V$ has an extremal ray of type $D$.
    Then there is a contraction morphism $\varphi : V\to Y$, which is a del Pezzo fibration of degree $d$ with relative Picard number $1$.
    The curve $Y$ is isomorphic to the projective line $\P^1$, since $V$ is a weak Fano $3$-fold and $H^1(V,\O_V)=0$.

    We can classify these weak Fano $3$-folds as in (\ref{thm:p2})--(\ref{thm:deform}).
    We use the following notation and abbreviation.
    Let $\E=\bigoplus^n_{i=0}\O(a_i)$ be a locally free sheaf of rank $n+1$ on $\P^1$, and $\P(\E)$ the corresponding $\P^n$-bundle over $\P^1$, denoted by $\P[a_0,a_1,\dots,a_n]$.
    Let $H$ and $F$ be the tautological divisor and a fiber of $\pi:\P[a_0,a_1,\dots,a_n] \to\P^1$, respectively; let $H_V$, $F_V$ be the restrictions of $H$, $F$ to a subvariety $V\subset X$.
    On the tables in Theorems(\ref{thm:quad})--(\ref{thm:df5}), the symbol $(a,b)$ indicates a member of the linear system $|aH+bF|$ on $\P[a_0,a_1,\dots,a_n]$; the same symbol also means a member of the linear system $|aH_V+bF_V|$ on $V$.
    The column ``$e$'' in the tables is the number of $K$-trivial curves on general member $V$; the column ``$R'$'' is the type of the extremal ray on the flop $V'$ of $V$; the column ``$V'\to{\rm cont}_{R'}(V')$'' is the contraction morphism of the ray, and ``$S_d$-fibration'' in the column means a del Pezzo fibration of degree $d$.
\def\dr[#1](#2)(#3)#4#5#6#7{
\refstepcounter{subsubsection}(\arabic{section}.\arabic{subsection}.\arabic{subsubsection})&\P[#1]\supset(#2)&(#3)&#4&#5&#6&#7\\}

    As mentioned in Introduction, the extremal rays of type $D$ are classified into three subtypes: of type $D_1$, $D_2$, or $D_3$.
    It is natural to treat three cases according to the subtypes.
%
%
    The case of type $D_3$ is very easy (see (\ref{subsubsec:p2})) and we get

\Thm(p2){
    Every weak Fano $3$-fold $V$ with only finite $K$-trivial curves, such that $V$ has an extremal ray of type $D_3$, is isomorphic to one of the following $\P^2$-bundle over $\P^1$:
\Titem
$V=\P[0^3]$,
which is a Fano $3$-fold isomorphic to $\P^2\times\P^1$;
\Titem
$V=\P[0^2,1]$,
which is a Fano $3$-fold obtained by the blowing-up of $\P^3$ along a line $Z\subset\P^3$; \quad and
\Titem
$V=\P[0,1^2]$,
which is a weak Fano $3$-fold with only one flopping curve $s_0$, and the $(-F)$-flop $V'$ has an extremal ray of type D${}_3$, i.e., $V'$ is again $\P^2$-bundle over $\P^1$.
\vspace{3pt}\\ {}\quad
In each case, $(-K_V)^3=54$.
}


    In the case of type $D_2$, the contraction morphism $\varphi : V\to \P^1$ is a del Pezzo fibration of degree $8$, i.e., a quadric bundle over $\P^1$.
    In this case, we get

\Thm(quad){
    Every weak Fano $3$-fold $V$ with only finite $K$-trivial curves, such that $V$ has an extremal ray of type $D_2$, is isomorphic to a general member of the linear system $\Lambda$ on the $\P^3$-bundle $X=\P[a_0,\dots,a_3]$ over $\P^1$, where $X$ and $\Lambda$ are one of the following:
$$\begin{array}{ccc@{\vd}c@{}ccl}
&X\supset V\in\Lambda&-K_V&(-K_V)^3&e&R'&V'\to{\rm cont}_{R'}(V')\\
\dr[0^3,1](2,0)(2,1){40}{0}{E_1}{V\!\!=\!\!V'\to\Q^3\supset {}_0C_2}
\dr[0,1^3](2,-1)(2,0){40}{1}{C_2}{\P^1\mbox{-bundle}\bigm/\P^2}
\dr[0^4](2,1)(2,1){32}{0}{E_1}{V\!\!=\!\!V'\to\P^3\supset {}_1C_4}
\dr[0^2,1^2](2,0)(2,0){32}{2}{D_2}{\Q^2\mbox{-bundle}\bigm/\P^1}
\dr[0^3,1](2,1)(2,0){24}{4}{E_2}{V'\to B_4\ni P}
\dr[0,1^3](2,0)(2,-1){24}{0}{C_1}{V\!\!=\!\!V'\to\P^2\supset\Delta_4}
\dr[0^4](2,2)(2,0){16}{8}{D_2}{\Q^2\mbox{-bundle}\bigm/\P^1}
\dr[0^2,1^2](2,1)(2,-1){16}{1}{D_1}{S_4\mbox{-fibration}\bigm/\P^1}
\dr[0,1^3](2,1)(2,-2){8}{18}{D_2}{\Q^2\mbox{-bundle}\bigm/\P^1}
\end{array}$$
}

\Sup(quads){
    In the above list, the $3$-fold $V$ has the following properties in each case.
    Let $H$ and $F$ be the tautological divisor and a fiber of $X\to\P^1$, and $H_V$ and $F_V$ their restrictions to $V$.
    If $e>0$, then $V$ has the $(-F_V)$-flop $V'$, and $H_{V'}$, $F_{V'}$ denote the birational transforms of $H_V$, $F_V$.
\Titem
    $(-K_V)^3=40$; $V$ is the Fano $3$-fold which is the blowing-up of a quadric $3$-fold $\Q^3$ along a smooth conic; $V\to\Q^3$ is defined by $|H_V|$, and the exceptional divisor is the unique member of $|H_V-F_V|$; $V$ is No.$29$ on Table 3 in \cite{MMT}.
\Titem
    $(-K_V)^3=40$; $V$ has only one flopping curve $s_0$ as a section; the $(-F_V)$-flop $V'$ has a $\P^1$-bundle structure over $\P^2$; $V'$ contains the strict transform of $s_0$ as a subsection; $V'\to\P^2$ is defined by $|H_{V'}-F_{V'}|$; the anti-canonical model $\ol{V}$ is a singular del Pezzo $3$-fold of degree $5$ with only one ODP.
\Titem
    $(-K_V)^3=32$; $V$ is the Fano $3$-fold which is the blowing-up of $\P^3$ along a quartic elliptic curve; $V\to\P^3$ is defined by $|H_V|$, and the exceptional divisor is the unique member of $|2H_V-F_V|$; $V$ is No.$25$ on Table 3 in \cite{MMT}.
\Titem
    $(-K_V)^3=32$; general $V$ has exactly two flopping curves as sections; the $(-F_V)$-flop $V'$ has a quadric bundle structure again; $V'\to\P^1$ is defined by $|H_{V'}-F_{V'}|$; $\ol{V}$ is a complete intersection of two quadrics in $\P^5$, and has two ODP's.
\Titem
    $(-K_V)^3=24$; general $V$ has exactly four flopping curves as sections; the $(-F_V)$-flop $V'$ is the blowing-up at a point of $B_4$, where $B_4$ is a del Pezzo $3$-fold of degree $4$, a complete intersection of two quadrics in $\P^5$; $V'\to B_4$ defined by $|2H_{V'}-F_{V'}|$, and the exceptional divisor is a unique member of $|H_{V'}-F_{V'}|$; $\ol{V}$ is a cubic $3$-fold with $4$ ODP's.
\Titem
    $(-K_V)^3=24$; $V$ is a Fano $3$-fold and has a conic bundle structure over $\P^2$ with the discriminant locus of degree $4$; $V\to\P^2$ is defined by $|H_V-F_V|$; $V$ is No.$18$ on Table 3 in \cite{MMT}.
\Titem
    $(-K_V)^3=16$; general $V$ has exactly eight flopping curves as sections; the $(-F_V)$-flop $V'$ has a quadric bundle structure again; $V'\to\P^1$ is defined by $|2H_{V'}-F_{V'}|$; $\ol{V}$ is a double covering of $\P^3$ branched along the quartic surface with $8$ ODP's.
\Titem\Label{subsubsec:quad}
    $(-K_V)^3=16$; $V$ has only one flopping curve $C$ as a bisection; the $(-F_V)$-flop $V'$ has a structure of del Pezzo fibration of degree $4$ over $\P^1$; $V'$ contains the strict transform of $C$ as a section; $V'\to\P^1$ is defined by $|H_{V'}-F_{V'}|$; $\ol{V}$ is a singular Fano $3$-fold of index $1$ of genus $9$ with only one ODP.
\Titem
    $(-K_V)^3=8$; general $V$ has exactly eighteen flopping curves $C_1$, $C_2$, $s_1, \dots, s_{16}$; $C_1$ and $C_2$ are bisections and the $s_1, \dots, s_{16}$ are sections of the fibration; the $(-F_V)$-flop $V'$ has a quadric bundle structure again; $\ol{V}$ is a hypersurface of degree $6$ with $18$ ODP's in the weighted projective space $\P(1^3,2,3)$.
}
\\
Theorem(\ref{thm:quad}) and Supplement(\ref{thm:quads}) are treated in Section \ref{sec:quad}.


    In the case of type $D_1$, $\varphi : V\to \P^1$ is a del Pezzo fibration of degree $d$, $1\leq d\leq 6$, over $\P^1$.
    In this paper, we treat the cases of degree $1$, $2$, $3$, $4$, and $5$, and give Theorems(\ref{thm:df1}), (\ref{thm:df2}), (\ref{thm:df3}), (\ref{thm:df4}), and (\ref{thm:df5}), respectively. 
    These results are derived in Sections \ref{sec:df1}, \ref{sec:df2}, \ref{sec:df3}, \ref{sec:df4}, and \ref{sec:df5}, respectively.


\Thm(df1){
    Every weak Fano $3$-fold $V$ with only finite $K$-trivial curves, such that $V$ has an extremal ray of type $D_1$ of degree $1$, is isomorphic to a general member of the linear system $\Lambda$ on the weighted projective space bundle $\proj(\S)$ with weights $(1^2,2,3)$ over $\P^1$, where the graded $\O_{\p^1}$-algebra $\S$ and the linear system $\Lambda$ are one of the following:
$$\begin{array}{ccc@{\vd}cccl}
&\proj(\S)\supset V\in\Lambda&-K_V&(-K_V)^3&e&R'&V'\to{\rm cont}_{R'}(V')\\
\dr[0,1;0;0](6,0)(1,1){4}{0}{E_1}{V\!\!=\!\!V'\to B_1\supset {}_1C_1}
\dr[0^2;-2;-3](6,6)(1,0){2}{1}{D_1}{S_1\mbox{-fibration}\bigm/\P^1}
\end{array}$$
    Here, $\P[a_0,a_1;a_2;a_3]$ denotes the weighted projective space bundle $\proj(\S)$ over $\P^1$ for the graded $\O_{\p^1}$-algebra $\S = \bigoplus_{d \geq 0} \S_d$ such that $\S_d = \bigoplus_{i+2j+3k=d} S^i(\O(a_0)\op\O(a_1))\otimes S^j(\O(a_2))\otimes S^k(\O(a_3))$, and, $(a,b)$ denotes a member of $|aH+bF|$ on $\proj(\S)$ (or the restriction to $V$) for the tautological $\Q$-divisor $H$ and a fiber $F$ of $\pi : \proj(\S)\to\P^1$ containing $\varphi : V\to\P^1$.
}

\Sup(df1s){
    In the above list, the $3$-fold $V$ has the following properties in each case.
\Titem
    $(-K_V)^3 = 4$; $V$ is the Fano $3$-fold which is the blowing-up of $B_1$ along an elliptic curve $C$ with $(-K_{B_1}\cdot C) = 2$, where $B_1$ is a Fano $3$-fold of index $2$, a hypersurface of degree $6$ in the weighted projective space $\P(1^3,2,3)$; $V$ is No.$1$ on Table 3 in \cite{MMT}.
\Titem
    $(-K_V)^3 = 2$; $V$ has only one flopping curve as a section; the $(-F_V)$-flop $V'$ has a structure of del Pezzo fibration of degree $1$; the (pluri-)anti-canonical model $\ol{V}$ is the complete intersection of the smooth del Pezzo 4-fold of degree 1 and a singular weighted hypersurface of degree 2 in the weighted projective space $\P(1^4,2,3)$; $\ol{V}$ has only one ODP.}


\Thm(df2){
    Every weak Fano $3$-fold $V$ with only finite $K$-trivial curves, such that $V$ has an extremal ray of type $D_1$ of degree $2$, is isomorphic to a double covering of the $\P^2$-bundle $\P(\E)$ over $\P^1$ with branch locus $B$, where $\E$ and $B$ are one of the following:
$$\begin{array}{ccc@{\vd}cccl}
&\P(\E)\supset B&-K_V&(-K_V)^3&e&R'&V'\to{\rm cont}_{R'}(V')\\
\dr[0^2,1](4,0)(1,1){8}{0}{E_1}{V\!\!=\!\!V'\to B_2\supset {}_1C_2}
\dr[0^3](4,2)(1,1){6}{0}{C_1}{V\!\!=\!\!V'\to\P^2\supset\Delta_8}
\dr[0,1,2](4,-2)(1,0){6}{1}{E_3}{V'\to W'\ni P}
\dr[0,1^2](4,0)(1,0){4}{2}{D_1}{S_2\mbox{-fibration}\bigm/\P^1}
\end{array}$$
}

\Sup(df2s){
    In the above list, the $3$-fold $V$ has the following properties in each case.
\Titem
    $(-K_V)^3=8$; $V$ is the Fano $3$-fold which is the blowing-up of $B_2$ along an elliptic curve $C$ such that $(C\cdot -K_{B_2})=4$, where $B_2$ is a Fano $3$-fold of index $2$, a double covering of $\P^3$ branched along a smooth quartic surface; $V$ is No.$3$ on Table 3 in \cite{MMT}.
\Titem
    $(-K_V)^3=6$; $V$ is a Fano $3$-fold and has a conic bundle structure over $\P^2$ with the discriminant locus of degree $8$; $V$ is No.$2$ on Table 3 in \cite{MMT}.
\Titem
    $(-K_V)^3=6$; $V$ has only one flopping curve as a section; the $(-F_V)$-flop $V'$ has a morphism to $W'$ corresponding the extremal ray of type $E_3$; $W'$ is a singular Fano $3$-fold $B_1'$ of index $2$ of degree $1$ with only one double point; the anti-canonical model $\ol{V}$ is the double covering of the cone $C(\Sigma_1)$ over the rational ruled surface $\Sigma_1\cong\P[0,1]$ branched along $B$, where $B$ is an irreducible component of the quartic hypersurface section $B+F_1+F_2$ through the vertex of $C(\Sigma_1)\subset\P^5$, and $F_i=C(f)\subset C(\Sigma_1)$ is the cone over a fiber $f$ of $\Sigma_1$; $\ol{V}$ has only one ODP.
\Titem
    $(-K_V)^3=4$; general $V$ has exactly two flopping curves as sections; the $(-F_V)$-flop $V'$ has a structure of del Pezzo fibration of degree $2$ over $\P^1$ again; $\ol{V}$ is the double covering of the singular quadric $3$-fold $\Q_0\subset\P^4$ branched along a quartic hypersurface section on $\Q_0$; $\ol{V}$ has two ODP's.
}


\Thm(df3){
    Every weak Fano $3$-fold $V$ with only finite $K$-trivial curves, such that $V$ has an extremal ray of type $D_1$ of degree $3$, is isomorphic to a general member of the linear system $\Lambda$ on the $\P^3$-bundle $\P(\E)$ over $\P^1$, where $\E$ and $\Lambda$ are one of the following:
$$\begin{array}{ccc@{\vd}cccl}
&\P(\E)\supset V\in\Lambda&-K_V&(-K_V)^3&e&R'&V'\to{\rm cont}_{R'}(V')\\
\dr[0^3,1](3,0)(1,1){12}{0}{E_1}{V\!\!=\!\!V'\to B_3\supset {}_1C_3}
\dr[0^4](3,1)(1,1){10}{0}{E_1}{V\!\!=\!\!V'\to\P^3\supset {}_{10}C_9}
\dr[0,1^2,2](3,-2)(1,0){10}{1}{E_1}{V'\to B_2\supset {}_0C_1}
\dr[0,1^3](3,-1)(1,0){8}{1}{C_1}{V'\to\P^2\supset\Delta_7}
\dr[0^2,1^2](3,0)(1,0){6}{3}{D_1}{S_3\mbox{-fibration}\bigm/\P^1}
\dr[0^3,1](3,1)(1,0){4}{9}{E_5}{V'\to W'\ni P}
\dr[0^4](3,2)(1,0){2}{27}{D_1}{S_3\mbox{-fibration}\bigm/\P^1}
\end{array}$$
}

\Sup(df3s){
    In the above list, the $3$-fold $V$ has the following properties in each case.
\Titem
    $(-K_V)^3=12$; $V$ is the Fano $3$-fold which is the blowing-up of a cubic $3$-fold $B_3$ along a smooth plane elliptic curve; $V$ is No.$5$ on Table 3 in \cite{MMT}.
\Titem
    $(-K_V)^3=10$; $V$ is the Fano $3$-fold which is the blowing-up of $\P^3$ along a smooth curve of degree $9$ and of genus $10$; $V$ is No.$4$ on Table 3 in \cite{MMT}.
\Titem
    $(-K_V)^3=10$; $V$ has only one flopping curve as a section; the $(-F_V)$-flop $V'$ is the blowing-up of $B_2$ along a line; the anti-canonical model $\ol{V}$ is a divisor of the cone $C(\P[0^2,1])\subset\P^7$ over the $\P^2$-bundle $\P[0^2,1]$, and is the irreducible component of a cubic hypersurface $\ol{V}+P_1+P_2$ through the vertex, where $P_i$ is the cone $C(f)\subset C(\P[0^2,1])$ over a fiber $f$ of $\P[0^2,1]\to\P^1$; the singularity of $\ol{V}$ is one ODP, the vertex of the cone.
\Titem
    $(-K_V)^3=8$; $V$ has only one flopping curve as a section; the $(-F_V)$-flop $V'$ has a conic bundle structure over $\P^2$ with the discriminant locus of degree $7$; $\ol{V}$ is embedded in the cone $C(\P^2\times\P^1)\subset\P^6$, and is the irreducible component of a cubic hypersurface $\ol{V}+P$ through the vertex, where $P=C(\P^2)\subset C(\P^2\times\P^1)$ is the cone over a fiber of $\P^2\times\P^1\to\P^1$; the singularity of $\ol{V}$ is one ODP, the vertex of the cone.
\Titem
    $(-K_V)^3=6$; general $V$ has exactly three flopping curves as sections; the $(-F_V)$-flop $V'$ has a structure of del Pezzo fibration of degree $3$ over $\P^1$ again; $\ol{V}\subset\P^5$ is the complete intersection of a cubic $4$-fold and a singular quadric $4$-fold having singularities along a line; the singularities of $\ol{V}$ consist three points, the intersection of the cubic $4$-fold and the singular locus of the singular quadric $4$-fold.
\Titem
    $(-K_V)^3=4$; general $V$ has exactly nine flopping curves as sections; the $(-F_V)$-flop $V'$ has a morphism to a singular $3$-fold corresponding the extremal ray of type $E_5$; the anti-canonical model $\ol{V}$ is a singular quartic $3$-fold in $\P^4$ which contains a plane and has $9$ ODP's on the plane.
\Titem
    $(-K_V)^3=2$; general $V$ has exactly twenty-seven flopping curves as sections; the $(-F_V)$-flop $V'$ has a structure of del Pezzo fibration of degree $3$ over $\P^1$ again; $\ol{V}$ is the double covering of $\P^3$ branched along the sextic surface $S$ defined by $C_1^2-4C_0C_2$ for three cubic forms $C_i$ with $27$ ODP's corresponding to $C_0=C_1=C_2=0$; the singularities of $\ol{V}$ correspond to the singularities of $S$.
}


\Thm(df4){
    Every weak Fano $3$-fold $V$ with only finite $K$-trivial curves, such that $V$ has an extremal ray of type $D_1$ of degree $4$, is isomorphic to a complete intersection $W_1 \cap W_2$ of general members $W_1$ and $W_2$ of the linear systems $\Lambda_1$ and $\Lambda_2$, respectively, on the $\P^4$-bundle $\P(\E)$ over $\P^1$, where $\E$, $\Lambda_1$ and $\Lambda_2$ are one of the following:
\def\dr[#1](#2)(#3)(#4)#5#6#7#8{
\refstepcounter{subsubsection}(\arabic{section}.\arabic{subsection}.\arabic{subsubsection})&\P[#1]\!\supset\!(#2)\!\cap\!(#3)&(#4)&#5&#6&#7&#8\\}
$$\begin{array}{l@{\!\!}c@{\!}c@{\!\!\vd\!\!}c@{\!\!}c@{\:\:}c@{\;}l}
&\P(\E)\supset V=W_1\cap W_2&-K_V&(-K_V)^3&e&R'&V'\to{\rm cont}_{R'}(V')\\
\dr[0^4,1](2,0)(2,0)(1,1){16}{0}{E_1}{V\!\!=\!\!V'\to B_4\supset {}_1C_4}
\dr[0,1^2,2^2](2,-2)(2,-2)(1,0){16}{1}{D_2}{\Q^2\mbox{-bundle}\bigm/\P^1}
\dr[0^5](2,0)(2,1)(1,1){14}{0}{E_1}{V\!\!=\!\!V'\to\Q^3\supset {}_5C_8}
\dr[0,1^3,2](2,-2)(2,-1)(1,0){14}{1}{E_1}{V'\to B_3\supset {}_0C_2}
\dr[0,1^4](2,-1)(2,-1)(1,0){12}{1}{E_1}{V'\to \P^3\supset {}_7C_8}
\dr[0^2,1^3](2,-1)(2,0)(1,0){10}{4}{C_1}{V'\to\P^2\supset\Delta_6}
\dr[0^3,1^2](2,0)(2,0)(1,0){8}{4}{D_1}{S_4\mbox{-fibration}\bigm/\P^1}
\dr[0^4,1](2,0)(2,1)(1,0){6}{8}{E_3}{V'\to W'\ni P}
\dr[0^5](2,0)(2,2)(1,0){4}{16}{D_1}{S_4\mbox{-fibration}\bigm/\P^1}
\dr[0^5](2,1)(2,1)(1,0){4}{16}{D_1}{S_4\mbox{-fibration}\bigm/\P^1}
\dr[0,1^4](2,-1)(2,0)(1,-1){2}{36}{D_1}{S_4\mbox{-fibration}\bigm/\P^1}
\end{array}$$
}

\Sup(df4s){
    In the above list, the $3$-fold $V$ has the following properties in each case.
\Titem
    $(-K_V)^3=16$; $V$ is the Fano $3$-fold which is the blowing-up of $B_4$ along a quartic elliptic curve; $V$ is No.$10$ on Table 3 in \cite{MMT}.
\Titem
    $(-K_V)^3=16$; $V$ has only one flopping curve as a section; the $(-F_V)$-flop $V'$ has a quadric bundle structure(this is the inverse of (\ref{subsubsec:quad})).
\Titem
    $(-K_V)^3=14$; $V$ is the Fano $3$-fold which is the blowing-up of a quadric $3$-fold $\Q^3$ along a smooth curve of degree $8$ and of genus $5$; $V$ is No.$8$ on Table 3 in \cite{MMT}.
\Titem
    $(-K_V)^3=14$; $V$ has only one flopping curve as a section; the $(-F_V)$-flop $V'$ is the blowing-up of a cubic $3$-fold $B_3$ along a smooth conic; the anti-canonical model $\ol{V}$ is embedded in $C(\P[0^3,1])\subset\P^9$, and is the intersection $\ol{W_1}\cap\ol{W_2}$ of the divisors $\ol{W_i}\subset C(\P[0^3,1])$ through the vertex, where $\ol{W_1}+P_1+P_2$, $\ol{W_2}+P_3$ are quadric hypersurface sections of $C(\P[0^3,1])$ for the cone $P_i=C(\P^3)\subset C(\P[0^3,1])$ over a fiber of $\P[0^3,1]\to\P^1$; the singularity of $\ol{V}$ is one ODP corresponding to the vertex of $C(\P[0^3,1])$.
\Titem
    $(-K_V)^3=12$; $V$ has only one flopping curve as a section; the $(-F_V)$-flop $V'$ is the blowing-up of $\P^3$ along a smooth curve of degree $8$ and of genus $7$; $\ol{V}$ is embedded in $C(\P^3\times\P^1)\subset\P^8$, and is the intersection $\ol{V}=\ol{W_1}\cap\ol{W_2}$ of the divisors $\ol{W_i}$ of $C(\P^3\times\P^1)$ through the vertex, where $\ol{W_i}+P$ is a quadric hypersurface section of $C(\P^3\times\P^1)$ for the cone $P=C(\P^3)\subset C(\P^3\times\P^1)$ over a fiber of $\P^3\times\P^1\to\P^1$; the singularity of $\ol{V}$ is one ODP corresponding to the vertex of $C(\P^3\times\P^1)$.
\Titem
    $(-K_V)^3=10$; general $V$ has exactly four flopping curves as sections; the $(-F_V)$-flop $V'$ has a conic bundle structure over $\P^2$ with the discriminant locus of degree $6$; $\ol{V}=Gr(5,2)\cap H_1\cap H_2\cap Q$ for two hyperplanes $H_1$, $H_2$ and a quadric hypersurface $Q$, where $Gr(5,2)\cap H_1$ has a singular locus $\Sigma\cong\P^2$, and $H_2$ and $Q$ are general; the singularities of $\ol{V}$ are two ODP's, the intersection $\Sigma\cap H_2\cap Q$.
\Titem
    $(-K_V)^3=8$; general $V$ has exactly four flopping curves as sections; the $(-F_V)$-flop $V'$ has a structure of del Pezzo fibration of degree $4$ over $\P^1$ again; $\ol{V}=Q_0\cap Q_1\cap Q_2\subset\P^6$ for three quadrics, only $Q_0$ has $2$-dimensional singular locus $\Sigma=\P^2$ and the others are smooth; the singularities of $\ol{V}$ consist of four ODP's, the intersection $\Sigma\cap Q_1\cap Q_2$.
\Titem
    $(-K_V)^3=6$; general $V$ has exactly eight flopping curves as sections; the $(-F_V)$-flop $V'$ has a morphism to a singular $3$-fold corresponding the extremal ray of type $E_3$ (or $E_4$); $\ol{V}=\ol{W_1}\cap \ol{W_2}\subset\P^5$ for a smooth quadric $\ol{W_1}$ and a singular cubic $4$-fold $\ol{W_2}$ containing $\P^3$ with $8$ singular points on the $\P^3$; the singularities of $\ol{V}$ come from the singularities of $\ol{W_2}$.
\Titem
    $(-K_V)^3=4$; general $V$ has exactly sixteen flopping curves as sections; the $(-F_V)$-flop $V'$ has a structure of del Pezzo fibration of degree $4$ over $\P^1$ again; $\ol{V}$ is the double covering of $\Q^3\subset\P^4$ branched along a surface of degree $8$; the singularities of $\ol{V}$ consist of $16$ ODP's coming from the singularities of the branching surface.
\Titem
    $(-K_V)^3=4$; general $V$ has exactly sixteen flopping curves as sections; the $(-F_V)$-flop $V'$ has a structure of del Pezzo fibration of degree $4$ over $\P^1$ again; $\ol{V}$ is the quartic $3$-folds defined by $Q_0Q_3-Q_1Q_2=0$  for $4$ quadrics; the singularities of $\ol{V}$ consist of $16$ ODP's defined by $Q_0=Q_1=Q_2=Q_3=0$.
\Titem
    $(-K_V)^3=2$; general $V$ has exactly thirty-six flopping curves $s_1, \dots, s_{32}$, $C_1, \dots, C_4$; $C_1, \dots, C_4$ are bisections and $s_1, \dots, s_{32}$ are sections of the fibration; the $(-F_V)$-flop $V'$ has a structure of del Pezzo fibration of degree $4$ over $\P^1$ again; $\ol{V}$ is the double covering of $\P^3$ branched along a sextic surface $S$ having $36$ ODP's; the singularities of $\ol{V}$ come from the singularities of $S$.
}


\Thm(df5){
    Every weak Fano $3$-fold $V$ with only finite $K$-trivial curves, such that $V$ has an extremal ray of type $D_1$ of degree $5$, is isomorphic to the subvariety defined by the Pfaffians of the $4\times4$ diagonal minors of a $5\times5$ skew-symmetric matrix $M = (m_{ij})$ with $m_{ij} \in H^0(\P(\E), \O(H + w_{ij}F))$ for $1 \leq i < j \leq 5$ on $\P^5$-bundle $\P(\E)$ over $\P^1$.
    Here $H$ and $F$ are the tautological divisor and a fiber of the $\P^5$-bundle $\pi : \P(\E)\to\P^1$ containing $\varphi : V\to\P^1$, respectively, and there is a sequence of five integers $k_1 \leq k_2 \leq k_3 \leq k_4 \leq k_5$ such that $w_{ij} = (k_1 + k_2 + k_3 + k_4 + k_5)/2 - k_i - k_j$.
    More precisely, $\E$ and the sequence of integers $k_i$ are one of the following:
\def\dr[#1](#2)(#3)#4#5#6#7{
\refstepcounter{subsubsection}(\arabic{section}.\arabic{subsection}.\arabic{subsubsection})&\P[#1]&(#2)&(#3)&#4&#5&#6&#7\\}
$$\begin{array}{l@{}c@{\:\:}c@{\:\:}c@{\vd}c@{}ccl}
&\P(\E)&(k_1,\dots,k_5)&-K_V&(-K_V)^3&e&R'&V'\to{\rm cont}_{R'}(V')\\
\dr[0,1^2,2^3](-3^2,-2^3)(1,0){22}{1}{C_2}{\P^1\mbox{-bundle}\bigm/\P^2}
\dr[0^5,1](0^5)(1,1){20}{0}{E_1}{V\!\!=\!\!V'\to B_5\supset {}_1C_5}
\dr[0,1^4,2](-2^3,-1^2)(1,0){18}{1}{E_1}{V'\to B_4\supset {}_0C_3}
\dr[0,1^5](-2,-1^4)(1,0){16}{1}{E_1}{V'\to\Q^3\supset {}_3C_7}
\dr[0^2,1^4](-1^4,0)(1,0){14}{2}{E_1}{V'\to\P^3\supset {}_4C_7}
\dr[0^3,1^3](-1^2,0^3)(1,0){12}{3}{C_1}{V'\to\P^2\supset\Delta_5}
\dr[0^4,1^2](0^5)(1,0){10}{5}{D_1}{S_5\mbox{-fibration}\bigm/\P^1}
\dr[0^5,1](0^3,1^2)(1,0){8}{8}{E_1}{V'\to V_{12}\supset {}_0C_1}
\dr[0^6](0,1^4)(1,0){6}{13}{E_1}{V'\to\P^3\supset {}_8C_9}
\dr[0,1^5](-1^4,0)(1,-1){4}{22}{E_1}{V'\to\Q^3\supset {}_9C_{11}}
\dr[0^2,1^4](-1^2,0^3)(1,-1){2}{43}{D_1}{S_5\mbox{-fibration}\bigm/\P^1}
\end{array}$$
}

\Sup(df5s){
    In the above list, the $3$-fold $V$ has the following properties in each case.
\Titem
    $(-K_V)^3=22$; $V$ has only one flopping curve $s_0$ as a section; the $(-F_V)$-flop $V'$ has a $\P^1$-bundle structure over $\P^2$; $V'$ contains the strict transform of $s_0$ as a subsection; the anti-canonical model $\ol{V}\subset\P^{13}$ is an intersection of quadrics with one ODP.
\Titem
    $(-K_V)^3=20$; $V$ is the Fano $3$-fold which is the blowing-up of $B_5$ along a quintic elliptic curve, where $B_5$ is a Fano $3$-fold of index $2$ defined by the Pfaffian of the $4\times4$ diagonal minors of a $5\times5$ skew-symmetric matrix whose elements are members of $H^0(\P^5,\O(1))$; $V$ is No.$14$ on Table 3 in \cite{MMT}.
\Titem
    $(-K_V)^3=18$; $V$ has only one flopping curve $s_0$ as a section; the $(-F_V)$-flop $V'$ is the blowing-up of $B_4$ along a rational cubic curve; $\ol{V}\subset\P^{11}$ is an intersection of quadrics with one ODP.
\Titem
    $(-K_V)^3=16$; $V$ has only one flopping curve as a section; the $(-F_V)$-flop $V'$ is the blowing-up of a quadric $3$-fold $\Q^3$ along a smooth curve of degree $7$ and of genus $3$; $\ol{V}\subset\P^{10}$ is an intersection of quadrics with one ODP.
\Titem
    $(-K_V)^3=14$; general $V$ has exactly two flopping curves as sections; the $(-F_V)$-flop $V'$ is the blowing-up of $\P^3$ along a smooth curve of degree $7$ and of genus $4$; $\ol{V}\subset\P^{9}$ is an intersection of quadrics with two ODP's.
\Titem
    $(-K_V)^3=12$; general $V$ has exactly three flopping curves as sections; the $(-F_V)$-flop $V'$ has a conic bundle structure over $\P^2$ with the discriminant locus of degree $5$; $\ol{V}\subset\P^{8}$ is an intersection of quadrics with three ODP's.
\Titem
    $(-K_V)^3=10$; general $V$ has exactly five flopping curves as sections; the $(-F_V)$-flop $V'$ has a structure of del Pezzo fibration of degree $5$ over $\P^1$ again; $\ol{V}=Gr(5,2)\cap H_1\cap H_2\cap Q\subset \P^9$, where $H_1$ and $H_2$ are general hyperplanes, and $Q$ is a quadric hypersurface of rank $4$; the singularities of $\ol{V}$ come from the singular locus $\Sigma\cong\P^5$ of $Q$, and are the intersection $\Sigma\cap Gr(5,2)\cap H_1\cap H_2$, five ODP's.
\Titem
    $(-K_V)^3=8$; general $V$ has exactly eight flopping curves as sections; the $(-F_V)$-flop $V'$ is the blowing-up of $V_{12}$ along a line, where $V_{12}$ is a Fano $3$-fold of index $1$ with genus $7$; $\ol{V}=Q_1\cap Q_2\cap Q_3\subset\P^6$, where $Q_i$ are smooth quadric hypersurfaces in the special position; the singularities of $\ol{V}$ are eight ODP's corresponding to the minimal sections of $V\to\P^1$.
\Titem
    $(-K_V)^3=6$; general $V$ has exactly thirteen flopping curves as sections; the $(-F_V)$-flop $V'$ is the blowing-up of $\P^3$ along a smooth curve of degree $9$ and of genus $8$; $\ol{V}$ is a complete intersection of a smooth quadric and a smooth cubic $4$-fold with thirteen ODP's as its singularities.
\Titem
    $(-K_V)^3=4$; general $V$ has exactly twenty-two flopping curves $s_1, \dots, s_{21}, C_1$; $C_1$ is a bisection and the others are sections of the fibration; the $(-F_V)$-flop $V'$ is the blowing-up of a quadric $3$-fold $\Q^3$ along a smooth curve of degree $11$ and of genus $9$; $\ol{V}$ is a quartic $3$-fold with $22$ ODP's.
\Titem
    $(-K_V)^3=2$; general $V$ has exactly $125$ flopping curves as sections; the $(-F_V)$-flop $V'$ has a structure of del Pezzo fibration of degree $5$ over $\P^1$ again; $\ol{V}$ is the double covering of $\P^3$ branched along a sextic surface with $43$ ODP's.
}


    Finally, we mention that different classes in the above theorems are not deformation equivalent.

\Thm(deform){Every weak Fano $3$-fold in the list of Theorems $(\ref{thm:p2})$--$(\ref{thm:df5})$ does not deform to one of the others in the list.
}

{\it Proof}.~
    Since the degree $(-K_V)^3$ is a deformation invariant, it has only to check the inequivalence for the case of the same degree.
    We calculate the deformation invariant $d(V)$ introduced in \cite{MM} (7.33):
$$d(V) = \det \left|\matrix{(-K_V\cdot H_V^2) & (-K_V\cdot H_V\cdot F_V) \cr (-K_V\cdot H_V\cdot F_V) & (-K_V\cdot F_V^2)}\right|.$$
    Since $(-K_V\cdot F_V^2) = 0$, we have $d(V)=-(-K_V\cdot H_V\cdot F_V)^2=-16$ if $V$ has an extremal ray of type $D_2$, and $d(V)=-d^2$ if of type $D_1$ of degree $d$.
    We recall that the flopping curve $C_0$ in the central fiber $V_0$ must deform out to the general fiber $V_t$ in a family ${\cal V}\to\Delta$, because, for the embedding $\iota : \P^1\to C_0\subset V_0$, the dimension of the deformation space ${\rm Hom}(\P^1,{\cal V})$ at $\iota$ is
\begin{eqnarray*}
\dim_{[\iota]} {\rm Hom}(\P^1,{\cal V}) \hspace{-15pt}
&&\geq \dim H^0(C_0, \iota^\ast T_{\cal V}) - \dim H^1(C_0, \iota^\ast T_{\cal V}) \\
&&= -(\iota_\ast(\P^1)\cdot K_{V_0}) + (1-g(\P^1))\dim {\cal V} \\
&&= 4 \\
&&> \dim Aut(\P^1).
\end{eqnarray*}

    Thus, the rest is to check the inequvalence for each of the following three sets of cases:
    (i) $\{$ (\ref{thm:quad}.7), (\ref{thm:quad}.8), (\ref{thm:df4}.2) $\}$, 
    (ii) $\{$ (\ref{thm:quad}.9), (\ref{thm:df4}.7) $\}$, and 
    (iii) $\{$ (\ref{thm:df4}.9), (\ref{thm:df4}.10) $\}$. 
    For the set (ii), the inequivalence between (\ref{thm:quad}.9) and (\ref{thm:df4}.7) is obtained by a similar argument as (7.32) in \cite{MM}:
    $V$ in the case (\ref{thm:quad}.9) has a divisor $H_V$ with $(H_V)^3 = 7 \not\equiv 0 ({\rm mod} 4)$, while, for $V$ in the case (\ref{thm:df4}.7), we have $(aH_V+bF_V)^3 = 4a^2(2a+3b) \equiv 0 ({\rm mod} 4)$ for any $a, b\in \Z$.
    For the set (i), we can check the inequivalence as follows.
    We can see that (\ref{thm:df4}.2) is inequivalent to the others, because $(aH_V+bF_V)^3 = 4a^2(4a+3b) \equiv 0 ({\rm mod} 4)$ for any $a, b\in \Z$ in the case (\ref{thm:df4}.2), while $(H_V)^3 \not\equiv 0 ({\rm mod} 4)$ in the other cases, as similar to (ii).
    We will show that $V$ deforms only to the same embedding type in each case of (\ref{thm:quad}.7) and (\ref{thm:quad}.8).
    Let ${\cal V}\to\Delta$ be a deformation of $V=V_0$, and $V_n$ the fiber over the $n$-th infinitesimal neighbourhood of the center $0\in\Delta$.
    Since $-K_V$ is nef and big, it extends to $V_n$, and the ambient projective space bundle $X=\P(\varphi_*(-K_V))\to\P^1$ also extends to $V_n$.
    Hence ${\cal V}\to\Delta$ is embedded in $X\times\P^1\to\P^1$.
    The central fiber $V=V_0$ is a divisor linearly equivalent to $|2H+kF|$ in $X$.
    In our case, $2H+kF$ is ample and extends to $X_n$ as $2H_{X_n}+kF_{X_n}$.
    Consequently, the type (\ref{thm:quad}.7) (resp. (\ref{thm:quad}.8)) deforms to the type (\ref{thm:quad}.7) (resp. (\ref{thm:quad}.8)).
    Thus we check the inequivalence for the set (i).
    For the set (iii), the inequivalence derives as similar to (i). \qed

\section{Preliminaries} \label{sec:prel}

\newpar\Label{subsec:setup}
    We begin with results on del Pezzo fibrations, which says that any del Pezzo fibration $\varphi:V\to\P^1$ is canonically in the diagram:
$$\begin{array}{ccl}
  V & \Longarrow{\psi}{} & X \vspace{5pt}\\
  \Downarrow{\varphi} & \Swarrow{\pi} \vspace{5pt}\\
  \P^1
\end{array}$$
where $X$ is a (weighted) projective space bundle over $\P^1$.
    The del Pezzo fibrations arising from the contraction of extremal ray are of degree $d$, $1 \leq d \leq 6$, or $8$, $9$.
    In this paper, we treat the cases of degree $d\ne6$.
    It is trivial in the case of degree $d=9$; this is the case when $V$ is a $\P^2$-bundle.
    If $d=8$, $V$ is a quadric bundle.
    From the general theory of quadric bundles(e.g. Beauville\cite{beau}), we can represent any quadric bundle as a hypersurface in a $\P^3$-bundle.
    More precisely,

\Thm(beau){{\rm (Beauville\cite{beau}, p.321)}~
    Any quadric surface bundle $V$ over $\P^1$ is embedded in $X=\P[0,a_1,a_2,a_3]\to\P^1$ as a member of $|2H+kF|$ for some integer $k$.
}

    Similarly, from the results about del Pezzo fibrations over curves (Fujita\cite{fujita}), we have

\Thm(dp){
    Let $V$ be a del Pezzo fibration over $\P^1$.
\\ \noindent
{\rm (i)}~
    If $V$ is of degree $1$, then $V$ is embedded in the weighted projective space bundle $X=\P[0,a;b;c]\to\P^1$ as a member of $|6H+kF|$ for some integer $k$, where $\P[0,a;b;c]$ denotes $\proj(\S)$ for the graded $\O_{\p^1}$-algebra $\S = \bigoplus_{d \geq 0} \S_d$ such that $\S_d = \bigoplus_{i+2j+3k=d} S^i(\O\op\O(a))\otimes S^j(\O(b))\otimes S^k(\O(c))$ for some integers $a$, $b$ and $c$.
\\ \noindent
{\rm (ii)}~
    If $V$ is of degree $2$, then $V$ is a double covering of $X=\P[0,a_1,a_2]\to\P^1$ branched along the member $B$ of $|4H+2kF|$ for some integer $k$.
\\ \noindent
{\rm (iii)}~
    If $V$ is of degree $3$, then $V$ is embedded in $X=\P[0,a_1,a_2,a_3]\to\P^1$ as the member of $|3H+kF|$ for some integer $k$.
\\ \noindent
{\rm (iv)}~
    If $V$ is of degree $4$, then $V$ is embedded in $X=\P[0,a_1,a_2,a_3,a_4]\to\P^1$ as the complete intersection of two members of $|2H+k_1F|$ and $|2H+k_2F|$ for some integers $k_1$, $k_2$.
\\ \noindent
{\rm (v)}~
    If $V$ is of degree $5$, then $V$ is embedded in $X=\P[0,a_1,a_2,a_3,a_4,a_5]\to\P^1$ as the subvariety defined by the Pfaffian of the $4\times4$ diagonal minors of a $5\times5$ skew-symmetric matrix $M = (m_{ij})$ such that its entries $m_{ij}$ are members of $|H+w_{ij}F|$ for some integers $w_{ij}$.
}

    These are also obtained from Reid's result (\cite{delPezzo}, p.697) that anti-canonical linear system of any (possibly nonnormal) del Pezzo surface of degree $d\geq 3$ is very ample and the image is ideal-theoretically an intersection of quadrics if $d\geq 4$, and that any one of degree $d\leq 4$ is a weighted hypersurface or complete intersection $S_6 \subset \P(1^2,2,3)$, $S_4 \subset \P(1^3,2)$, $S_3 \subset \P^3$, or $S_{2,2} \subset \P^4$.
    In the case of degree $d=2$, we can also consider $V$ as the hypersurface of relative degree $4$ in the weighted projective space bundle with weights $(1^3,2)$.

\newpar
    From this point on we can fix the situation described in the above theorems with the following notation.

    Let $X$ be the $\P^n$-bundle $\P(\E)=\P[0,a_1,\dots,a_n]$ over $\P^1$ for a locally free sheaf of rank $n+1$ $\E=\bigoplus^n_{i=0}\O(a_i)$, $0=a_0\leq a_1\leq\cdots\leq a_n$, on $\P^1$.
    Let $H$ and $F$ be the tautological divisor and a fiber of $\pi:X=\P[0,a_1,\dots,a_n] \to\P^1$, respectively.
    Then, the set of numerical equivalence classes of divisors and 1-cycles on $X$ are
$$
  N^1(X)=\R[H]\op\R[F],\quad
  N_1(X)=\R[l]\op\R[s_0],
$$
respectively.
    Here $l$ is a line in a fiber and $s_0$ is the minimal section associated to the exact sequence $\E\to\O\to 0$; their classes form the dual basis of $[H]$ and $[F]$ with respect to intersection number.
    The closure $\ol{NE}(X)$ of the convex cone generated by numerical classes of effective 1-cycles of $X$ (the Mori cone of $X$) is
$$
  \ol{NE}(X)=\R_+[l]+\R_+[s_0] \subset N_1(X).
$$
    Here $\R_+$ denotes the set of non-negative real numbers.

    The canonical divisor of the $\P^n$-bundle $X=\P[0,a_1,\dots,a_n]\to\P^1$ is
\Eq(can){
K_X\sim -(n+1)H+(-2+\sum_{i=0}^n a_i)F.
}\\
\Titem{
    When $n=2$, the anti-canonical divisor is $-K_X\sim 3H+(2-a_1-a_2)F$.
    If $2-a_1-a_2\geq0$, then $-K_V$ is nef.
    Conversely assume that $-K_V$ is nef.
    Since $(-K_X\cdot s_0)=2-a_1-a_2$ for the section $s_0$ of $X$ associated to the exact sequence $\O\op\O(a_1)\op\O(a_2)\to\O\to0$ of sheaves on $\P^1$, the pair $(a_1,a_2)$ is equal to one of $(0,0)$, $(0,1)$, $(0,2)$, and $(1,1)$.
    In the case $(a_1,a_2)=(0,2)$, there exists a family of $K$-trivial curves on $X$, which shows Theorem(\ref{thm:p2}). \label{subsubsec:p2}}

    Let
$$\begin{array}{ccl}
  V & \Longarrow{\psi}{} & X \vspace{5pt}\\
  \Downarrow{\varphi} & \Swarrow{\pi} \vspace{5pt}\\
  \P^1
\end{array}$$
be the diagram stated in (\ref{subsec:setup}): here $\psi$ is an isomorphism in the case of degree 9, an embedding in the case of degree 3, 4, 5 or 8, and a double covering in the case of degree 2, respectively.
    The case of degree 1 is slightly different, and is stated in Section \ref{sec:df1} precisely.

    In addition, we assume that the Picard number $\rho(V)$ of $V$ equals $2$; this assumption is natural because the extremal ray contraction of type $D$ has the relative Picard number $1$.
    Then, the restriction map $\psi^\ast:N^1(X)\to N^1(V)$ is an isomorphism and so is $\psi_\ast:N_1(V)\to N_1(X)$.
    Under the isomorphism $\psi_\ast$, $N_1(V)$ is identified with $N_1(X)$; the Mori cone $\ol{NE}(V)$ of $V$ maps into $\ol{NE}(X)$ but not surjectively in general; the edges of $\ol{NE}(V)$ are generated by $[l]$ and $[C_V]$ for some rational curve $C_V\subset V$ with $(\psi^\ast H\cdot C_V)†0$.
$$\begin{array}{rcccl}
  & N_1(V) & \Longarrow{\sim}{\psi_\ast} & N_1(X) & \\
  \vspace{0.5\jot} & \cup & & \cup & \\ \vspace{1\jot}
  \R_+[l] + \R_+[C_V] = \hspace{-.8em} & \ol{NE}(V) & \Longarrow{}{} & \ol{NE}(X) & \hspace{-.8em}=\R_+[l]+\R_+[s_0]
\end{array}$$
    We denote $H_V = \psi^\ast H$ and $F_V = \psi^\ast F$.
    For the sake of convenience, we set $\mu(Z) = (H_V\cdot Z)_X/(F_V\cdot Z)_X$ ($= (H\cdot Z)/(F\cdot Z)$) for each $1$-cycle $Z$ on $V\subset X$, and call it a {\it slope} of $Z$.
    Moreover we set $\mu(R)=\mu(Z)$ for any ray $R=\R_+[Z]\subset N_1(V)$, and call it a {\it slope} of $R$.

\newpar
    If $V$ has only finite $K$-trivial curves, there is a flop $\chi:V\cdots\to V'$ corresponding to the edge $R_0\subset\ol{NE}(V)$ generated by the $K$-trivial curves.
    We have the canonical isomorphism $\chi^\ast:\Pic(V') \to \Pic(V)$, and we can identify $N^1(V')$ with $N^1(V)$ by $\chi^\ast:N^1(V') \to N^1(V)$; put $H_{V'} = \chi^\ast{H_V}$, $F_{V'} = \chi^\ast{F_V}$.
    As stated above, we can consider the slope $\mu(Z') = (H_{V'}\cdot Z')/(F_{V'}\cdot Z')$ of 1-cycle $Z'$ on $V'$.
    If $Z'$ is the strict transform of $Z\subset V$, we have the following expression of the value $\mu(Z')$.

    Let $C_1,\dots,C_k\subset V$ be all the curves having the same slope $\mu(R_0)$, and put $a = (H_V\cdot Z)$, $b = (F_V\cdot Z)$, $\alpha_i = (H_V\cdot C_i)$, $\beta_i = (F_V\cdot C_i)$, $n_i = \sharp\{Z\cap C_i\}$ for each $i$.
    Then
\Eq(slope){
\mu(Z') = \dfrac{a + \sum_{i=1}^k n_i \alpha_i}{b + \sum_{i=1}^k n_i \beta_i}.}\\
    This formula is useful in fixing an extremal ray of $\ol{NE}(V')$.

\newpar \Label{subsec:ineq}
    Consider the sub-$\P^{n-2}$-bundle $T=\P[0,a_1,\dots,a_{n-2}]$ in $X$ associated to the exact sequence
$$
  \E = \bigoplus^n_{i=0} \O(a_i) \longrightarrow
  \bigoplus^{n-2}_{i=0} \O(a_i) \longrightarrow 0
$$
of sheaves, and denote the pull-back $\psi^\ast T$ in $V$ by $D$.
    The first inequality
\Eq(nef){
  (D\cdot -K_V)\geq 0}\\
obtained from the nefness of $-K_V$
    Since $H_V$ is nef, the second inequality
\Eq(eff){
  (D\cdot H_V)\geq 0}\\
holds.
    The third inequality
\Eq(big){
  (-K_V)^3>0}\\
follows from the bigness of $-K_V$.

    In Sections \ref{sec:quad}--\ref{sec:df1}, we restrict $\E$ and $V$ through these three inequalities (\ref{eqn:nef})--(\ref{eqn:big}), and then we check the geometrical realizability in order to classify the weak Fano 3-folds with an extremal ray of type $D$.
    Sections \ref{sec:quad}--\ref{sec:df1} has the same frame: first we state the situation in (?.1), restrict $\E$ and $V$ in (?.2)--(?.3), and list the possibilities in the table (?.4); next, in (?.5), we exclude some cases in the table (?.4), and finally, in (?.6), we construct the geometrical objects for each of the other cases.
    In (?.6), we devote ourselves to verifying that $V$ is smooth, that the Picard number $\rho(V)$ of $V$ is two, and that the number of $K$-trivial curves on $V$ is finite.

\newpar
    The following lemma assures that $D$ in (\ref{subsec:ineq}) can be used as a substitute for the edge curve $C_V$ of the Mori cone $\ol{NE}(V)$.
\Titem \setcounter{thm}{\value{subsection}} \Label{thm:slope} {\bf Lemma.}~{\it
    Let $X = \P[0,a^k,a_{k+1},\dots,a_n]$ be a $\P^n$-bundle over $\P^1$ for a sequence of integers $0=a_0 \leq a_1 = \dots = a_k (= a) < a_{k+1} \leq \dots \leq a_n$.
    Assume that $C$ is a rational curve in $X$ with slope $\mu(C)\leq a$.
    Then, the curve $C$ is on a ruled surface $S_\lambda = \P[0,a]$ determined by $\lambda \in \P^{k-1}$.}

{\it Proof}.~
    Let $D$ be the normalization of $C$ and let $\mu : D \to X$ and $\sigma = \pi \circ \mu : D \to \P^1$ be the natural morphisms.
    Denote $(F \cdot C)_X$ and $(H \cdot C)_X$ respectively by $d$ and $b$.
    Then, we have ${\rm deg }\sigma = d$ and $\sigma^\ast \E \to \O_C(b) \to 0$.
    Since $\sigma^\ast \O(1) = \O_C(d)$, we have an exact sequence
\Eq(seq){
  \bigoplus^n_{i=0} \O_C(da_i) = \O_C \op \O_C(da)^{\op k} \op \bigoplus^n_{i=k+1} \O_C(da_i) \longrightarrow \O_C(b) \longrightarrow 0.}

    If $0\leq b<da$, then the sequence (\ref{eqn:seq}) factors through $\O_C$ and the curve $C$ is a multiple of the unique minimal section $s_0 = \P[0] \subset X$ associated to $\E \to \O \to 0$.
    In this case, the assertion is trivial.

    We now assume that $0 \leq b = da$.
    Then the sequence (\ref{eqn:seq}) factors through $\O_C \op \O_C(b)^{\op k}$.
    By choosing a basis of $\O_C(b)^{\op k}$ such that $\O_C \op \O_C(b)^{\op k} \to \O_C(b) \to 0$ factors through $\O_C \op \O_C(b)$, the morphism $\mu$ factors through $S_\lambda = \P[0,a]$, where $\lambda \in \P^{k-1}$ indicates the choice of the basis.
\qed

\newpar
    We end off this section by recalling the Lefschetz theorem on the Picard group with its corollaries.

\def\Leff{\rm Leff}

\Thm(PicIso){{\rm (\cite{H} p.178)}~ 
    Let $X$ be a complete non-singular variety, and let $Y$ be a closed subscheme.
    Assume
\\
{\rm (i)}~ $\Leff(X,Y)$,
\\
{\rm (ii)}~ $Y$ meets every effective divisor on $X$, and
\\
{\rm (iii)}~ $H^i(Y,{\cal I}^n/{\cal I}^{n+1}) = 0$ for $i=1$, $2$, and all $n \geq 1$, where $\cal I$ is the sheaf of ideals of $Y$.
\\
    Then the natural map $\Pic(X) \to \Pic(Y)$ is an isomorphism.}

    Here $\Leff(X,Y)$ is the effective Lefschetz condition, this holds as in the following case:

\Thm(Leff){{\rm (\cite{H} p.172)}~ 
    Let $X$ be a non-singular subvariety of $\P^N_k$, and let $Y$ be a closed subscheme of codimension $r$, which is a complete intersection(i.e., $Y=X \cap H_1 \cap \dots \cap H_r$ where the $H_i$ are hypersurfaces in $\P^N_k$).
    Assume that $\dim Y \geq 2$.
    Then we have $\Leff(X,Y)$.}

\newpar
    Now we derive some corollaries of (\ref{thm:PicIso}), used in the following sections.

\Cor(ampleLeff){
    Let $X$ be a non-singular projective variety, and let $Y$ be an ample effective divisors on $X$.
    Assume that $\dim Y \geq 3$.
    Then the natural map $\Pic(X) \to \Pic(Y)$ is an isomorphism.}

\noindent{\it Proof}.~ This is obvious from (\ref{thm:Leff}) and (\ref{thm:PicIso}).
\qed

\Cor(blowup){
    Let $X'$ be a non-singular projective variety, let $Y'=\bigcap_i D'_i$ be a complete intersection of $n$ ample effective divisors $D'_i \subset X'$, and let $Z'$ be a non-singular subvariety of $X'$.
    Let $f : X \to X'$ be the blowing-up of $X'$ along $Z'$, and let $Y \subset X$ be the strict transform of $Y'$.
    Assume that $Z' \cap Y'$ is non-singular and irreducible, and that $\dim Y \geq 3$.
    Then the natural map $\Pic(X) \to \Pic(Y)$ is an isomorphism.}

\noindent{\it Proof}.~
    We have $\Pic(X') \cong \Pic(Y')$ by (\ref{thm:ampleLeff}), and $\Pic(X) \cong f^\ast \bigl(\Pic(X')\bigr) \op \Z$.
    From the condition on $Z'\cap Y'$, we have $\Pic(Y) \cong f^\ast \bigl(\Pic(Y')\bigr) \op \Z$.
\qed

    Applying the above corollary to the morphism $f=\phi_{|H|} : X=\P[0^n,1] \to \P^{n+1}$ associated to the tautological divisor $H$ of $X \to \P^1$, we have the following two corollaries.

\Cor(0n1){
    Let $X$ be the $\P^n$-bundle $\P[0^n,1] \to \P^1$, and let $Y$ be a complete intersection of general members of linear systems $|m_iH|$ on $X$ for positive integers $m_i$.
    Assume that $\dim Y \geq 3$.
    Then $\Pic(Y) \cong \Z^{\op2}$.} \qed

\Cor(0n11){
    Let $X$ be the $\P^n$-bundle $\P[0^{n-1},1^2] \to \P^1$, and let $Y$ be a complete intersection of general members of linear systems $|m_iH|$ on $X$ for positive integers $m_i$.
    Assume that $\dim Y \geq 3$.
    Then $\Pic(Y) \cong \Z^{\op2}$.}

\noindent{\it Proof}.~
    Let $W=\P[0^{n+1},1]$ be a $\P^{n+1}$-bundle containing $X$ as an ample effective divisor.
    The morphism $\phi_{|H|}$ is the restriction of $\phi_{|H_W|} : W \to \P^{n+2}$ associated to the tautological divisor $H_W$ of $W \to \P^1$.
    There is a complete intersection $V'$ of hypersurfaces in $\P^{n+2}$ such that $Y=V \cap X$, where $V$ is a pulled-back of $V'$.
    (\ref{thm:ampleLeff}) derives $\Pic(V) \cong \Pic(Y)$ from the ampleness of $X\subset W$, and (\ref{thm:blowup}) deduces $\Pic(V) \cong \Z^{\op2}$ from the generality of $V'$.
\qed

\Cor(01n){
    Let $X$ be a $\P^n$-bundle $\P[0,1^n] \to \P^1$, and let $Y$ be a general member of linear system $|mH|$ on $X$ for a positive integer $m$.
    Assume that $\dim Y \geq 3$.
    Then $\Pic(Y) \cong \Z^{\op2}$.}

\noindent{\it Proof}.~
    Modifying the argument of proof of (\ref{thm:Leff}) (described in \cite{H}) slightly, we obtain $\Leff(X,Y)$ in this case too.
    The conditions (ii) and (iii) of (\ref{thm:PicIso}) are checked immediately.
    Thus $\Pic(Y) \cong \Z^{\op2}$ since $\Pic(X) \cong \Z^{\op2}$. \qed

\def\On{\O_{\p^{n-l}}}

\Cor(0k1n){
    Let $X=\P[0^l,1^{n-l+1}]$ be a $\P^n$-bundle over $\P^1$ for $n>l>0$, and let $H$, $F$ be the tautological divisor and a general fiber of $X\to\P^1$, respectively.
    Let $Y$ be a general members of the linear systems $|dH-F|$ on $X$ for $d>0$.
    Assume that $\dim Y \geq 3$.
    Then $\Pic(Y)\cong\Z^{\op 2}$.}

\noindent{\it Proof}.~
    Consider the graded ring $R = \C[x_0, x_1, \dots, x_{l-1}, y_l, y_{l+1}, \dots, y_n, t_0, t_1]$ with bidegree $\deg x_i = (1,0)$, $\deg y_i = (1,-1)$ and $\deg t_i = (0,1)$ as the bihomogeneous coordinate ring of $X$.
    Then, each member $Y$ of $|dH-F|$ is defined by a bihomogeneous polynomial $F\in R$ of bidegree $\deg F = (d,-1)$; the polynomial $F$ is a linear combination of monomials $x^Iy^Jt^K$ for $|I|+|J|=d$ and $|K|=|J|-1$, where, for the set $I=\{i_0, \dots, i_{l-1}\}$ of non-negative integers, $x^I$ and $|I|$ denote the monomial ${x_0}^{i_0}\cdots{x_{l-1}}^{i_{l-1}}$ and the sum $i_0+\cdots+i_{l-1}$, etc.
    It is obvious that $Y$ is smooth for general $F$.
    Let $L=\P[0^1]$ be a sub-$\P^{l-1}$-bundle of $X$ associated to the natural surjection $\O^{\op l}\op\O(1)^{\op(n-l+1)}\to\O^{\op l}$.
    Then, $L$ is the common zero locus of $y_l, y_{l+1}, \dots, y_n$, and $Y$ for any $F$ contains $L$.
    Consider the $\P^{l+1}$-bundle $X' = \P(\On^{\op l}\op\On(1)^{\op2}) \to \P^{n-l}$ having the bihomogeneous coordinate ring $R' = \C[x_0', x_1' \dots, x_{l-1}', z_l', z_{l+1}', s_0', s_1', \dots, s_{n-l}']$ with bidegree $\deg x_i' = (1,0)$, $\deg z_i' = (1,-1)$ and $\deg s_i' = (0,1)$.
    Let $H'$ and $F'$ be the tautological divisor and a fiber of $X'\to\P^{n-l}$.
    Let $L'=\P(\On^{\op l})$ be the sub-$\P^{l+1}$-bundle of $X'$ associated to the natural surjection $\On^{\op l}\op\On(1)^{\op2}\to\On^{\op l}$.
    Then, $L$ is the common zero locus of $z_l'$ and $z_{l+1}'$.
    The ring isomorphism $R\to R'$ determined by $(x_i', z_{l+i}', s_i')\to(x_i, t_i, y_{l+i})$ induces a birational map $f: X\cdots\to X'$.
    The birational map $f: X\cdots\to X'$ is isomorphic in codimension 1, and is decomposed into the blowing-up and blowing-down $X\gets\cdots X''\cdots\to X'$ with centers $L\gets E''\to L'$, where $E'' \cong \P^{n-l}\times\P^{l-1}\times\P^1$ is the exceptional divisor for the blowing-up and the blowing-down.
    Under $f$, the members of $|dH-F|$ map to the members of $|(d-1)H'+F'|$, because the monomials $x^Iy^Jt^K$, $|I|+|J|=d$ and $|K|=|J|-1$, are the image of ${x'}^I{z'}^K{s'}^J$, $|I|+|K|=d-1$ and $|J|=|K|+1$, by the ring isomorphism.
    Moreover, the image of the smooth member $Y\in|dH-F|$ is a smooth member $Y'\in|(d-1)H'+F'|$, and the restriction $f|_Y$ of $f: X\cdots\to X'$ is isomorphic in codimension 1, hence $\Pic(Y)\cong\Pic(Y')$.
    Since $(d-1)H'+F'$ is ample in $X'$, (\ref{thm:ampleLeff}) shows $\Pic(Y')\cong\Pic(X')\cong\Z^{\op2}$. \qed

\section{Quadric bundles} \label{sec:quad}

\newpar \Label{subsec:setupq}
    In this section, we will derive Theorem(\ref{thm:quad}) and Supplement(\ref{thm:quads}).
    Let $V$ be a quadric bundle in $X=\P[0,a_1,a_2,a_3]$, $0\leq a_1\leq a_2\leq a_3$, such that $V\in |2H+kF|$ (see (\ref{subsec:setup})):
$$\begin{array}{ccl}
  V & \Longarrow{\psi}{} & X = \P[0,a_1,a_2,a_3] \vspace{5pt}\\
  \Downarrow{\varphi} & \Swarrow{\pi} & \qquad V \in |2H + kF| \vspace{5pt}\\
  \P^1
\end{array}$$
    Assume that $V$ is a weak Fano $3$-fold with only finite $K$-trivial curves and with Picard number $\rho(V)=2$.

\newpar
    Under the situation (\ref{subsec:setupq}), the anti-canonical divisor of $V$ is
$$
  -K_V\sim 2H_V-(\sum_{i=1}^3a_i+k-2)F_V
$$
from (\ref{eqn:can}), the intersection numbers are
$$
  (H_V^3)=2\sum_{i=1}^3a_i+k, \quad (H_V^2\cdot F_V)=2, \quad
  {\rm and}\quad (l\cdot H_V)=(s_0\cdot F_V)=1,
$$
and the curve $D$ described in (\ref{subsec:ineq}) is
$$
  D \equiv (H_V-a_2F_V)\cdot(H_V-a_3F_V) \equiv H_V^2-(a_2+a_3)H_V\cdot F_V,
$$
hence
$$
  D \equiv (2a_1+k)[l]+2[s_0].
$$
    The three inequalities (\ref{eqn:nef}), (\ref{eqn:eff}) and (\ref{eqn:big}) described in the previous section imply
\Eq(nefq){
  2(a_1+2-a_2-a_3) \geq 0,}
\Eq(effq){
  2a_1+k \geq 0,}\\
and
\Eq(bigq){
  8(6-\sum_{i=1}^3a_i-2k) > 0,}\\
respectively.

\newpar
    It follows $0\leq a_1\leq a_2\leq a_3\leq 2$ from (\ref{eqn:nefq}), and hence $0\leq \sum_{i=1}^3a_i\leq 6$.
    Furthermore, we have the following seven possibilities of triples $(a_1,a_2,a_3)$:
\Eq(000){
  (a_1, a_2, a_3) = (0, 0, 0),}
\Eq(001){
  (a_1, a_2, a_3) = (0, 0, 1),}
\Eq(002){
  (a_1, a_2, a_3) = (0, 0, 2),}
\Eq(011){
  (a_1, a_2, a_3) = (0, 1, 1),}
\Eq(111){
  (a_1, a_2, a_3) = (1, 1, 1),}
\Eq(112){
  \phantom{\mbox{\quad and}}(a_1, a_2, a_3) = (1, 1, 2),\mbox{\quad and}}
\Eq(222){
  (a_1, a_2, a_3) = (2, 2, 2).}\\
    In each case (\ref{eqn:002}), (\ref{eqn:112}), or (\ref{eqn:222}), the equality holds in (\ref{eqn:nefq}), hence $(-K_V\cdot D)=0$ and $D$ is a $K$-trivial curve.
    On the other hand, the ruled surface $S$, associated to the exact sequence $\E\to\O\op\O(a_1)\to 0$ of sheaves, has a family of dimension $\geq1$ in each case.
    This contradicts the assumption that $V$ has only finite $K$-trivial curves.
    Thus these cases are excluded.

    Based on the description in (\ref{eqn:effq}) and (\ref{eqn:bigq}) we can restrict the value of $k$ as follows:
$$\begin{array}{ll}
  k = 0, 1, \mbox{ or } 2, & \mbox{ in case (\ref{eqn:000});} \\
  k = 0, 1, \mbox{ or } 2, & \mbox{ in case (\ref{eqn:001});} \\
  k = 0, \mbox{ or } 1, & \mbox{ in case (\ref{eqn:011});\quad and} \\
  k = -2, -1, 0, \mbox{ or } 1, \quad\quad & \mbox{ in case (\ref{eqn:111}).}
\end{array}$$

\newpar
    Thus we obtain the possibilities of numbers $a_1,a_2,a_3,k$, and $(-K_V^3)$ as in Table \ref{tbl:quad}.
\begin{table}[htbp]
\caption{}\label{tbl:quad}
$$\begin{array}{c@{\qquad}cccc@{\qquad}l@{\qquad}l@{\qquad}rc}
{\rm Nos.}& a_1& a_2& a_3&  k &\qquad V&\quad -K_V &  D \qquad & (-K_V^3) \\
  1 &  1 &  1 &  1 & -2 & 2H-2F & 2H_V+ F_V &   2s_0&  56 \\
  2 &  0 &  0 &  0 &  0 & 2H    & 2H_V+2F_V &   2s_0&  48 \\
  3 &  0 &  0 &  1 &  0 & 2H    & 2H_V+ F_V &   2s_0&  40 \\
  4 &  1 &  1 &  1 & -1 & 2H- F & 2H_V      & l+2s_0&  40 \\
  5 &  0 &  0 &  0 &  1 & 2H+ F & 2H_V+ F_V & l+2s_0&  32 \\
  6 &  0 &  1 &  1 &  0 & 2H    & 2H_V      &   2s_0&  32 \\
  7 &  0 &  0 &  1 &  1 & 2H+ F & 2H_V      & l+2s_0&  24 \\
  8 &  1 &  1 &  1 &  0 & 2H    & 2H_V- F_V &2l+2s_0&  24 \\
  9 &  0 &  0 &  0 &  2 & 2H+2F & 2H_V      &2l+2s_0&  16 \\
 10 &  0 &  1 &  1 &  1 & 2H+ F & 2H_V- F_V & l+2s_0&  16 \\
 11 &  0 &  0 &  1 &  2 & 2H+2F & 2H_V- F_V &2l+2s_0&   8 \\
 12 &  1 &  1 &  1 &  1 & 2H+ F & 2H_V-2F_V &3l+2s_0&   8
\end{array}$$
\end{table}
    The columns \lq\lq$V$\rq\rq and \lq\lq$-K_V$\rq\rq in Table \ref{tbl:quad} denote a linear equivalence class of $V$ and $-K_V$, respectively, and the column \lq\lq$D$\rq\rq a numerical class of $D$.
    We omit the symbol $[ \quad ]$ for the sake of simplicity.

    In the rest of this section, we consider the realization of each possibility in Table \ref{tbl:quad}.

\newpar
    We now exclude three possibilities, Nos.$1$, $2$, and $11$, in Table \ref{tbl:quad}.
\ex{1}
    This is the case that $V\subset X=\P[0,1^3]$ is a member of $|2H-2F|$.
    But this case is to be excluded, because any member of the linear system $|2H-2F|$ has singularities along the minimal section $s_0=\P[0] \subset \P[0,1^3]$.
\ex{2}
    This is the case that $V\subset X=\P[0^4]$ is a member of $|2H|$.
    In this case, $V$ is a $\P^1\times\P^1$-bundle and its Picard number is three.
    Hence this is an excluded case.
\ex{11}
    In this case, $V\subset X=\P[0^3,1]$ is a member of $|2H+2F|$, and $-K_V\sim 2H_V-F_V$.
    Let $L=\P[0^3]\cong\P^2\times\P^1$ be a sub-$\P^2$-bundle of $X$, and $S=L\cap V$.
    If $V$ is sufficiently general, $S$ is a smooth surface.
    Then, the linear system $|-K_V|\big|_S$ defines a fiber structure on $S$ because of
$$
  \bigl((-K_V)^2 \cdot S\bigr)
    =4(H_V^3) - 8(H_V^2 \cdot F_V) = 0.
$$
    Therefore there exists a $1$-dimensional family of $K$-trivial curves on $S\subset V$, which contradicts the assumption of the finiteness of $K$-trivial curves.

\newpar \Label{subsec:quad}
    Next, we construct the quadric bundle with our conditions in the rest of the cases in Table \ref{tbl:quad}.
    To do this, we have only to confirm that there exists a smooth member $V$ of the linear system $|2H-kF|$, that $V$ has only finite $K$-trivial curves, and that the Picard number $\rho(V)$ of $V$ is two.
    We thus show Theorem(\ref{thm:quad}) and Supplement(\ref{thm:quads}).
\cons{3}
    In this case, we have $X=\P[0^3,1]$.
    The general member $V$ of $|2H|$ is a smooth Fano $3$-fold with $-K_V\sim 2H_V+F_V$ and $\rho(V)=2$ by (\ref{thm:0n1}).
    Let $L=\P[0^3] \cong \P^2\times\P^1$ be a sub-$\P^2$-bundle of $X$.
    The first projection $p : L\to\P^2$ maps the intersection $L\cap V$ onto a smooth conic $C$.
    For any $c\in C$, $s=p^{-1}(c)\cong\P^1$ and $(s\cdot H_V)=0$, hence the class $[s]$ generates an edge of the Mori cone $\ol{NE}(V)$ of $V$.
    The edge is an extremal ray of $E_1$-type, and the corresponding linear system is $|H_V|$.
Since $\dim|H_V|=5$ and $(H_V^3)=2$,
the contraction morphism is a blowing-up of $\Q^3$ along a smooth conic.
    We have $\dim|H_V-F_V|=1$ and the unique member of $|H_V-F_V|$ is the exceptional divisor for $V\to\Q^3$.
    This is (\ref{thm:quad}.1).
\cons{4}
    We have $X=\P[0,1^3]$ and $V\sim 2H-F$.
    The general member $V$ of $|2H-F|$ is smooth, $\rho(V)=2$ by (\ref{thm:0k1n}), and $-K_V\sim 2H_V$.
    Let $s_0=\P[0] \subset X=\P[0,1^3]$ be the minimal section.
    Then, $V$ contains $s_0$ because $(V\cdot s_0)_X=-1$; $\mu(s_0)=0$; and the class $[s_0]$ generates an edge $R={\Bbb R}_+[s_0]\subset\ol{NE}(V)$.
    Since $(R\cdot -K_V)=0$, the edge $R$ is not an extremal ray, and the irreducible curves whose classes are in $R$ should be flopped.
    Any curve in $V$ can be considered a curve in $X$, and $X$ has only one irreducible reduced curve whose class is contained in $R$, which is $s_0$, hence $V$ has also only one $K$-trivial curve $s_0$.
    After flopping along $s_0$, we have again a weak Fano $3$-fold $V'$ with an extremal ray $R'$.
We can see that there is no curve $C'$ with a slope $\mu(C')<1$ as follows: ~
    let $a=(H_V\cdot C), b=(F_V\cdot C)$, and $n_0=\#\{C\cap s_0\}$ for any effective curve $C\subset V$;
    then $(H_V-F_V\cdot C)=a-b\geq n_0$ because $Bs|H_V-F_V|=s_0$;
    hence the strict transform $C'$ of $C$ has a slope $\mu(C')=a/(b+n_0)\geq 1$ from (\ref{eqn:slope}).
    The ray $R'$ is generated by the classes of the strict transforms of $C_\lambda$ which are intersection curves between $V$ and $S_\lambda\subset X$.
    Here $S_\lambda=\P[0,1] \subset X=\P[0,1^3]$ are ruled surfaces parameterized by $\lambda\in\P^2$.
    Thus $R'$ is of $C$-type(indeed $C_1$-type) and $V'$ has a conic bundle structure defined by $|H_{V'}-F_{V'}|$.
    The anti-canonical model $\ol{V}$ is a del Pezzo $3$-fold of degree $5$ having only one ODP, because $\dim|H_V|=5$ and $(H_V)^3=5$.
    This is (\ref{thm:quad}.2).
\cons{5}
    In this case, $X=\P[0^4] \cong \P^3\times\P^1$ and $V\sim 2H+F$ is a smooth Fano $3$-fold with $-K_V\sim 2H_V+F_V$ of $\rho(V)=2$ by (\ref{thm:ampleLeff}).
    The $3$-fold $V$ is defined by an equation $q_0(x)t_0+q_1(x)t_1$, where $x=[x_0:x_1:x_2:x_3]$ and $t=[t_0:t_1]$ are the homogeneous coordinates of $\P^3$ and $\P^1$, respectively, and $q_i(x)$ are quadratic forms in $x$.
    Let $C$ be the set $\bigl\{ \alpha=[\alpha_0:\alpha_1:\alpha_2:\alpha_3] \bigm| q_0(\alpha)=q_1(\alpha)=0 \bigr\}$; the set $C$ is a quartic elliptic curve in $\P^3$ if $V$ is general.
    For any $\alpha\in C$, $s_\alpha = \{\alpha\}\times\P^1 \subset V$ has a slope $\mu(s_\alpha)=0$, hence its class generates an extremal ray $R$ of $NE(V)$.
    Every curve whose class is in $R$ is a section in $V$ and is $s_\alpha$ as above.
    The linear system corresponding to $R$ is $|H_V|$ with $\dim|H_V|=3$ and $(H_V^3)=1$.
    Thus the contraction morphism is the restriction $p|_V : V\to\P^3$ of the first projection $p : X\to\P^3$, i.e., a blowing-up $\P^3$ along a quartic elliptic curve $C$.
    The exceptional divisor for $p|_V$ is the unique member of $|2H_V-F_V|$.
    This is (\ref{thm:quad}.3).
\cons{6}
    In this case, $\E=\O^{\op 2}\op\O(1)^{\op 2}$, $X=\P(\E)=\P[0^2,1^2]$ and $V\sim 2H$ is a smooth weak Fano $3$-fold with $-K_V\sim 2H_V$ of $\rho(V)=2$ by (\ref{thm:0n11}).
    Let $S = \P[0^2] \cong\P^1\times\P^1$ be a ruled  surface in $X$.
    The intersection $V\cap S$ is a pair of sections, i.e., $s_0+s_1$, and these classes generate an edge $R\subset\ol{NE}(V)$ such that $\mu(R)=0$.
    Any curve whose class is in $R$ is nothing but a section in $S$.
    Thus $V$ has only two $K$-trivial curves $s_0, s_1$ and has the $(-F_V)$-flop $V'$, which is again a weak Fano $3$-fold with a unique extremal ray $R'\subset\ol{NE}(V')$.
    The bihomogeneous coordinates analysis implies that $\ol{V}\subset\P^5$ is the intersection of a smooth quadric $\Q^4$ and a singular quadric $\Q_0^4$ of rank $3$.

    Now we will find curves generating the ray $R'$.
    We fix a trivial sheaf ${\cal L}\cong\O$ on $\P^1$ such that an exact sequence $\E\to{\cal L}\to 0$ associates to the section $s_0$.
    Consider a family of ruled surfaces $T_\lambda$ associated to $\E\to{\cal L}\op\O(1)\to 0$ parameterized by $\lambda\in\P^1$.
    Then we have a family of curves $C_\lambda=T_\lambda\cap V$ on $V$.
    Since $\cal L$ associates to $s_0$, each curve $C_\lambda$ is decomposed into $s_0+D_\lambda$.
    Each $D_\lambda$ meets $s_0$ at a point from $(D_\lambda\cdot s_0)_{T_\lambda}=1$.
    Considering the possibility for decomposition of $D_\lambda$, we can see that there are four points on $\P^1$ where $D_\lambda=D_{\lambda,0}+l_0$.
    Here $l_0$ is a line in a fiber $\P^3$.
    Since $D_\lambda\equiv 2[l]+[s_0]$ and $D_{\lambda,0}\equiv [l]+[s_0]$, the strict transforms $D_\lambda'$, $D_{\lambda,0}'$ by the $(-F_V)$-flop are numerically equivalent to $2[l']+2[s_0']$, $[l']+[s_0']$, respectively, therefore $\mu(D_\lambda')=\mu(D_{\lambda,1}')=1$.
    Similarly as in (\ref{subsec:quad}.2), we can see that $\mu(C')\geq 1$ for any curve $C'\subset V'$ by using (\ref{eqn:slope}) for $|H_V-F_V|$ with $Bs|H_V-F_V|=s_0+s_1$.
    Thus the extremal ray $R'$ is generated by $D_\lambda'$ (or $D_{\lambda,0}'$).
    Therefore $R'$ is of type D$_2$, the corresponding linear system is $|H_V'-F_V'|$, and $V'$ has a quadric bundle structure again.
    This is (\ref{thm:quad}.4).
\cons{7}
    One has $X=\P[0^3,1]$.
    The general member $V$ of $|2H+F|$ is smooth, $-K_V\sim 2H_V$, and $\rho(V)=2$ by (\ref{thm:ampleLeff}).
    Let $L=\P[0^3] \cong\P^2\times\P^1$ be a sub-$\P^2$-bundle in $X$.
    The intersection $V\cap L$ is defined by an equation $q_0(x)t_0+q_1(x)t_1$, where $x=[x_0:x_1:x_2]$ and $t=[t_0:t_1]$ are the homogeneous coordinates of $\P^2$ and $\P^1$, respectively.
    For any $\alpha\in A=\bigl\{ \alpha=[\alpha_0:\alpha_1:\alpha_2] \bigm| q_0(\alpha)=q_1(\alpha)=0 \bigr\}$, $s_\alpha=\{\alpha\}\times\P^1$ is in $V$ and $\mu(s_\alpha)=0$, hence the class $[s_\alpha]$ generates an edge $R\subset\ol{NE}(V)$.
    Similarly as in (\ref{subsec:quad}.3), each curve whose class is in $R$ is the same $s_\alpha$ as stated above.
    For general $V$, $\sharp A=4$, hence $V$ has just four $K$-trivial curves $s_1,\dots,s_4$ and has the $(-F_V)$-flop $V'$.
    The bihomogeneous coordinates analysis implies that $\ol{V}\subset\P^4$ is a cubic $3$-fold with $4$ ODP's on a line in $\ol{V}$.

    We will fix an extremal ray of $\ol{NE}(V')$ in a manner similar to (\ref{subsec:quad}.2).
    Using (\ref{eqn:slope}) for $|2H_V-F_V|$ with $Bs|2H_V-F_V|=\bigcup_{i=1}^4s_i$, we have $\mu(C')\geq 1/2$ for any curve $C'\subset V'$.
    On the other hand, for general sub-$\P^1$-bundle $S=\P[0^2]\cong\P^1\times\P^1$, the intersection curve $C=V\cap S$ has a slope $\mu(C)=1/2$ and does not meet the sections $\bigcup_{i=1}^4s_i$.
    Hence the extremal ray $R'$ is generated by the image of $C$, and $\mu(R')=\mu(C')=1/2$.
    The contraction morphism of $R'$ is defined by a linear system $|2H_{V'}-F_{V'}|$; $\dim|2H_{V'}-F_{V'}|=5$ and $(2H_{V'}-F_{V'})^3=4$.
    Thus $V'$ has an extremal ray of type $E_2$, and has a morphism to $B_4$ which is a blowing-up at a point of $B_4$, where $B_4$ is a Fano $3$-fold of index $2$, a complete intersection of two quadrics in $\P^5$.
    Moreover, the exceptional divisor $E'$ for the morphism is the unique member of $|H_{V'}-F_{V'}|$, and this is the strict transform of $V\cap L$ by the $(-F_V)$-flop.
    This is (\ref{thm:quad}.5).
\cons{8}
    In this case, we have $X=\P[0,1^3]$, and $V$ is a member of $|2H|$ with $-K_V\sim 2H_V-F_V$ and $\rho(V)=2$ by (\ref{thm:01n}).
    The minimal section $s_0=\P[0] \subset X$ is not contained in $V$, because if it was in $V$ then $-K_V$ would not be nef by $(-K_V\cdot s_0)=-1$.
    The general member $V\in |2H|$ is smooth and does not contain $s_0$.
    Now consider a family of ruled surfaces $S_\lambda=\P[1,0] \subset X$ parameterized by $\lambda\in\P^2$.
    Each $S_\lambda$ contains the minimal section $s_0$, and we have
  $\Pic(S_\lambda) \cong \Z[s_0]\op\Z[l] \cong \Z[H_\lambda]\op\Z[F_\lambda]$ for restrictions $H_\lambda,F_\lambda$ of $H,F$ to $S_\lambda$.
    Here $l\sim F_\lambda$ and $H_\lambda\sim s_0+l$ on $S_\lambda$.
    If the intersection curve $D_\lambda=V\cap S_\lambda\equiv 2[s_0]+2[l]$ is decomposed into $D_{\lambda,1}+D_{\lambda,2}$, then one of the followings holds:
$$
  D_{\lambda,1} \equiv 2[s_0]+[l]\quad{\rm and}\quad D_{\lambda,2}\equiv [l], \quad{\rm or}$$
$$
  D_{\lambda,1} \equiv D_{\lambda,2}\equiv [s_0]+[l],
$$
since $s_0$ cannot be a component of $D_\lambda$.
    The first case cannot occur because $D_{\lambda,1}$ has $s_0$ as an irreducible component by $(D_{\lambda,1}\cdot s_0)_{S_\lambda}=-1$.
    The second case occurs for $\lambda$ on a plane curve $C_4$ of degree 4; each $D_{\lambda,i}$ is no longer decomposable.
    The slope of an irreducible $D_\lambda$ or of $D_{\lambda,i}$ in the second case is $\mu(D_\lambda)=\mu(D_{\lambda,i})=1$.
    Any curve of slope $\leq 1$ is in a ruled surface $S_\lambda$ by (\ref{thm:slope}), thus $\ol{NE}(V)$ has an edge $R$ of slope 1 generated by $[D_{\lambda,i}]$ (or $[D_\lambda]$).
    Since $-K_V\sim 2H_V-F_V$, $R$ is an extremal ray and $V$ is a Fano $3$-fold.
    The contraction morphism associated to $R$ is defined by $|H_V-F_V|$ with $\dim|H_V-F_V|=3$.
    Hence $V$ has a conic bundle structure over $\P^2$ with a discriminant locus $C_4$ of degree $4$.
    This is (\ref{thm:quad}.6).
\cons{9}
    One has $X=\P[0^4] \cong \P^3\times\P^1$.
    Each member $V$ of $|2H+2F|$ is defined by an equation $q_0(x)t_0^2+q_1(x)t_0t_1+q_2(x)t_1^2$ as in (\ref{subsec:quad}.3); $\rho(V)=2$ by (\ref{thm:ampleLeff}).
    Consider the set $A = \bigl\{ \alpha=[\alpha_0:\alpha_1:\alpha_2:\alpha_3] \bigm| q_0(\alpha)=q_1(\alpha)=q_2(\alpha)=0 \bigr\}$ in $\P^3$, which consists of eight points for general $V$.
    For each $\alpha\in A$, the corresponding section $s_\alpha=\{\alpha\}\times\P^1$ is in $V$ and $\mu(s_\alpha)=0$, hence the class generates an edge of $\ol{NE}(V)$.
    Similarly as in (\ref{subsec:quad}.3) one can see that the $K$-trivial curves on $V$ are the same $s_\alpha$ as stated above.
    Hence, for the general member $V\in |2H+2F|$, $V$ is a smooth weak Fano $3$-fold with just eight $K$-trivial curves and has the $(-F_V)$-flop $V'$.
    The variety $V'$ is again a weak Fano $3$-fold with only one extremal ray $R'$ of $\ol{NE}(V')$.
    To find $\mu(R')=1/2$, we use (\ref{eqn:slope}) as in (\ref{subsec:quad}.5) for $|2H_V-F_V|$.
    Then, we can see that the type of $R'$ is $D_2$ and that the contraction morphism is defined by $|2H_{V'}-F_{V'}|$.
    The bihomogeneous coordinates analysis implies that $\ol{V}$ is the double covering of $\P^3$ branched along a quartic surface defined by $q_1^2-q_0q_2$ with $8$ ODP's, where $q_i$ are quadratic forms.
    This is (\ref{thm:quad}.7).
\cons{10}
    We have $X=\P[0^2,1^2]$.
    The general member $V$ of $|2H+F|$ is a smooth variety with $-K_V\sim 2H_V-F_V$, and does not contain minimal sections $s=\P[0]$; $\rho(V)=2$ by (\ref{thm:ampleLeff}).
    Let $S=\P[0^2] \cong \P^1\times\P^1$ be a ruled surface in $X$ and consider the intersection curve $C=V\cap S$; $C$ is an irreducible curve with $(-K_V\cdot C)=0$ and $\mu(C)=1/2$.
    Using a similar argument as in (\ref{thm:slope}) shows that this curve is the only $K$-trivial curve in $V$.
    Hence $V$ is a weak Fano $3$-fold with only one $K$-trivial curve and has the $(-F_V)$-flop $V'$; the $3$-fold $V'$ is again a weak Fano $3$-fold, and has the unique extremal ray $R'$ of $\ol{NE}(V')$.
    We have $\mu(R')=1$, because the image of the intersection curve $V\cap\P[0,1]\subset\P[0^2,1^2]$ is of slope $=1$, and because every curve in $V'$ is of slope $\geq1$ by using (\ref{eqn:slope}) for $|H_V-F_V|$ with $Bs|H_V-F_V|=C$.
    Hence the contraction morphism is defined by $|H_V'-F_V'|$ with $\dim|H_V'-F_V'|=1$ and $(H_V'-F_V')\cdot(-K_V')^2=4$, and $R'$ is of type $D_1$.
    Therefore $V'$ has a del Pezzo fibration of degree $4$.
    $\ol{V}$ is a singular Fano $3$-fold of index $1$ of genus $9$ with only one ODP.
    This is (\ref{thm:quad}.8).
\cons{12}
    In this case, $X=\P[0,1^3]$.
    The general member $V$ of $|2H+F|$ is a smooth variety with $-K_V\sim 2H_V-2F_V$, $\rho(V)=2$, and does not contain the minimal section $s_0$.
    We consider a family of ruled surfaces $S_\lambda=\P[0,1] \subset X$ parameterized by $\lambda\in \P^2$.
    Similarly as (\ref{subsec:quad}.6), $D_\lambda=V\cap S_\lambda \equiv 2[s_0]+3[l]$ may be decomposed into
$$
  D_{\lambda,1} \equiv  [s_0]+ [l] \quad {\rm and} \quad D_{\lambda,2} \equiv  [s_0]+2[l], \quad{\rm or}
$$
$$
  D_{\lambda,1} \equiv 2[s_0]+2[l] \quad {\rm and} \quad D_{\lambda,2} \equiv [l].
$$
    Direct calculation using the bihomogeneuos coordinates shows that the first case occurs at just sixteen $\lambda$'s, and that the second case occurs at just two $\lambda$'s; $(-K_V\cdot D_{\lambda,1})=0$ in each case.
    (The second case is easy but the first case is slightly complicated.)
    In order to find $K$-trivial curves in $V$, we have only to consider the curves stated above for the sake of (\ref{thm:slope}).
    Hence $V$ has just sixteen $K$-trivial curves as sections and two $K$-trivial curves as bisections, and has the $(-F_V)$-flop $V'$ with a unique extremal ray.
    The bihomogeneous coordinates calculus give us that there exist the members of $|4H_V|$ having $5$ fibers, i.e., the members of $|4H_V-5F_V|$, and that $Bs|4H_V-5F_V|$ is the union of all the $K$-trivial curves above.
    Similarly as in the other cases, we see that the type of the ray is $D_2$ and that $V'$ has a quadric bundle structure defined by $|4H_{V'}-5F_{V'}|$.
    Moreover, $\ol{V}$ is a hypersurface of degree $6$ with $18$ ODP's in the weighted projective space $\P(1^3,2,3)$.
    This is (\ref{thm:quad}.9).

\section{Del Pezzo fibrations of degree 2} \label{sec:df2}

\newpar \Label{subsec:setup2}
    In this section, we treat del Pezzo fibrations of degree $2$ over $\P^1$, and derive Theorem(\ref{thm:df2}) and Supplement(\ref{thm:df2s}).
    Let $V$ be a del Pezzo fibration of degree $2$ over $\P^1$ which is a double covering of $X=\P[0,a_1,a_2]$, $0\leq a_1\leq a_2$, branched along $B\sim 4H+2kF$ (see (\ref{subsec:setup})):
$$\begin{array}{ccl}
  V & \Longarrow{\psi}{2:1} & X = \P[0,a_1,a_2] \vspace{5pt}\\
  \Downarrow{\varphi} & \Swarrow{\pi} & \qquad B \sim 4H + 2kF \vspace{5pt}\\
  \P^1
\end{array}$$
    Assume that $V$ is a weak Fano $3$-fold with only finite $K$-trivial curves and with $\rho(V)=2$.

\newpar
    Under the situation (\ref{subsec:setup2}), the anti-canonical divisor of $V$ is
$$
  -K_V\sim H_V-(\sum_{i=1}^2a_i+k-2)F_V,
$$
from (\ref{eqn:can}).
    The argument in Section \ref{sec:prel} implies the following two inequalities:
\Eq(nef2){
  -(a_1+a_2+k-2) \geq 0,}
\Eq(big2){
  2(6-2a_1-2a_2-3k) > 0,}\\
from (\ref{eqn:nef}) and (\ref{eqn:big}), respectively((\ref{eqn:eff}) gives trivial conditions).
    If $a_2+2k<0$, then every member of $|4H+2kF|$ has singularities along $s_0$.
    Hence, if the branch $B$ is smooth, it needs to hold that
$$
  a_2+2k \geq 0.
$$

\newpar
   One has $3a_1-4 \leq 2a_1+a_2-4 \leq 0$  from (\ref{eqn:nef2}), hence $a_1=0$ or $1$.
    Note that the equality in (\ref{eqn:nef2}) does not hold when $a_1=0$.
    Indeed, if the equality held, there were a family of $K$-trivial curves in $V$.

\newpar
    Thus we have the six possibilities of numbers $a_1,a_2,k$, and $(-K_V^3)$ as in Table \ref{tbl:df2}.
\begin{table}[htbp]
\caption{}\label{tbl:df2}
$$\begin{array}{c@{\qquad}ccc@{\qquad}l@{\qquad}l@{\qquad}c}
Nos.& a_1& a_2&  k &\qquad B&\quad -K_V & (-K_V^3) \\
  1 &  0 &  0 &  0 & 4H     & H_V+2F_V  &  12 \\
  2 &  0 &  2 & -1 & 4H-2F  & H_V+ F_V  &  10 \\
  3 &  0 &  1 &  0 & 4H     & H_V+ F_V  &   8 \\
  4 &  0 &  0 &  1 & 4H+2F  & H_V+ F_V  &   6 \\
  5 &  1 &  2 & -1 & 4H-2F  & H_V       &   6 \\
  6 &  1 &  1 &  0 & 4H     & H_V       &   4
\end{array}$$
\end{table}
    In the rest of this section, we consider the realization of each possibility in Table \ref{tbl:df2}.

\newpar
    We now exclude two possibilities Nos.1 and 2 on Table \ref{tbl:df2}.
\ex{1}
    This is the case that $V$ is the double covering of $X=\P[0^3] \cong \P^2\times\P^1$ branched along $B\sim 4H$.
    In this case, $V$ is isomorphic to a product of a del Pezzo surface of degree $2$ by $\P^1$, hence its Picard number is $9$, which is an excluded case.
\ex{2}
    In this case, $\E=\O^{\op2}\op\O(2)$ and $V$ is the double covering of $X=\P[0^2,2]$ branched along $B\sim 4H-2F$.
    Since $H^0\bigl(X, \O_X(4H-2F)\bigr) = H^0\bigl(\P^1, S^4(\E)\otimes\O(-2)\bigr)$, every member of $|4H-2F|$ has an irreducible component $L\sim H-2F$, which contradicts that the branch divisor is smooth.
    Thus this case is excluded.

\newpar \Label{subsec:df2}
    Next, we construct the del Pezzo fibrations with each condition in the rest cases on Table \ref{tbl:df2}, and show Theorem(\ref{thm:df2}) and Supplement(\ref{thm:df2s}).
\cons{3}
    In this case, we have $X=\P[0^2,1]$.
    General member $B\in |4H|$ is smooth, hence the double covering $V$ of $X$ branched along $B$ is smooth and a del Pezzo fibration of degree $2$.
    Since $-K_V\sim H_V+F_V$ is ample, $V$ is a Fano $3$-fold.
    Consider $X$ as an ample effective divisor of $Y=\P[0^4] \cong \P^3\times\P^1$, and $B$ as the intersection of $X$ with $D$, where $D$ is a member of $|4H_Y|$.
    Then $V$ is an ample effective divisor of $W=B_2\times\P^1$ which is the double covering of $Y$ branched along $D$, where $B_2$ is the double covering of $\P^3$ branched along a quartic curve $C_4\subset \P^3$, $\rho(B_2)=1$.
    Hence $\rho(V)=\rho(W)=2$ by (\ref{thm:ampleLeff}).
    The extremal ray $R$ not corresponding to $\varphi : V\to\P^1$ is generated by $\sigma$, where $\sigma$ is the strict transform of the section $s=\P[0]$ of $\pi : X=\P[0^2,1] \to \P^1$.
    The contraction morphism corresponding to $R$ is defined by $|H_V|$ and maps to $\P^3$ with degree $2$.
    By taking the Stein factorization, this morphism factors through the double covering $B_2$ of $\P^3$ branched along a quartic surface;
    the exceptional divisor $E$ of the morphism is a unique member of $|H_V-F_V|$;
    and $E$ contracts to an elliptic curve on $B_2$.
    This is (\ref{thm:df2}.1).
\cons{4}
    In this case, $V$ is the double covering of $X=\P[0^3] \cong \P^2\times\P^1$ branched along $B\sim 4H+2F$.
    General member $B\in|4H+2F|$ is smooth, hence $V$ is smooth and is a del Pezzo fibration of degree $2$.
    We have $\rho(V)=\rho(X)=2$, because the branch locus $B$ is ample (cf. \cite{MM} p.508).
    Since $-K_V\sim H_V+F_V$, $V$ is a Fano $3$-fold.
    The extremal ray $R$ not corresponding to $\varphi$ is generated by $\sigma$ as above.
    The contraction morphism is defined by $|H_V|$, and maps $V$ onto $\P^2$.
    Thus $V$ has a conic bundle structure over $\P^2$ with discriminant locus $\Delta$ of degree $8$.
    This is (\ref{thm:df2}.2).
\cons{5}
    In this case, $X=\P[0,1,2]$.
    One can see that general member $B\in|4H-2F|$ is smooth by using the local coordinate system on $X$, hence the double covering $V$ of $X$ branched along $B$ is a weak Fano $3$-fold with $-K_V\sim H_V$ and is a del Pezzo fibration of degree $2$.
    Let $\sigma_0$ be the strict transform of the unique minimal section $s_0\subset X$.
    Then $B$ contains $s_0$ by $(s_0\cdot B)=-2$, hence $\sigma_0\cong s_0\cong\P^1$; $\sigma_0$ is only one $K$-trivial curve, and is contained as a section.
    Thus $V$ has the $(-F_V)$-flop $V'$ with an extremal ray $R'$; $\mu(R')=1$ by a similar argument to (\ref{subsec:quad}.2).
    The strict transform $E'$ of the unique member $E\in|H_V-2F_V|$ is the double covering of $\P^2$ branched along a conic, hence is a quadric surface(not necessary smooth).
    Consider the morphism defined by $|6H_{V'}-6F_{V'}|$.
    The surface $E'$ is the exceptional locus of the morphism and contracted to a double point; the image of the morphism is a singular Fano $3$-fold $W'$ of index $2$ of degree $1$ with a double point.
    In particular, $\rho(V)=\rho(V')=\rho(W')+1=2$.

    More precisely, we consider as follows. 
    In this case, $V$ can be regarded as a hypersurface in the weighted projective space bundle $X=\proj\S$ associated to $\S=\bigoplus_{d\geq0} \S_d$ such that $\S_d=\bigoplus_{i+2j=d} S^i(\O\op\O(1)\op\O(2))\otimes S^j(\O(1))$.
    Note that $X$ has the bihomogeneous coordinates ring $R=\C[x,y,z,u,t_0,t_1]$ with bidegrees $\deg x=(1,0)$, $\deg y=(1,-1)$, $\deg z=(1,-2)$, $\deg u=(2,-1)$, and $\deg t_j=(0,1)$, and that $V$ is defined by a bihomogeneous polynomial $F\in R$ of bidegree $(4,-2)$.
    Choosing the generators of the bigraded ring $R$, we may assume that $F=u^2+x^3z+cx^2y^2+xf(y,z,t_j)+g(y,z,t_j)$ with general $f, g\in R$ of $\deg f=(3,-2)$, $\deg g=(4,-2)$, because $V$ is general.
    Straightforward calculation shows that $V$ is smooth and that $V$ has the unique flopping curve $s_0=\{y=z=u=0\}$.
    The birational map $\chi:X\cdots\to Y'$ defined by the linear system $|6H_X-6F_X|$ is determined by the ratio of monomials with bidegree $(6,-6)$.
    Consider the graded ring $R'=\C[x'_0,x'_1,x'_2,y',z']$ with degrees $\deg x'_i=1$, $\deg y'=2$ and $\deg z'=3$.
    Then the map $\chi$ is corresponding to the ring homomorphism $R'\to R$ defined by $(x'_0,x'_1,x'_2,y',z')\mapsto(y,zt_0,zt_1,xz,zu)$, and $Y'={\rm Proj}\:R'$ is the weighted projective space with weights $(1,1,1,2,3)$.
    Moreover, $V$ is mapped onto $W'\subset Y'$ defined by the weighted homogeneous polynomial $F'={z'}^2+{y'}^3+c{y'}^2{x'_0}^2+y'f'(x'_i)+g'(x'_i)$ of degree $6$, i.e., a del Pezzo $3$-fold of degree $1$.
    We can see that $f'(x'_i)=\sum_{k=0}^3 f'_k(x'_1,x'_2){x'_0}^k$ and $g'(x'_i)=\sum_{k=0}^4 g'_k(x'_1,x'_2){x'_0}^k$, and that $W'$ has a double point at $[x'_0:x'_1:x'_2:y':z']=[1:0:0:0:0]$.
    Therefore, the rational map $\chi|_V:V\cdots\to W'$ is the composite of the flop $V\cdots\to V'$ along $s_0$ and the reverse of the blowing-up of $W'$ at the double point.
    By the description using the bihomogeneous coordinates above, we can see that the anti-canonical model $\ol{V}$ is the double covering of the cone $C(\Sigma_1)\subset\P^5$ over the rational ruled surface $\Sigma_1\cong\P[0,1]$ branched along $B$.
    The branch locus $B$ is an irreducible component of the quartic hypersurface section $B+F_1+F_2$ through the vertex of $C(\Sigma_1)$, where $F_i=C(f)\subset C(\Sigma_1)$ is the cone over a fiber $f$ of $\Sigma_1$.
    The singularity of $\ol{V}$ is one ODP, the inverse image of the vertex of the cone.
    This case is (\ref{thm:df2}.3).
\cons{6}
    In this case, $X=\P[0,1^2]$.
    General member $B\in|4H|$ is smooth, hence $V$ is smooth and a del Pezzo fibration of degree $2$.
    Since $-K_V\sim H_V$, $V$ is a weak Fano $3$-fold.
    As similar to (\ref{subsec:df2}.1), consider $X$ as a complete intersection of ample divisors in $Y=\P[0^5] \cong \P^4\times\P^1$, and $V$ as the same in a double covering $W$ of $Y$.
    Here $W$ is the direct product of a del Pezzo $4$-fold of degree $2$ and $\P^1$; $\rho(W)=2$; hence $\rho(V)=2$ by (\ref{thm:ampleLeff}).
    Let $\sigma_0$ be the pulled-back of the unique minimal section $s_0\subset X$.
    Because $(s_0\cdot B)=0$, $\sigma_0$ is decomposed into disjoint two curves each of which is isomorphic to $\P^1$.
    These curves are nothing but all of $K$-trivial curves in $V$, thus $V$ has the $(-F_V)$-flop $V'$, which is also a weak Fano $3$-fold with an extremal ray $R'$.
    As similar to (\ref{subsec:quad}.2), $\mu(R')=1$ by (\ref{eqn:slope}); hence the contraction morphism associated to $R'$ is defined by $|H_{V'}-F_{V'}|$; $V'$ is a del Pezzo fibration of degree $2$ over $\P^1$.
    The anti-canonical model $\ol{V}$ is the double covering of the singular quadric $3$-fold $\Q_0\subset\P^4$ branched along a quartic hypersurface section on $\Q_0$.
    The singularities of $\ol{V}$ consist of two ODP's, which are the inverse image of the vertex of $\Q_0$.
    This is (\ref{thm:df2}.4).

\section{Del Pezzo fibrations of degree 3} \label{sec:df3}

\newpar \Label{subsec:setup3}
    In this section, we treat del Pezzo fibrations of degree $3$ over $\P^1$,
i.e., cubic bundles over $\P^1$, and derive Theorem(\ref{thm:df3}) and Supplement(\ref{thm:df3s}).
    Let $V$ be a del Pezzo fibrations of degree $3$ over $\P^1$ embedded in $X=\P[0,a_1,a_2,a_3]$, $0\leq a_1\leq a_2\leq a_3$, such that $V\in |3H+kF|$ (see (\ref{subsec:setup})):
$$\begin{array}{ccl}
  V & \Longarrow{\psi}{} & X = \P[0,a_1,a_2,a_3] \vspace{5pt}\\
  \Downarrow{\varphi} & \Swarrow{\pi} & \qquad V \in |3H + kF| \vspace{5pt}\\
  \P^1
\end{array}$$
    Assume that $V$ is a weak Fano $3$-fold with only finite $K$-trivial curves and with $\rho(V)=2$.

\newpar
    Under the situation (\ref{subsec:setup3}), the anti-canonical divisor of $V$ is
$$
  -K_V\sim H_V-(\sum_{i=1}^3a_i+k-2)F_V
$$
by (\ref{eqn:can}).
    The argument in Section \ref{sec:prel} implies the following three inequalities:
\Eq(nef3){
  6-3a_2-3a_3-2k \geq 0,}
\Eq(eff3){
  3a_1+k \geq 0,}
\Eq(big3){
  9-3a_1-3a_2-3a_3-4k > 0,}\\
from (\ref{eqn:nef}), (\ref{eqn:eff}) and (\ref{eqn:big}), respectively.

\newpar
    From (\ref{eqn:nef3}) and from $a_1\leq a_2\leq a_3$, we have $3\geq 3a_1+k$.
    Hence we treat only the following four cases from (\ref{eqn:eff3}):
\Eq(3a0){
  k=-3a_1\phantom{{}+1} \quad {\rm and} \quad 2+2a_1\geq a_2+a_3, \qquad\quad}
\Eq(3a1){
  k=-3a_1+1 \quad {\rm and} \quad 4+6a_1\geq3a_2+3a_3, \qquad}
\Eq(3a2){
  k=-3a_1+2 \quad {\rm and} \quad 2+6a_1\geq3a_2+3a_3, \quad {\rm or}}
\Eq(3a3){ 
  k=-3a_1+3 \quad {\rm and} \quad 2a_1\geq a_2+a_3. \qquad\qquad\quad}

    In the case (\ref{eqn:3a0}), the possible values of a quadruple $(a_1, a_2, a_3, k)$ are enumerated as follows using a non-negative integer $b$:
$$\begin{array}{r@{\:}l}
  (a_1,a_2,a_3,k) &= (b,b,  b,  -3b), \\
  (a_1,a_2,a_3,k) &= (b,b,  b+1,-3b), \\
  (a_1,a_2,a_3,k) &= (b,b,  b+2,-3b), \quad {\rm and} \\
  (a_1,a_2,a_3,k) &= (b,b+1,b+1,-3b).
\end{array}$$
    If $b$ was positive, every member $V\in|3H+kF|$ had singularities along the minimal section $s_0=\P[0] \subset X$, except for the case $(a_1,a_2,a_3,k)=(1,1,3,-3)$.
    Hence $b=0$ with the one exception.

    In the case (\ref{eqn:3a1}), for any non-negative integer $b$, the following two cases occur:
$$\begin{array}{r@{\:}l}
  (a_1,a_2,a_3,k) &= (b,b,b,  1-3b), \quad {\rm and}\\
  (a_1,a_2,a_3,k) &= (b,b,b+1,1-3b).
\end{array}$$
    A similar argument as above shows that $b=0$ in the first case and that $b=0$ or $1$ in the second case.

    In the cases (\ref{eqn:3a2}) and (\ref{eqn:3a3}), we have $a_1=a_2=a_3$ and let the number $b$.
    From the smoothness of $V$, we have that $b=0$ or $1$.
    We can exclude the cases from (\ref{eqn:3a3}) because we have $b>1$ by (\ref{eqn:big3}).

\newpar
    Thus we obtain ten possibilities of numbers $a_1,a_2,a_3,k$, and $(-K_V^3)$ as in Table \ref{tbl:df3}.
\begin{table}[htbp]
\caption{}\label{tbl:df3}
$$\begin{array}{c@{\qquad}cccc@{\qquad}l@{\qquad}l@{\qquad}rc}
Nos.& a_1& a_2& a_3&  k &\qquad V&\quad -K_V &  D \qquad & (-K_V^3) \\
  1 &  0 &  0 &  0 &  0 & 3H    &  H_V+2F_V &   3s_0&  18 \\
  2 &  0 &  0 &  1 &  0 & 3H    &  H_V+ F_V &   3s_0&  12 \\
  3 &  1 &  1 &  3 & -3 & 3H-3F &  H_V      &   3s_0&  12 \\
  4 &  0 &  0 &  0 &  1 & 3H+ F &  H_V+ F_V & l+3s_0&  10 \\
  5 &  1 &  1 &  2 & -2 & 3H-2F &  H_V      & l+3s_0&  10 \\
  6 &  1 &  1 &  1 & -1 & 3H- F &  H_V      &2l+3s_0&   8 \\
  7 &  0 &  0 &  2 &  0 & 3H    &  H_V      &   3s_0&   6 \\
  8 &  0 &  1 &  1 &  0 & 3H    &  H_V      &   3s_0&   6 \\
  9 &  0 &  0 &  1 &  1 & 3H+ F &  H_V      & l+3s_0&   4 \\
 10 &  0 &  0 &  0 &  2 & 3H+2F &  H_V      &2l+3s_0&   2
\end{array}$$
\end{table}
    Here Nos.$1$, $2$, $3$, $7$, $8$ come from (\ref{eqn:3a0}), Nos.$4$, $5$, $9$ from (\ref{eqn:3a1}), and Nos.$6$, $10$ from (\ref{eqn:3a2}).

\newpar
    We first exclude three possibilities Nos.$1$, $3$, and $7$ on Table \ref{tbl:df3}.
\ex{1}
    This is the case that $V\in |3H|$ in $X=\P[0^4] \cong\P^3\times\P^1$.
    This $V$ is a trivial cubic bundle, i.e., the product of a del Pezzo surface of degree $3$ and $\P^1$, and $\rho(V)=8$.
    Hence this is an excluded case.
\ex{3}
    This is the case that $X=\P[0,1^2,3]$, $V\in |3H-3F|$. 
    The restriction of $V$ to $\P[0,1^2]\subset X$ is the union of three $\P[0,1]$'s.
    Each $\P[0,1]$ does not represent by the linear combination of $H_V$ and $F_V$ on $V$.
    Hence, $\rho(V)>2$ in this case.
\ex{7}
    In this case, $X=\P[0^3,2]$, $V\in |3H|$, and $-K_V\sim H_V$.
    For any intersection curve $C_\lambda=V\cap S_\lambda$, one has $(-K_V\cdot C_\lambda)=0$, where $S_\lambda=\P[0,2] \cong \P^1\times\P^1$ are ruled surfaces in $X$, parameterized by $\lambda\in\P^2$.
    Thus, there is a family of $K$-trivial curves $C_\lambda\subset V$, which contradicts our assumption.

\newpar \Label{subsec:df3}
    Now we construct the del Pezzo fibrations of degree $3$ with each condition in the rest cases on Table \ref{tbl:df3}, and show Theorem(\ref{thm:df3}) and Supplement(\ref{thm:df3s}).
\cons{2}
    In this case, we have that $X=\P[0^3,1]$.
    General member $V\in|3H|$ is a smooth Fano $3$-fold with $-K_V\sim H_V+F_V$, and $\rho(V)=2$ by (\ref{thm:0n1}).
    Let $L=\P[0^3] \cong \P^2\times\P^1$ be a sub-$\P^2$-bundle of $X$.
    Then the first projection $p : L\to\P^2$ maps the intersection $L\cap V$ onto a smooth cubic $C$ and each fiber $s=p^{-1}(c)$, $c\in C$, is $\P^1$.
    Moreover, $\mu(s)=0$ and the class $[s]$ generates an extremal ray of $\ol{NE}(V)$.
    The corresponding morphism is defined by $|H_V|$ with $\dim|H_V|=4$ and $(H_V)^3=3$.
    The morphism maps $V$ onto $B_3$ contracting the above fibers $s$; this is a blowing-up of $B_3$ along a smooth cubic.
    This is (\ref{thm:df3}.1).
\cons{4}
    In this case, $X=\P[0^4] \cong \P^3\times\P^1$.
    General member $V\in|3H+F|$ is smooth with $-K_V\sim H_V+F_V$, and $\rho(V)=2$ by (\ref{thm:ampleLeff}).
    Hence $V$ is a Fano 3-fold and the section of $\pi : V\to\P^1$ generates an extremal ray $R$ of $\ol{NE}(V)$.
    After fixing the coordinates $[x_0:x_1:x_2:x_3]\times[t_0:t_1]$ on $X=\P^3\times\P^1$, $V$ is defined by $c_0(x)t_0+c_1(x)t_1$ for cubic forms $c_0(x)$ and $c_1(x)$ in $x$.
    Then, there is a family of sections $s_\alpha$ parameterized by $\alpha\in C=\bigl\{ \alpha=[\alpha_0:\alpha_1:\alpha_2:\alpha_3] \bigm| c_0(\alpha)=c_1(\alpha)=0 \bigr\}$, where $C$ is a space curve of degree $9$ and of genus $10$.
    The contraction morphism corresponding to $R$ is defined by $|H_V|$ with $\dim|H_V|=3$ and $(H_V)^3=1$, and hence this is a morphism onto $\P^3$ which is the blowing-up along a curve $C$ of degree $9$ and genus $10$.
    This is (\ref{thm:df3}.2).
\cons{5}
    We have $X=\P[0,1^2,2]$.
    There is the minimal section $s_0=\P[0]$ on $X$.
    Let $V$ be a general member of $|3H-2F|$.
    From $(V\cdot s_0)_X=-2$, $V$ contains the section $s_0$.
    Straightforward calculation using the local coordinates on $X$ shows that $V$ is smooth and that the section $s_0$ is the unique $K$-trivial curve on $V$ because of $-K_V\sim H_V$.
    Therefore, $V$ is a weak Fano 3-fold, and has the $(-F_V)$-flop $V'$, which is again a weak Fano 3-fold with only one extremal ray $R'$.
    Similar as in (\ref{subsec:quad}.2), applying (\ref{eqn:slope}) for $|H_V-F_V|$, we see that $\mu(R')=1$.
    The contraction morphism is defined by the linear system $|H_{V'}-F_{V'}|$ with $\dim|H_{V'}-F_{V'}|=\dim|H_V-F_V|=3$ and $(H_{V'}-F_{V'})^3=2$.
    This gives a morphism onto $\P^3$ with degree $2$, and contracts a divisor $D\sim H_{V'}-2F_{V'}$ to a line on $\P^3$.
    As similar to (\ref{subsec:df2}.1), this morphism factors through a Fano 3-fold $B_2$ of index $2$; the exceptional divisor is the unique member of $|H_{V'}-2F_{V'}|$, and contracts to a line on $B_2$.
    Note that the composed birational map $V\cdots\to B_2$ is defined by the linear system $|2H_V-2F_V|$, and given explicitly by using the bihomogeneous coordinates as (\ref{subsec:df2}.3).
    By choosing a line on $B_2$ and following the reverse process, $V$ with the above structure can be constructed, hence $\rho(V)=2$.
    Using the bihomogeneous coordinates, we can see that the anti-canonical model $\ol{V}$ is a divisor of the cone $C(\P[0^2,1])\subset\P^7$ over the $\P^2$-bundle $\P[0^2,1]$, and is the irreducible component of a cubic hypersurface $\ol{V}+P_1+P_2$ through the vertex of $C(\P[0^2,1])$.
    Here, $P_i$ is the cone $C(f)\subset C(\P[0^2,1])$ over a fiber $f$ of $\P[0^2,1]\to\P^1$.
    The singularity of $\ol{V}$ is one ODP, the vertex of $C(\P[0^2,1])$.
    This obtains (\ref{thm:df3}.3).
\cons{6}
    In this case, $X=\P[0,1^3]$.
    Let $s_0$ be the minimal section of $X$ and $V$ a general member of $|3H-F|$.
    Similarly as (\ref{subsec:df3}.3), $V$ is smooth and contains $s_0$ as its unique $K$-trivial curve.
    Hence there is the $(-F_V)$-flop $V'$.
    From (\ref{eqn:slope}), the extremal ray $R'$ of $\ol{NE}(V')$ is of slope $\mu(R')=1$.
    The contraction morphism is defined by $|H_{V'}-F_{V'}|$, and this is a conic bundle with the discriminant locus $\Delta$ of degree $7$.
    We have $\rho(V)=2$ by (\ref{thm:0k1n}).
    The anti-canonical model $\ol{V}$ is a divisor of the cone $C(\P^2\times\P^1)\subset\P^6$, and is the irreducible component of a cubic hypersurface $\ol{V}+P$ through the vertex of $C(\P^2\times\P^1)$.
    Here, $P=C(\P^2)\subset C(\P^2\times\P^1)$ is the cone over a fiber of $\P^2\times\P^1\to\P^1$.
    The singularity of $\ol{V}$ is one ODP, the vertex of the cone.
    This is (\ref{thm:df3}.4).
\cons{8}
    One has $X=\P[0^2,1^2]$.
    General member $V\in |3H|$ is smooth with $-K_V\sim H_V$, and $\rho(V)=2$ by (\ref{thm:0n11}).
    Considering the intersection $V\cap S$ of $V$ with a ruled surface $S=\P[0^2] \cong \P^1\times\P^1$ in $X$, one can see that there are three $K$-trivial curves on $V$, and hence that $V$ has the $(-F_V)$-flop $V'$.
    One obtains that the extremal ray $R'\subset\ol{NE}(V')$ is of slope $\mu(R')=1$ by (\ref{eqn:slope}) for $|H_V-F_V|$, and that $V'$ has a del Pezzo fibration of degree $3$, which is defined by $|H_{V'}-F_{V'}|$.
    The anti-canonical model $\ol{V}\subset\P^5$ is the complete intersection of a cubic $4$-fold and a singular quadric $4$-fold $\Q_0$ of rank $4$, i.e., a singular quadric $4$-fold having singularities along a line.
    The singularities of $\ol{V}$ consict three points, the intersection of the cubic $4$-fold and the singular locus of $\Q_0$.
    This is (\ref{thm:df3}.5).
\cons{9}
    In this case, we have $X=\P[0^3,1]$.
    General member $V\in|3H+F|$ is smooth, $-K_V\sim H_V$, and $\rho(V)=2$ by (\ref{thm:ampleLeff}).
    Let $L=\P[0^3] \cong\P^2\times\P^1$ be a sub-$\P^2$-bundle of $X$.
    Since $V\cap L$ can be defined by the equation $c_0(x)t_0+c_1(x)t_1$ for the coordinates $[x_0:x_1:x_2]\times[t_0:t_1]$ on $L$, there are nine minimal sections on $V$ corresponding to $c_0=c_1=0$, and the $K$-trivial curves on $V$ are nothing but these sections.
    Since $Bs|3H_V-F_V|$ is the union of the nine minimal sections, the unique extremal ray $R'$ of the $(-F_V)$-flop $V'$ of $V$ is of slope $1/3$ by (\ref{eqn:slope}) for $|3H-F|$.

    We now consider the ruled surfaces $S_\gamma=\P[0^2] \cong \P^1\times\P^1$ in $X$ parameterized by $\gamma\in\P^2$.
    The locus swept out by these $S_\gamma$ is the sub-$\P^2$-bundle $L$.
    The strict transforms $C_\gamma'$ of curves $C_\gamma=S_\gamma\cap V$ are to be contracted by the contraction morphism associated to $R'$, because $\mu(C_\gamma')=1/3$.
    Each curve $C_\gamma$ meets a fixed curve $C_{\gamma_0}$, since the intersection $S_\gamma\cap S_{\gamma_0}$ is a line and meets $V$ at a point.
    Denote by $E$ the locus swept out by $C_\gamma$, which is the intersection $L\cap V$ linearly equivalent to $H_V-F_V$, and denote by $E'$ its strict transform on $V'$.
    Then, $E'$ is the exceptional devisor on $V'$ contracted to a point, because each curve $C_\gamma$ contracts to a point and meets the fixed curve $C_{\gamma_0}$.
    Each intersection $F_V\cap L=F\cap L\cap V$ is a cubic curve, so the surface $E$ has an elliptic fiber structure.
    The elliptic fiber of $E$ meets the nine minimal sections.
    These sections are contracted by the $(-F_V)$-flop.
    Now consider the first projection $p_1 : L=\P^2\times\P^1\to\P^2$.
    The restriction $p_1|_E : E\to\P^2$ is surjective.
    The images of the elliptic fibers intersects each other at only the nine points which are the images of the minimal sections, hence $p_1|_E(E)=\P^2$ is isomorphic to the image $E'$ by the $(-F_V)$-flop.
    Since $E'$ is linearly equivalent to $H_{V'}-F_{V'}$, the normal bundle of $E'\subset V'$ is linearly equivalent to $(H_{V'}-F_{V'})|_{E'}$.
    Because the pulled-back of a line in $E'\cong\P^2$ by $p_1|_E$ is linearly equivalent to $3s_0+l$ in $E$,
$$
  \deg\bigl((H_{V'}\!-\!F_{V'})|_{E'}\bigr)=
  \bigl((H_V\!-\!F_V)|_E\cdot(3s_0+l)\bigr)_E=
  \bigl((H_V\!-\!F_V)\cdot(3s_0+l)\big)=-2.
$$
    Therefore the extremal ray $R'$ is actually of type $E_5$.
    The anti-canonical model $\ol{V}$ is a singular quartic $3$-fold in $\P^4$ which contains a plane and nine ODP's on the plane.
    This case obtains (\ref{thm:df3}.6).
\cons{10}
    In this case, $X=\P[0^4] \cong \P^3\times\P^1$.
    General member $V\in|3H+2F|$ is smooth and $\rho(V)=2$.
    Using the standard coordinates $(x,t)$ of $X$ as in (\ref{subsec:df3}.2), $V$ is defined by the equation $c_0(x)t_0^2+c_1(x)t_0t_1+c_2(x)t_1^2$.
    There exist $27$ minimal sections on $V$, and $K$-trivial curves on $V$ are only these sections.
    Since $Bs|3H_V-F_V|$ is the union of the $27$ minimal sections, similar argument as others shows that the $(-F_V)$-flop $V'$ of $V$ is a del Pezzo fibration of degree $3$ over $\P^1$, and that the fibration is defined by $|3H_{V'}-F_{V'}|$.
    The anti-canonical model $\ol{V}$ is the double covering of $\P^3$ branched along the surface $S$ defined by $c_1^2-4c_0c_2$ with $27$ ODP's corresponding to $c_0=c_1=c_2=0$.
    The singularities of $\ol{V}$ correspond to the singularities of $S$.
    This is (\ref{thm:df3}.7).

\section{Del Pezzo fibrations of degree 4} \label{sec:df4}

\newpar\Label{subsec:setup4}
    In this section, we treat del Pezzo fibrations of degree $4$ over $\P^1$, and derive Theorem(\ref{thm:df4}) and Supplement(\ref{thm:df4s}).
    Let $V$ be a del Pezzo fibration of degree $4$ over $\P^1$ embedded in $X=\P[0,a_1,a_2,a_3,a_4]$, $0\leq a_1\leq a_2\leq a_3\leq a_4$, such that $V=W_1\cap W_2$ where $W_i\in |2H+k_iF|$ for $i=1,2$ (see (\ref{subsec:setup})):
$$\begin{array}{rcl}
  W_1\cap W_2 = V & \Longarrow{\psi}{} & X = \P[0,a_1,a_2,a_3,a_4] \vspace{5pt}\\
  \Downarrow{\varphi} & \Swarrow{\pi} & \qquad W_i \in |2H + k_iF| \vspace{5pt}\\
  \P^1
\end{array}$$
    Assume that $V$ is a weak Fano 3-fold with only finite $K$-trivial curves and with $\rho(V)=2$.

\newpar
    Under the situation (\ref{subsec:setup4}), the anti-canonical divisor of $V$ is
$$
  -K_V\sim H_V-(\sum_{i=1}^4a_i+k_1+k_2-2)F_V.
$$
    The argument in Section \ref{sec:prel} implies the following three inequalities:
\Eq(nef4){
  2(4-k_1-k_2-2(a_3+a_4)) \geq 0,}
\Eq(eff4){
  2(2(a_1+a_2)+k_1+k_2) \geq 0,}
\Eq(big4){
  2(12-\sum_{i=1}^4a_i-5(k_1+k_2)) > 0,}\\
from (\ref{eqn:nef}), (\ref{eqn:eff}) and (\ref{eqn:big}), respectively.
    If $a_4+k_i<0$ then $W_i$ contains $s_0=\P[0]$ as its singular locus.
Thus
\Eq(a4ki){
  a_4+k_i\geq 0 \qquad{\rm for}\quad i=1,2.}

\newpar
    From (\ref{eqn:nef4}) and (\ref{eqn:eff4}) we have
\Eq(aka){
  -2(a_1+a_2)\leq k_1+k_2\leq 4-2(a_3+a_4),}\\
hence $(a_4-a_1)+(a_3-a_2)\leq 2$.
    Thus we have only to treat the following three cases.

    {\it Case of $a_4=a_1+2$}.~
    In this case, $a_2=a_3$.
    Hence the equality in (\ref{eqn:nef4}) holds; and the number of $K$-trivial curves on $V$ is infinite, which is a contradiction.

    {\it Case of $a_4=a_1+1$}.~
    Let $b=a_1$.
    Then
$$\begin{array}{r@{\:}l}
  (a_1,a_2,a_3,a_4) &= (b,b,  b,  b+1), \\
  (a_1,a_2,a_3,a_4) &= (b,b,  b+1,b+1), \quad{\rm and} \\
  (a_1,a_2,a_3,a_4) &= (b,b+1,b+1,b+1). \\
\end{array}$$
\def\k{\ol{k}}
    In each subcase, we can restrict the value of $\k=k_1+k_2$ by (\ref{eqn:aka}):
$$\begin{array}{r@{\:}l@{\quad}l}
  \k &=   -4b, 1-4b, \\
  \k &=   -4b, \quad{\rm and} \\
  \k &= -2-4b,-1-4b,
\end{array}$$
according to the three values of quadruples $(a_1,a_2,a_3,a_4)$.
    Using (\ref{eqn:a4ki}), we have eleven possibilities for $a_i$ and $k_i$ (see Table \ref{tbl:df4}).

    {\it Case of $a_4=a_1$}.~
    As above, we have $a_1=a_2=a_3=a_4=b$ and $\k=-4b,1-4b,2-4b$, or $3-4b$.
    By (\ref{eqn:big4}) and (\ref{eqn:a4ki}), we restrict the values of $a_i$ and $k_i$ as in Table \ref{tbl:df4}.

\newpar
    We obtain the following possibilities of numbers $a_1,a_2,a_3,a_4,k_1,k_2$, and $(-K_V^3)$ by symmetricity of $k_1$ and $k_2$ as in Table \ref{tbl:df4}.
\begin{table}[htbp]
\caption{}\label{tbl:df4}
$$\begin{array}{cc@{\:\:}c@{\:\:}c@{\:\:}ccclll@{\;}r@{\;}c}
Nos.&a_1&a_2&a_3&a_4&k_1&k_2& \quad W_1&\quad W_2& \quad -K_V &  D \quad  &(-K_V^3) \\
  1 & 0 & 0 & 0 & 0 & 0 & 0 & 2H       &2H       &H_V\!+\!2F_V&       4s_0&  24 \\
  2 & 1 & 1 & 1 & 2 &-2 &-2 & 2H\!-\!2F&2H\!-\!2F&H_V\!+\! F_V&       4s_0&  24 \\
  3 & 0 & 1 & 1 & 1 &-1 &-1 & 2H\!-\! F&2H\!-\! F&H_V\!+\! F_V&       4s_0&  20 \\
  4 & 0 & 0 & 0 & 1 &-1 & 1 & 2H\!-\! F&2H\!+\! F&H_V\!+\! F_V&       4s_0&  16 \\
  5 & 0 & 0 & 0 & 1 & 0 & 0 & 2H       &2H       &H_V\!+\! F_V&       4s_0&  16 \\
  6 & 1 & 1 & 2 & 2 &-2 &-2 & 2H\!-\!2F&2H\!-\!2F&H_V         &       4s_0&  16 \\
  7 & 0 & 0 & 0 & 0 & 0 & 1 & 2H       &2H\!+\! F&H_V\!+\! F_V&2l\!+\!4s_0&  14 \\
  8 & 1 & 1 & 1 & 2 &-2 &-1 & 2H\!-\!2F&2H\!-\! F&H_V         &2l\!+\!4s_0&  14 \\
  9 & 1 & 1 & 1 & 1 &-1 &-1 & 2H\!-\! F&2H\!-\! F&H_V         &4l\!+\!4s_0&  12 \\
 10 & 0 & 1 & 1 & 1 &-1 & 0 & 2H\!-\! F&2H       &H_V         &2l\!+\!4s_0&  10 \\
 11 & 0 & 0 & 1 & 1 &-1 & 1 & 2H\!-\! F&2H\!+\! F&H_V         &       4s_0&   8 \\
 12 & 0 & 0 & 1 & 1 & 0 & 0 & 2H       &2H       &H_V         &       4s_0&   8 \\
 13 & 0 & 0 & 0 & 1 &-1 & 2 & 2H\!-\! F&2H\!+\!2F&H_V         &2l\!+\!4s_0&   6 \\
 14 & 0 & 0 & 0 & 1 & 0 & 1 & 2H       &2H\!+\! F&H_V         &2l\!+\!4s_0&   6 \\
 15 & 0 & 0 & 0 & 0 & 0 & 2 & 2H       &2H\!+\!2F&H_V         &4l\!+\!4s_0&   4 \\
 16 & 0 & 0 & 0 & 0 & 1 & 1 & 2H\!+\! F&2H\!+\! F&H_V         &4l\!+\!4s_0&   4 \\
 17 & 1 & 1 & 1 & 1 &-1 & 0 & 2H\!-\! F&2H       &H_V\!-\! F_V&6l\!+\!4s_0&   2 \\
\end{array}$$
\end{table}
    Here Nos.$1$, $7$, $9$, $15$, $16$, and $17$ come from the case of $a_4 = a_1$, and the rest from the case of $a_4 = a_1 + 1$ in the previous subsection.
    In the rest of this section, we consider the realization of each possibility on Table \ref{tbl:df4}.

\newpar
    We first exclude some possibilities: Nos.$1$, $2$, $3$, $4$, $11$, and $13$ on Table \ref{tbl:df4}.
\ex{1}
    This is the case that $V=W_1\cap W_2\subset X=\P[0^5]$, where $W_1, W_2\in |2H|$.
    But this case is to be excluded, because $V$ is isomorphic to the direct product of a del Pezzo surface of degree $4$ and $\P^1$ and has the Picard number $7$.
\ex{2}
    This is the case that $V=W_1\cap W_2\subset X=\P[0,1^3,2]$, where $W_1, W_2\in |2H-2F|$.
    In this case, for any two members of $|2H-2F|$, the intersection $V$ has singularities along the section $s_0=\P[0] \subset X$.
    Hence this is an excluded case.
\ex{3}
    In this case, $V=W_1\cap W_2\subset X=\P[0^2,1^3]$, where $W_1, W_2\in |2H-F|$.
    Let $T=\P[0^2]\cong \P^1\times\P^1$ be a ruled surface in $X$.
    The surface $T$ is contained in $V$ as a divisor, and has the intersection numbers $(T\cdot s_0)=-1$, $(T\cdot l)=1$ in $V$.
    If $\Pic(V)$ was generated by $H_V$ and $F_V$, then $T$ was linearly equivalent to $H_V-F_V$ for the sake of intersection numbers.
    But each member of the linear system $|H_V-F_V|$ is not $T$.
    Thus $\Pic(V)$ is not of rank $2$.
    This is a case to be excluded.
    Indeed, this $V$ has Picard number $3$ and $T$ is an exceptional divisor for the contraction morphism of an extremal ray.
\exs{4, 11, {\it and} 13}
    In these cases, $V$ has singularities of dimension $\geq 1$, which contradicts our assumption for smoothness.

\newpar \Label{subsec:df4}
    We next construct the del Pezzo fibrations of degree $4$ with our conditions in the rest cases on Table \ref{tbl:df4}, and derive Theorem(\ref{thm:df4}) and Supplement(\ref{thm:df4s}).
\cons{5}
    In this case, we have $X=\P[0^4,1]$.
    The intersection $V=W_1\cap W_2$ is smooth for general $W_1, W_2\in|2H|$, is a Fano $3$-fold with $-K_V\sim H_V+F_V$, and has $\rho(V)=2$ by (\ref{thm:0n1}).
    As similar to (\ref{subsec:quad}.1), considering the intersection with sub-$\P^3$-bundle $L\subset X$, we can see the following:
    the extremal ray $R\subset\ol{NE}(V)$ of slope 0 is of type $E_1$;
    the contraction morphism associated to $R$ is defined by $|H_V|$ with $\dim|H_V|=5$ and $(H_V^3)=4$;
    the exceptional divisor for the morphism is $V\cap L$, which is the unique member of $|H_V-F_V|$, is contracted to a quartic elliptic curve.
    This is (\ref{thm:df4}.1).
\cons{6}
    This case is $V=W_1\cap W_2\subset X=\P[0,1^2,2^2]$ for general members $W_1$, $W_2$ of $|2H-2F|$.
    Straightforward calculation using the local coordinates shows that $V$ is smooth with $-K_V\sim H_V$, and that $V$ is a weak Fano $3$-fold with only one $K$-trivial curve $s_0$, where $s_0=\P[0] \subset X$ is the minimal section.
    The $(-F_V)$-flop $V'$ of $V$ is also a weak Fano $3$-fold with one extremal ray $R'$.
    To fix $R'$, consider the linear system $|H_V-2F_V|$.
    Every member of $|H_V-2F_V|$ comes from a member of $|H-2F|$ on $X$, hence $Bs|H_V-2F_V|=Bs|H-2F|\cap V=s_0$.
    For any curve $C$ on $V$, let $a=(H_V\cdot C)$, $b=(F_V\cdot C)$, and $n_0=\#\{C\cap s_0\}$ as in (\ref{subsec:quad}.2).
    Then, $(H_V-2F_V\cdot C)=a-2b\geq 2n_0$, because $V\cap L$ has singularities along $s_0$ for any member $L\in|H-2F|$.
    The strict transform $C'$ of $C$ has a slope $\mu(C')=a/(b+n_0)\geq 2$ from (\ref{eqn:slope}).
    A curve $D=V\cap\P[0,1,2]\equiv 4[s_0]+4[l]$ is decomposed into $D_0+2s_0$ such that $D_0\equiv2[s_0]+4[l]$ and $D_0\cap s_0=\emptyset$, hence the strict transform $D_0'$ of $D_0$ has a slope $\mu(D_0')=2$.
    Thus we have $\mu(R')=2$.
    The contraction morphism of $R'$ is defined by $|H_{V'}-2F_{V'}|$, and is a quadric bundle over $\P^1$.
    This is nothing but the reverse process of (\ref{subsec:quad}.8) and $\rho(V')=\rho(V)=2$.
    This is (\ref{thm:df4}.2).
\cons{7}
    In this case, $V=W_1\cap W_2\subset\P[0^5] \cong\P^4\times\P^1$, where $W_1\in|2H|$ and $W_2\in|2H+F|$.
    For general $W_1$ and $W_2$, $V=W_1\cap W_2$ is a smooth Fano $3$-fold with $-K_V\sim H_V+F_V$.
    Since $W_1$ is the product of a quadric $3$-fold $\Q^3$ and $\P^1$, and since $V$ is an ample divisor on $W_1$, $\rho(V)=\rho(W_1)=2$ by (\ref{thm:ampleLeff}).
    The first projection $W_1\to\Q^3$ maps $V$ onto $\Q^3$ and the locus of the minimal sections in $V$ onto the curve $C$ which is a complete intersection of three quadrics in $\P^4$.
    The extremal ray $R$ of $\ol{NE}(V)$ is of slope $\mu(R)=0$, and the contraction morphism is defined by $|H_V|$.
    The morphism is nothing but the restriction of the first projection to $V$, which is the blowing-up of $\Q^3$ along $C$.
    This is (\ref{thm:df4}.3).
\cons{8}
    We have $X=\P[0,1^3,2]$.
    Calculation using the local coordinates on $X$ shows that $V=W_1\cap W_2$ is smooth for general members $W_1\in|2H-2F|$, $W_2\in|2H-F|$.
    Since $-K_V\sim H_V$ and $s_0\subset V$, $V$ is a weak Fano $3$-fold with only one $K$-trivial curve, where $s_0=\P[0] \subset X$ is the minimal section.
    Considering the intersection with any curve and $H_V-F_V$, we can see that the extremal ray $R'$ of the $(-F_V)$-flop $V'$ is of slope $\geq1$ by (\ref{eqn:slope}).
    Let $T_\gamma=\P[0,1^2]$ be sub-$\P^2$-bundles parameterized by $\gamma\in\P^2$, and $C_\gamma=V\cap T_\gamma$ for any $\gamma\in\P^2$.
    Each curve $C_\gamma$ is decomposed into $2s_0+C_{\gamma,1}+C_{\gamma,2}$, where $C_{\gamma,1}\equiv C_{\gamma,2}\equiv [s_0]+[l]$.
    The strict transform $C_{\gamma,i}'$ of $C_{\gamma,i}$ is of slope $1$, hence $\mu(R')=\mu(C_{\gamma,i}')=1$.
    Therefore, the contraction morphism of $R'$ is defined by $|H_{V'}-F_{V'}|$ with $\dim|H_{V'}-F_{V'}|=4$ and $(H_{V'}-F_{V'})^3=3$, the morphism contracts $C_{\gamma,i}'$
    The family of $C_\gamma$'s sweeps out a divisor $D$ in $V$, which is the unique member of $|H_V-2F_V|$, and maps to a conic on $B_3$.
    Hence $R'$ is of $E_1$-type, and the contraction morphism is the blowing-up of $B_3$ along a smooth conic.
    Following the reverse process from a cubic $3$-fold $B_3$ with a conic, we have $\rho(V)=2$ as similar to (\ref{subsec:df3}.3).

    The anti-canonical model $\ol{V}$ is embedded in $C(\P[0^3,1])\subset\P^9$, and is the intersection $\ol{W_1}\cap\ol{W_2}$ of the divisors $\ol{W_i}\subset C(\P[0^3,1])$ through the vertex.
    Here, $\ol{W_1}+P_1+P_2$, $\ol{W_2}+P_3$ are quadric hypersurface sections of $C(\P[0^3,1])$ for the cone $P_i=C(\P^3)\subset C(\P[0^3,1])$ over a fiber of $\P[0^3,1]\to\P^1$.
    The singularity of $\ol{V}$ is one ODP corresponding to the vertex of $C(\P[0^3,1])$.
    This is (\ref{thm:df4}.4).
\cons{9}
    We have $V=W_1\cap W_2\subset X=\P[0,1^4]$, and $W_1$, $W_2\in |2H-F|$.
    For general members $W_1, W_2\in|2H-F|$, $V$ is a smooth weak Fano $3$-fold with $-K_V\sim H_V$.
    The unique minimal section $s_0=\P[0] \subset X$ is contained in $V$, and is only one $K$-trivial curve of $V$.
    After flopping along $s_0$, we obtain again a weak Fano $3$-fold $V'$ with an extremal ray $R'$ of slope $\geq1$ by (\ref{eqn:slope}).

    Now we find a curve generating $R'$.
    For any sub-$\P^2$-bundle $T_\gamma=\P[0,1^2]$ for $\gamma\in Gr(\P^1\!\subset\!\P^3)$, the intersection curve $C_\gamma=V\cap T_\gamma$ is decomposed into $s_0+C_{\gamma,0}$, $C_{\gamma,0}\equiv 3[s_0]+4[l]$.
    If $C_{\gamma,0}$ is irreducible then it does not meet $s_0$.
    Consider the case that $C_{\gamma,0}$ is reducible and decomposed into $l_0+D_\gamma$, $D_\gamma\equiv 3[s_0]+3[l]$.
    By Lemma(\ref{thm:slope}), $D_\gamma$ is contained in some sub-$\P^1$-bundle $S_\lambda=\P[0,1]$.
    Hence we have only to find parameters $\lambda\in\P^3$ where $S_\lambda\cap V$ are curves $s_0+D_\lambda$.
    These parameters form a curve $A$ which is a complete intersection of two cubic surfaces in $\P^3$.
    The curve $A$ is actually decomposed into a line $Z$ and a curve $C$.
    Here $D_\lambda\equiv[s_0]$ for $\lambda\in Z$ and $D_\lambda\equiv[s_0]+[l]$ for $\lambda\in C$.
    Hence the latter curves are desired ones, and $\mu(R')=1$.
    Therefore, the contraction morphism is defined by $|H_{V'}-F_{V'}|$ with $\dim|H_{V'}-F_{V'}|=3$, $(H_{V'}-F_{V'})^3=1$.
Since the curve $C$ is of degree $8$ and of genus $7$ in $\P^3$, the morphism is a blowing-up of $\P^3$ along ${}_7C_8$.
    We have $\rho(V)=2$ by (\ref{thm:0k1n}).

    The anti-canonical model $\ol{V}$ is embedded in the cone $C(\P^3\times\P^1)\subset\P^8$ over $\P^3\times\P^1$, and is the intersection $\ol{W_1}\cap\ol{W_2}$ of two divisors $\ol{W_i}$ through the vertex of the cone.
    Here $\ol{W_i}+P$ is a quadric hypersurface section of $C(\P^3\times\P^1)$ such that $P=C(\P^3)\subset C(\P^3\times\P^1)$ is the cone over a fiber of $\P^3\times\P^1\to\P^1$.
    The singularity of $\ol{V}$ is one ODP corresponding to the vertex of $C(\P^3\times\P^1)$.
    This is (\ref{thm:df4}.5).
\cons{10}
    In this case, we have $X=\P[0^2,1^3]$, and $V=W_1\cap W_2$ for $W_1\in|2H-F|$ and $W_2\in|2H|$.
    Since $V$ is smooth for general members, $V$ is a weak Fano $3$-fold with $-K_V\sim H_F$, and has two $K$-trivial curves $s_0, s_1$ as the minimal sections of $V$.
    Thus there is the $(-F_V)$-flop $V'$.
    Now consider $T_\alpha=\P[0^2,1]$ and $C_\alpha=V\cap T_\alpha$ parameterized by $\alpha\in\P^2$.
    Since $T_\alpha$ contains $s_0$ and $s_1$, the curve $C_\alpha$ is decomposed into $s_0+s_1+D_\alpha$, $D_\alpha\equiv 2[s_0]+2[l]$.
    For general $\alpha\in\P^2$, $D_\alpha$ is irreducible and these three curves are mutually disjoint.
    The strict transform of $D_\alpha$ generates the extremal ray $R'\subset\ol{NE}(V')$, because (\ref{eqn:slope}) shows $\mu(R')\geq 1$.
    The contraction morphism of $R'$ is defined by $|H_{V'}-F_{V'}|$ with $\dim|H_{V'}-F_{V'}|=2$, and is a conic bundle.
    The discriminant locus $\Delta$ corresponds to the locus of $\alpha\in\P^2$ such that $D_\alpha$ is reducible, and $\Delta$ is a curve of degree $6$.
    Similarly to the proof of (\ref{thm:0k1n}), we have $\rho(V)=\rho(V')=2$.

    The anti-canonical model $\ol{V}$ is $Gr(5,2)\cap H_1\cap H_2\cap Q$, where $Gr(5,2)\cap H_1$ has a singular locus $\Sigma=\P^2$, and $H_2$, $Q$ are a general hyperplane and a general quadric hypersurface.
    The singularities of $\ol{V}$ are two ODP's, the intersection $\Sigma\cap H_2\cap Q$.
    This is (\ref{thm:df4}.6).
\cons{12}
    In this case, $X=\P[0^3,1^2]$, and $V=W_1\cap W_2$ for $W_1, W_2\in|2H|$.
    For general members $W_1, W_2\in|2H|$, $V$ is smooth, hence is a weak Fano $3$-fold with $-K_V\sim H_V$, and $\rho(V)=2$ by (\ref{thm:0n11}).
    There are four minimal sections $s_1,\dots,s_4\subset V$ as the $K$-trivial curves, and the $(-F_V)$-flop $V'$ of $V$ is also a weak Fano $3$-fold with the unique extremal ray $R'$.
    Since $Bs|H_V-F_V|=\bigcup_{i=1}^4s_i$, $\mu(R')\geq1$ by (\ref{eqn:slope}) as similar to other cases.
    For general sub-$\P^2$-bundle $T=\P[0^2,1]$, the intersection curve $C=T\cap V$ is of slope $1$, and generates $R'$.
    Thus the contraction morphism of $R'$ is defined by $|H_{V'}-F_{V'}|$ with $\dim|H_{V'}-F_{V'}|=1$.
    Because $\bigl(-K_{V'}|_{L'}\bigr)^2=4$ for $L'\in|H_V'-F_V|$, the morphism is a del Pezzo fibration of degree $4$.

    The anti-canonical model $\ol{V}$ is the complete intersection $Q_0\cap Q_1\cap Q_2$ of three quadric $5$-folds in $\P^6$, where $Q_0$ is defined by a quadric form of rank $4$, i.e., has $2$-dimensional singular locus $\Sigma=\P^2$, and $Q_1$, $Q_2$ are general.
    The singularities of $\ol{V}$ consist of four ODP's, the intersection $\Sigma\cap Q_1\cap Q_2$.
    This is (\ref{thm:df4}.7).
\cons{14}
    In this case, $X=\P[0^4,1]$ and $V=W_1\cap W_2$ for $W_1\in|2H|$, $W_2\in|2H+F|$.
    For general members $W_1\in|2H|$ and $W_2\in|2H+F|$, $V$ is smooth, hence is a weak Fano $3$-fold with $-K_V\sim H_V$, and $\rho(V)=2$ by (\ref{thm:0n1}) and (\ref{thm:ampleLeff}).
    There are eight minimal sections $s_1,\dots,s_8\subset V$, and the $(-F_V)$-flop $V'$ of $V$ is also a weak Fano $3$-fold with only one extremal ray $R'$.
    Members of $|2H_V-F_V|$ comes from reducible ones of $|2H_V|$ which has an irreducible component linearly equivalent to $F_V$.
    We can see that $Bs|2H_V-F_V|=\bigcup_{i=1}^8s_i$ and $\mu(R')\geq1/2$ by (\ref{eqn:slope}).
    The curve $C=V\cap T$ for general sub-$\P^2$-bundle $T=\P\bigl[0^3]$ is numerically equivalent to $4[s_0]+2[l]$; $\mu(C)=1/2$; hence $\mu(R')=1/2$.
    This shows that the strict transform $E'$ of the unique member $E\in|H_V-F_V|$ is the exceptional locus of the contraction morphism of $R'$.
    Moreover we can check that $E'$ is a quadric surface and the ray is of type $E_3$ (or $E_4$).

    The anti-canonical model $\ol{V}$ is the complete intersection $\ol{W_1}\cap \ol{W_2}\subset\P^5$ of a smooth quadric $\ol{W_1}$ and a singular cubic $4$-fold $\ol{W_2}$ which contains $\P^3$ and has eight singular points on the $\P^3$.
    The singularities of $\ol{V}$ come from the singularities of $\ol{W_2}$.
    This is (\ref{thm:df4}.8).
\cons{15}
    In this case, $X=\P[0^5] \cong\P^4\times\P^1$ and $V=W_1\cap W_2$ for $W_1\in|2H|$, $W_2\in|2H+2F|$.
    General $V$ is a smooth weak Fano $3$-fold with $-K_V\sim H_V$, and has sixteen $K$-trivial curves $s_1,\dots,s_{16}$.
    There is the $(-F_V)$-flop $V'$, a weak Fano $3$-fold with one extremal ray.
    Similarly as (\ref{subsec:df4}.3), we have $\rho(V)=2$.
    By using the homogeneous coordinates, the equations defining $V$ are $q_3(x)$ and $q_0(x)t_0^2+q_1(x)t_0t_1+q_2(x)t_1^2$ for $4$ quadrics $q_i$.
    The flopping curves correspond to $\bigl\{ x\in\P^3 \bigm| Q_i(x)=0, i=0,1,2,3 \bigr\}$.
    The divisor defined by $q_0(x)\alpha_0^2+q_1(x)\alpha_0\alpha_1+q_2(x)\alpha_1^2$ on $V$ for any $\alpha=[\alpha_0:\alpha_1]\in\P^1$ is decomposed into $D_\alpha+F_\alpha$, $D_\alpha\in|2H_V-F_V|$ and $F_\alpha\in|F_V|$.
    Since $\dim|2H_V-F_V|=1$, there is $\alpha\in\P^1$ such that $D=D_\alpha$ for each member $D\in|2H_V-F_V|$.
    Therefore, $Bs|2H_V-F_V|=\bigcap_{\alpha\in\p^1} D_\alpha=\bigl\{ x\in\P^3 \bigm| q_i(x)=0, i=0,1,2,3 \bigr\}=\bigcup_{i=1}^{16}s_i$, and the extremal ray $R'\subset\ol{NE}(V')$ is of slope $\mu(R')\geq1/2$ by (\ref{eqn:slope}).
    We consider a family of sub-$\P^2$-bundles $T_\gamma=\P[0^3]$ containing $s_0$ and $s_1$, and of intersection curves $C_\gamma=V\cap T_\gamma$, parameterized by $\gamma\in\P^2$.
    The curve $C_\gamma-s_0-s_1$ is effective and irreducible for general $\gamma\in\P^2$, and its strict transform is of slope $1/2$, whence $\mu(R')=1/2$.
    Consequently, the contraction morphism of $R'$ is defined by $|2H_{V'}-F_{V'}|$ and is a del Pezzo fibration of degree $4$.
    Moreover, $\ol{V}$ is the double covering of the quadric $3$-fold $\bigl\{ x\in\P^4 \bigm| q_3(x)=0 \bigr\}$ branched along $\bigl\{ x\in\P^4 \bigm| q_3(x)=q_0(x)q_2(x)-q_1(x)^2=0 \bigr\}$, and the singularities of $\ol{V}$ come from the singularities $\bigl\{ x\in\P^4 \bigm| q_0(x)=q_1(x)=q_2(x)=q_3(x)=0 \bigr\}$ of the branching surface, which is (\ref{thm:df4}.9).
\cons{16}
    In this case, $X=\P[0^5] \cong\P^4\times\P^1$ and $V=W_1\cap W_2$ for $W_1, W_2\in|2H+F|$.
    This is very similar to (\ref{subsec:df4}.9).
    We can find the curve of slope $1/2$, and show that the extremal ray $R'$ is of slope $1/2$ and that the contraction morphism is a del Pezzo fibration of degree $4$.
    The anti-canonical model $\ol{V}$ is the quartic $3$-fold in $\P^4$ defined by $q_0q_3-q_1q_2=0$ for $4$ quadrics $q_i$, and its singularities are $16$ points defined by $q_0=q_1=q_2=q_3=0$.
    This is (\ref{thm:df4}.10).
\cons{17}
    This case is $V=W_1\cap W_2\subset X=\P[0,1^4]$, where $W_1\in|2H-F|$ and $W_2\in|2H|$.
    General $V$ is smooth and $-K_V\sim H_V-F_V$.
    The divisor $W_2$ in $X$ does not meet the minimal section $s_0=\P[0]\subset X$, hence so is $V$.
    Since $\rho(W_2)=2$ by (\ref{thm:01n}) and since $V$ is ample in $W_2$, we have $\rho(V)=2$ by (\ref{thm:ampleLeff}).
    Denote by $T_\lambda$ the ruled surface $\P[0,1]\subset X$ parameterized by $\lambda\in\P^3$.
    For a curve $C$ in $V$, if $\mu(C)\leq1$ then $C\subset T_\lambda$ for some $\lambda\in\P^3$ by (\ref{thm:slope}).
    To find $K$-trivial curves, we have only to check the cases that $V\cap T_\lambda$ is a curve.
    The bihomogeneous coordinates analysis shows that there are four bisections $C_1,\dots,C_4$ and $32$ sections $s_1,\dots,s_{32}$ as $K$-trivial curves.
    There is the $(-F_V)$-flop $V'$ of $V$, which is a weak Fano $3$-fold.

    We fix the unique extremal ray $R'$ of $\ol{NE}(V')$.
    Consider the linear system $|4H_V-5F_V|$.
    We can see $Bs|4H_V-5F_V|=\bigcup_{i=1}^4C_i\cup\bigcup_{j=1}^{32}s_j$, and $\mu(C')\geq 5/4$ for any curve $C'\in V'$ by (\ref{eqn:slope}).
    For any $L_\alpha=\P[0,1^2]$ containing $C_1$ parameterized by $\alpha\in\P^2$, each curve $D_\alpha=L_\alpha\cap V$ is decomposed into $D_{\alpha,0}+C_1$, $D_{\alpha,0}\equiv 2[s_0]+4[l]$.
    For some $\alpha\in\P^2$, $D_{\alpha,0}$ is decomposed into $D_{\alpha,1}+l_0$ such that $D_{\alpha,1}\equiv 2[s_0]+3[l]$.
    Hence the strict transform $D_{\alpha,1}'$ of $D_{\alpha,1}$ has a slope $\mu(D_{\alpha,1}')=5/4$, and the extremal ray has a slope $\mu(R')=5/4$.
    We can see that the type of $R'$ is $D_1$ and that the contraction morphism is a del Pezzo fibration of degree $4$ over $\P^1$.

    The anti-canonical model $\ol{V}$ is the double covering of $\P^3$ branched along a sextic surface $S$ defined by $(q_1^2-4q_0q_2)y_1^2 - 4(p_0^2q_2-p_0p_1q_1+p_1^2q_0)$ for quadrics $p_i=p_i(y)$ and $q_i=q_i(y)$ in $y_1, \dots, y_4$.
    These quadrics appear in the defining equations $x_0y_1+p_0t_0+p_1t_1 = x_0^2+q_0t_0^2+q_1t_0t_1+q_2t_1^2 = 0$ of $V$ by the bihomogeneous coordinates on $X$.
    The singularities of $\ol{V}$ are $36$ ODP's coming from the singularities of $S$, $\{ (y_i)\in S | y_1=p_0=p_1=0 \} \cup \{ (y_i)\in S | q_0y_1^2-p_0^2 = q_1y_1^2-2p_0p_1 = q_2y_1^2-p_1^2 = 0 \}$. 
    This obtains (\ref{thm:df4}.11).

\section{Del Pezzo fibrations of degree 5} \label{sec:df5}
\newcount\wm \newcount\wmt
\def\Putij(#1,#2){
  \wm=-\wmt
  \advance\wm #1 \advance\wm #2
  \wm=-\wm \number\wm
}
\def\WM(#1,#2,#3,#4,#5){
  \wmt=#1
  \advance\wmt #2 \advance\wmt #3
  \advance\wmt #4 \advance\wmt #5
  \divide \wmt by 2
  \left[\begin{array}{cccc}
    \Putij(#1,#2) & \Putij(#1,#3) & \Putij(#1,#4) & \Putij(#1,#5) \\
    & \Putij(#2,#3) & \Putij(#2,#4) & \Putij(#2,#5) \\
    & & \Putij(#3,#4) & \Putij(#3,#5) \\
    & & & \Putij(#4,#5)
  \end{array}\right]
}
\def\dual#1{\widehat#1}

\newpar\Label{subsec:setup5}
    In this section, we treat del Pezzo fibratoins of degree $5$ over $\P^1$, and derive Theorem(\ref{thm:df5}) and Supplement(\ref{thm:df5s}).
    Let $V$ be a del Pezzo fibration of degree $5$ over $\P^1$ which is a subvariety of $X = \P[0,a_1,a_2,a_3,a_4,a_5]$, $0=a_0\leq a_1\leq a_2\leq a_3\leq a_4\leq a_5$ (see (\ref{subsec:setup})):
$$\begin{array}{ccl}
  V & \Longarrow{\psi}{} & X = \P[0,a_1,a_2,a_3,a_4,a_5] \vspace{5pt}\\
  \Downarrow{\varphi} & \Swarrow{\pi} & \qquad M \vspace{5pt}\\
  \P^1
\end{array}$$
    Here $V$ is defined by the Pfaffian of the $4\times4$ diagonal minors of a skew-symmetric matrix
$$M = \left(\begin{array}{cccccc}
    0   &  m_{12} &  m_{13} &  m_{14} & m_{15} \\
-m_{12} &     0   &  m_{23} &  m_{24} & m_{25} \\
-m_{13} & -m_{23} &     0   &  m_{34} & m_{35} \\
-m_{14} & -m_{24} & -m_{34} &     0   & m_{45} \\
-m_{15} & -m_{25} & -m_{35} & -m_{45} &    0
\end{array}\right)$$
whose elements $m_{ij}$ are members of $H^0(X, \O(H + w_{ij}F))$ for some integers $w_{ij}$.
    Namely $V$ is defined by the following five equations in $X$ :
\begin{eqnarray*}
f_5 = m_{12}m_{34} - m_{13}m_{24} + m_{14}m_{23} &\in& H^0(X, \O(2H + k_5F)) \\
f_4 = m_{12}m_{35} - m_{13}m_{25} + m_{15}m_{23} &\in& H^0(X, \O(2H + k_4F)) \\
f_3 = m_{12}m_{45} - m_{14}m_{25} + m_{15}m_{24} &\in& H^0(X, \O(2H + k_3F)) \\
f_2 = m_{13}m_{45} - m_{14}m_{35} + m_{15}m_{34} &\in& H^0(X, \O(2H + k_2F)) \\
f_1 = m_{23}m_{45} - m_{24}m_{35} + m_{25}m_{34} &\in& H^0(X, \O(2H + k_1F))
\end{eqnarray*}
    Here we may assume that
\Eq(seq5){
k_1 \leq k_2 \leq k_3 \leq k_4 \leq k_5}\\
without loss of generality, and we can see that
\Eq(wij){
w_{ij} = t - k_i - k_j}\\
for $t = (k_1 + k_2 + k_3 + k_4 + k_5)/2$.
    In addition, assume that $V$ is a weak Fano 3-fold with only finite $K$-trivial curves and with $\rho(V)=2$.
    For convenience, $w(M)$ denotes the matrix of twist of each element of $M$:
$$w(M) = \left[\begin{array}{rrrrr}
  w_{12} & w_{13} & w_{14} & w_{15} \\
  & w_{23} & w_{24} & w_{25} \\
  & & w_{34} & w_{35} \\
  & & & w_{45}
\end{array}\right].$$

\newpar
    In the Chow ring ${\rm CH}(X)$, $V$ is equivalent to $5H^3 + 3tH^2F$.
    The anti-canonical divisor $V$ is
$$
  -K_V\sim H_V + (2 - \sum_{i=0}^5 a_i - t)F_V
$$
    Enumerating the intersection numbers, we have
\Eq(nef5){
  10 - 5(a_4+a_5) - 2t \geq 0,}
\Eq(eff5){
  5(a_1+a_2+a_3) + 3t \geq 0, \qquad {\rm and}}
\Eq(big5){
  30 - 12t - 10\sum_{i=0}^5 a_i > 0.}\\
from (\ref{eqn:nef}), (\ref{eqn:eff}) and (\ref{eqn:big}), respectively.

\newpar
    Sevral additional inequalities listed below are obtained from smoothness of $V$ and from $\Pic V \cong \Z^{\op2}$.
\Eq(w12){    a_0 + w_{12} \geq 0,}
\Eq(a5w35){  a_5 + w_{35} \geq 0,}
\Eq(a2w14){  a_2 + w_{14} \geq 0,}
\Eq(a4w25){  a_4 + w_{25} \geq 0,}
\Eq(a4w34){  a_4 + w_{34} \geq 0,}
\Eq(or){  a_0 + w_{14} \geq 0 \quad {\rm or} \quad a_5 + w_{45} \geq 0,}
\Eq(a1w14){ \phantom{\qquad {\rm and}} a_1 + w_{14} \geq 0, \qquad {\rm and}}
\Eq(a1w23){  a_1 + w_{23} \geq 0.}\\
    We will show these inequalities.

    The contrary assertion of each (\ref{eqn:w12})--(\ref{eqn:or}) implies that $V$ has singularities, and the contrary assertion of (\ref{eqn:a1w14}) or of (\ref{eqn:a1w23}) contradicts $\rho(V) = 2$.

\def\pf#1{\vspace{3pt}\noindent{\it Proof of} #1.}
    \pf{(\ref{eqn:w12})}
    If $a_0 + w_{12} < 0$, then we have $a_0 + w_{ij} < 0$ and $H^0(\P^1, \O(a_0 + w_{ij}) = 0$ for $1 \leq i < j \leq 5$.
    It follows that each $m_{ij} \in H^0(X, \O(H + w_{ij}F)) = \bigoplus_{l = 0}^{5} H^0(\P^1, \O(a_l + w_{ij}))$ vanishes along the minimal section $s_0=\P[0]$ of $X \to \P^1$, and $V$ has singularities along $s_0$.
    This contradict smoothness of $V$.

    \pf{(\ref{eqn:a5w35})}
    Assuming the contrary, we have $a_i + w_{35} < 0$ and $a_i + w_{45} < 0$ for $i = 0, \dots, 5$. It follows that $m_{35} = m_{45} = 0$, hence $f_1$ is reducible, which implies that $V$ is reducible.

    \pf{(\ref{eqn:a2w14})}
    If $a_2+ w_{14} < 0$, then we obtain that $H^0(\P^1, (\O\op\O(a_1)\op\O(a_2)) \otimes \O(w_{ij})) = 0$ for $j \geq 4$.
    It follows that $m_{ij}$ vanishes on $\P[0,a_1,a_2]$ for $j \geq 4$, and so does $f_i$ for $i = 1, 2, \dots, 5$.
    Therefore $V$ contains the sub-$\P^2$-bundle $\P[0,a_1,a_2]$ of $X$,  and $V$ is reducible, which is a contradiction.

    \pf{(\ref{eqn:a4w25})}
    Assume the contrary. It implies that $H^0(\P^1, (\bigoplus_{l=0}^{4} \O(a_l)) \otimes \O(w_{i5})) = 0$ for $i \geq 2$.
    Hence $m_{ij}$ vanishes on $\P[0,a_1,a_2,a_3,a_4]$ for $i \geq 2$, and so does $f_1$.
    It follows that $f_1$ is reducible and so is $V$.
    This is a contradiction.

    \pf{(\ref{eqn:a4w34})}
    Similarly to (\ref{eqn:a4w25}), assuming the contrary, we have that $f_1$ is reducible, which is a contradiction.

    \pf{(\ref{eqn:or})}
    We will derive a contradiction from the assumption that $a_0 + w_{14} < 0$ and $a_5 + w_{45} < 0$.
    The inequality $a_0 + w_{14} < 0$ gives $a_0 + w_{ij} < 0$ for $j = 4$ or $5$, hence $m_{ij}$ for $j =4$ or $5$ vanishes on the minimal section $s_0 = \P[0] \subset X = \P[0,a_1,a_2,a_3,a_4,a_5]$, and each $f_i$ vanishes on $s_0$.
    Moreover, $f_i$, for $i = 1$, $2$, and $3$, vanishes on $s_0$ with multiplicity $2$, becouse $a_5 + w_{45} < 0$ gives $m_{45} = 0$.
    Thus $V$ has singularities alog $s_0$, which is a contradiction.

    \pf{(\ref{eqn:a1w14}) {\it and} (\ref{eqn:a1w23})}
    Since $\Pic V \cong \Z^{\op2}$, $\Pic V$ is generated by the restrictions $H_V = H|_V$ and $F_V = F|_V$ of $H$ and $F$ to $V$.
    As a $2$-cycle in $X = \P[a_0,\dots,a_5]$, we have
\begin{eqnarray*}
H_V = &\hspace{-7pt}H \cdot (5H^3+3tH^2F) &\hspace{-7pt}= 5H^4 + 3tH^3F, \qquad {\rm and} \\
F_V = &\hspace{-7pt}F \cdot (5H^3+3tH^2F) &\hspace{-7pt}= 5H^3F,
\end{eqnarray*}
hence any member of $\Pic V$ is equivalent to $5aH^4 + (3ta + 5b)H^3F$ for some integers $a$ and $b$.

    If $a_1 + w_{14} < 0$, then we obtain that $H^0(\P^1, \O(a_1 + w_{i4})) = H^0(\P^1, \O(a_1 + w_{i5})) = 0$ for $i = 1$, $2$, and $3$.
    Hence $f_i$ vanishes on the ruled surface $S = \P[0,a_1]$ in $X$, and $S$ is contained in $V$ as a divisor.
    As a $2$-cycle in $X$, $S$ is equivalent to $(H-a_5F) \cdot (H-a_4F) \cdot (H-a_3F) \cdot (H-a_2F) = H^4 - (a_2 + a_3 + a_4 + a_5)H^3F$, but this cannot be expressed by a linear combination of the generators $H_V$ and $F_V$ of $\Pic V$.
    This is a contradiction, and (\ref{eqn:a1w14}) is obtained.

    If $a_1 + w_{23} < 0$, then we have $H^0(\P^1, \O(a_1 + w_{ij})) = 0$ for $2 \leq i < j \leq 5$.
    Each $f_i$ vanishes on the ruled surface $S = \P[0,a_1]$ in $V$.
    The argument similar to the above gives the inequality (\ref{eqn:a1w23}).

\newpar\Label{subsec:poss5}
    We will now enumerate the possibilities of values $a_i$ and $k_i$ under the conditions (\ref{eqn:seq5})--(\ref{eqn:a1w23}).
    We have $2t \geq 5k_1$ from the sequence (\ref{eqn:seq5}), hence we have $2 - k_1 \geq a_4 + a_5$ from (\ref{eqn:nef5}).
    Conditions (\ref{eqn:wij}), (\ref{eqn:a4w25}), and (\ref{eqn:a4w34}) imply $2a_4 \geq -k_1$.
    Gathering the above two inequalities and $a_4 \leq a_5$, we can see that there exist only three cases: $a_4+a_5 = 2-k_1$, $1-k_1$, and $-k_1$.
    Thus we treat these cases separetely, and give three Tables \ref{tbl:t1}--\ref{tbl:t3} of possibilities.

    First, we set $a_4+a_5 = 2-k_1$, then we have $2t = 5k_1$ and $k_i = k_1$ for all $i$, hence $w_{ij} = k_1$ for $1 \leq i < j \leq 5$.
   The condition (\ref{eqn:w12}) gives $k_1 \geq 0$, and the condition (\ref{eqn:big5}) fixes the values of $a_i$ and $k_i$ as in Table \ref{tbl:t1}.
\begin{table}[htbp]
\caption{}\Label{tbl:t1}
$$\begin{array}{c@{\qquad}ccccc@{\qquad}rrrrr}
{\rm Nos.} & a_1 & a_2 & a_3 & a_4 & a_5 & k_1 & k_2 & k_3 & k_4 & k_5 \\
1 & 0 & 0 & 0 & 0 & 2 &  0 &  0 &  0 &  0 &  0 \\
2 & 0 & 0 & 0 & 1 & 1 &  0 &  0 &  0 &  0 &  0
\end{array}$$
\end{table}

    Next, we treat the cases $a_4+a_5 = 1-k_1$, and $-k_1$.
    There are integers $0 \leq d_2 \leq d_3 \leq d_4 \leq d_5$ which denote the difference $d_i = k_i - k_1$ for $i = 2, \dots, 5$.
    Then (\ref{eqn:nef5}), (\ref{eqn:a5w35}), (\ref{eqn:a4w25}), and (\ref{eqn:a4w34}) imply
\Eq(nef51){
5 \geq\: d_2+d_3+d_4+d_5, \qquad \mbox{in the case} \quad a_4+a_5 = 1-k_1,}
\Eq(nef50){
10 \geq\: d_2+d_3+d_4+d_5, \qquad \mbox{in the case} \quad a_4+a_5 = -k_1,}
\Eq(a5){
2a_5+k_1 \geq (d_3-d_2) + (d_5-d_4),}
\Eq(a43){ \phantom{\quad {\rm and}}
2a_4+k_1 \geq (d_3-d_2) - (d_5-d_4), \quad {\rm and}}
\Eq(a45){
2a_4+k_1 \geq (d_5-d_4) - (d_3-d_2),}\\
respectively.

    In the case $a_4+a_5 = -k_1$, we have that $a_4=a_5=-k_1/2$ and that $k_1$ is an even integer less than or equal to $0$.
    Hence (\ref{eqn:a5}) shows that $d_2=d_3$ and $d_4=d_5$, and then (\ref{eqn:nef50}), (\ref{eqn:w12}) and (\ref{eqn:or}) gives that $k_1 = 0$, $-2$, or $-4$.
    If $k_1 = -4$, then we obtain that the differences $d_i$ between $k_i$ and $k_1$ equal $2$ for all $i$ or that $d_2=d_3=2$ and $d_4=d_5=3$.
    If $k_1 = -2$, then we obtain that $d_2=d_3=1 \leq d_4=d_5 \leq 4$, that $d_i=2$ for all $i$, or that $d_2=d_3=2$ and $d_4=d_5=3$.
    If $k_1 =0$, then $a_i = 0$ for all $i$.
    Thus the inequalities (\ref{eqn:a1w23}), (\ref{eqn:eff5}), and (\ref{eqn:big5}) fix the possible values of $a_i$ and $k_i$ as in Table \ref{tbl:t2}.
\begin{table}[htbp]
\caption{}\Label{tbl:t2}
$$\begin{array}{c@{\qquad}ccccc@{\qquad}rrrrr}
{\rm Nos.} & a_1 & a_2 & a_3 & a_4 & a_5 & k_1 & k_2 & k_3 & k_4 & k_5 \\
 3 & 2 & 2 & 2 & 2 & 2 & -4 & -2 & -2 & -2 & -2 \\
 4 & 1 & 1 & 2 & 2 & 2 & -4 & -2 & -2 & -1 & -1 \\
 5 & 1 & 1 & 1 & 2 & 2 & -4 & -2 & -2 & -1 & -1 \\
 6 & 1 & 1 & 1 & 1 & 1 & -2 & -1 & -1 & -1 & -1 \\
 7 & 1 & 1 & 1 & 1 & 1 & -2 & -1 & -1 &  0 &  0 \\
 8 & 0 & 1 & 1 & 1 & 1 & -2 & -1 & -1 &  0 &  0 \\
 9 & 0 & 1 & 1 & 1 & 1 & -2 & -1 & -1 &  1 &  1 \\
10 & 0 & 0 & 1 & 1 & 1 & -2 & -1 & -1 &  1 &  1 \\
11 & 0 & 0 & 0 & 1 & 1 & -2 & -1 & -1 &  2 &  2 \\
12 & 0 & 0 & 0 & 1 & 1 & -2 &  0 &  0 &  1 &  1 \\
13 & 0 & 0 & 0 & 0 & 0 &  0 &  0 &  0 &  0 &  0 \\
14 & 0 & 0 & 0 & 0 & 0 &  0 &  0 &  0 &  1 &  1 \\
15 & 0 & 0 & 0 & 0 & 0 &  0 &  0 &  0 &  2 &  2 \\
16 & 0 & 0 & 0 & 0 & 0 &  0 &  1 &  1 &  1 &  1
\end{array}$$
\end{table}

    Last, we set $a_4+a_5 = 1-k_1$.
    For even $k_1$, we obtain that $a_4 = -k_1/2$ and $a_5 = -k_1/2 + 1$ and that $d_2+d_3+d_4+d_5$ is even, while, for odd $k_1$, we obtain that $a_4 = a_5 = -(k_1+1)/2$ and that $d_2+d_3+d_4+d_5$ is odd.
    Hence the inequalities (\ref{eqn:nef51}), (\ref{eqn:a5}), (\ref{eqn:a43}), and (\ref{eqn:a45}) imply that $(d_2,d_3,d_4,d_5) = (0,0,0,0)$, $(0,0,1,1)$, $(0,0,2,2)$, $(1,1,1,1)$, or $(0,1,1,2)$ in case of even $k_1$, and that $(d_2,d_3,d_4,d_5) = (0,0,0,1)$, $(0,0,1,2)$, $(0,0,2,3)$, $(1,1,1,2)$, $(0,1,1,1)$, or $(0,1,2,2)$ in case of odd $k_1$.
    Furthermore (\ref{eqn:w12}) shows that $-1 \leq k_1 \leq 5$.
    However, it follows that $k_1 \ne -1$ from (\ref{eqn:big5}), and that $k_1 \ne 5$ from (\ref{eqn:or}).
    From other values of $k_1$, the other inequalities togethar fix the possible values of $a_i$ and $k_i$ as in Table \ref{tbl:t3}.
\begin{table}[htbp]
\caption{}\Label{tbl:t3}
$$\begin{array}{c@{\qquad}ccccc@{\qquad}rrrrr}
{\rm Nos.} & a_1 & a_2 & a_3 & a_4 & a_5 & k_1 & k_2 & k_3 & k_4 & k_5 \\
17 & 2 & 2 & 2 & 2 & 3 & -4 & -4 & -3 & -3 & -2 \\
18 & 1 & 2 & 2 & 2 & 3 & -4 & -4 & -3 & -3 & -2 \\
19 & 2 & 2 & 2 & 2 & 2 & -3 & -3 & -2 & -2 & -2 \\
20 & 1 & 2 & 2 & 2 & 2 & -3 & -3 & -2 & -2 & -2 \\
21 & 1 & 1 & 2 & 2 & 2 & -3 & -3 & -2 & -2 & -2 \\
22 & 1 & 1 & 2 & 2 & 2 & -3 & -2 & -2 & -2 & -1 \\
23 & 1 & 1 & 1 & 2 & 2 & -3 & -2 & -2 & -2 & -1 \\
24 & 1 & 1 & 1 & 1 & 2 & -2 & -2 & -2 & -1 & -1 \\
25 & 1 & 1 & 1 & 1 & 2 & -2 & -2 & -1 & -1 &  0 \\
26 & 1 & 1 & 1 & 1 & 2 & -2 & -1 & -1 & -1 & -1 \\
27 & 0 & 1 & 1 & 1 & 2 & -2 & -2 & -1 & -1 &  0 \\
28 & 1 & 1 & 1 & 1 & 1 & -1 & -1 & -1 & -1 &  0 \\
29 & 0 & 1 & 1 & 1 & 1 & -1 & -1 & -1 & -1 &  0 \\
30 & 0 & 1 & 1 & 1 & 1 & -1 & -1 & -1 &  0 &  1 \\
31 & 0 & 1 & 1 & 1 & 1 & -1 & -1 &  0 &  0 &  0 \\
32 & 0 & 0 & 1 & 1 & 1 & -1 & -1 & -1 &  0 &  1 \\
33 & 0 & 0 & 1 & 1 & 1 & -1 & -1 &  0 &  0 &  0 \\
34 & 0 & 0 & 0 & 1 & 1 & -1 & -1 & -1 &  1 &  2 \\
35 & 0 & 0 & 0 & 1 & 1 & -1 & -1 &  0 &  1 &  1 \\
36 & 0 & 0 & 0 & 1 & 1 & -1 &  0 &  0 &  0 &  1 \\
37 & 0 & 0 & 0 & 0 & 1 &  0 &  0 &  0 &  0 &  0 \\
38 & 0 & 0 & 0 & 0 & 1 &  0 &  0 &  0 &  1 &  1
\end{array}$$
\end{table}

\newpar\Label{subsec:exc5}
    Some of the values in Tables \ref{tbl:t1}--\ref{tbl:t3} cannot give weak Fano 3-folds $V$ with $\Pic V \cong \Z^{\op2}$ by the construction stated in (\ref{subsec:setup5}). We will now excluded these ones.
\ex{1}
    Since the sub-$\P^4$-bundle $Y = \P[0^5] \subset X = \P[0^5,2]$ is isomorphic to $\P^4 \times \P^1$, the intersection $V \cap Y$ is the product $C \times \P^1$, where $C$ is a curve in a fiber $F_V$.
    Thus, any section $s = \{ c \} \times \P^1$ is $K$-trivial curve because of $(-K_V \cdot s) = (H_V \cdot s) = 0$, hence there is a family of $K$-trivial curves in $V$, which contradicts to the assumption.
\exs{3, 4, 17, 19, 20, 22, {\it and} 26}
    In these cases, we have $w_{23} = t-k_2-k_3 = (k_1-k_2-k_3+k_4+k_5)/2 < 0$ and $2 - \sum_{i=0}^5 a_i - t < 0$.
    Hence, $w_{ij}$, $i\ne1$, are negative and $m_{ij}$, $i\ne1$, vanish on the minimal section $s_0 = \P[0] \subset X$.
    Thus, $s_0$ is contained in $V$.
    Since $(-K_V \cdot s_0) = (H_V +(2 - \sum_{i=0}^5 a_i - t)F_V \cdot s_0) < 0$, $-K_V$ is not nef.
    This contradicts the assumption.
\exs{5, 7, 8, 9, 10, 11, 14, 15, 30, 32, 34, {\it and} 35}
    For example, in the case No.5, each element $m_{ij}$ of the matrix $M$ defining $V$ is a member of $H^0(X, \O(H + w_{ij}F))$, where
$$w(M)=\WM(-4,-2,-2,-1,-1).$$
    Consider the exact sequence $\O \op \O(1)^{\op3} \op \O(2)^{\op2} \to \O(1)^{\op3} \to 0$ defined by $m_{15}, m_{25}, m_{35} \mapsto 0$.
    Then, $m_{15}$, $m_{25}$, and $m_{35}$ vanish on the sub-$\P^2$-bundle $\P[1^3] \subset X=\P[0,1^3,2^2]$ associated to the exact sequence above.
    Moreover, $m_{45}$ is zero constantly, hence the surface $S = (f_5) \cap \P[1^3]$ is contained in $V$.
    However, $S$ is equivalent to $2H^4 - 9H^3F$ which is not expressed by a linear combination of $H_V$ and $F_V$, hence $\Pic V$ is not isomorphic to $\Z^{\op2}$.
    Thus this case is to be excluded.

    Similarly, in the other cases, $m_{15}$, $m_{25}$, and $m_{35}$ vanish on $\P[a_i, a_j, a_k] \subset X$ for some $0 \leq i < j < k \leq 5$, and $V$ contains the surface $S = (f_5) \cap \P[a_i, a_j, a_k]$ equivalent to $2H^4 + (k_5 - 2(\sum_{l \ne i,j,k} a_l))H^3F$, which is a contradiction.
\exs{18, 23, 25, 27, {\it and} 36}
    For the case No.18, the twists of the elements of $M$ are
$$w(M)=\WM(-4,-4,-3,-3,-2).$$
    Since $m_{35}$, $m_{45}\in H^0(X, \O(H-3F)) \cong H^0(\P^1, \E(-3)) \cong H^0(\P^1, \O) \cong \C$, we may assume that $m_{45}=0$ by making a linear transform of $M$.
    Hence the argument similar to the case No.5 leads to a contradiction.

    For the case No.23, the twists of the elements of $M$ is
$$w(M)=\WM(-3,-2,-2,-2,-1).$$
    We can also assume that $m_{45}=0$ by linear tarnsform of $M$ if necessary, because $m_{25}$, $m_{35}$, $m_{45} \in H^0(X, \O(H-2F)) \cong H^0(\P^1, \E(-2)) \cong H^0(\P^1, \O)^{\op2} \cong \C^{\op2}$.
    In this case, we have a contradiction.

    The other cases are similar to these cases.
\ex{12}
    In this case, we have
$$w(M)=\WM(-2,0,0,1,1),$$
hence $m_{24}$, $m_{25}$, $m_{34}$, and $m_{35}$ vanish on $\P[0^4] \subset X=\P[0^4,1^2]$, and $m_{45}$ is constantly zero.
    Consider the sub-$\P^2$-bundle $\P[0^3] \subset X$ associated to the exact sequence $\E \to \O^{\op3} \to 0$ defined by $\O(1)^{\op2} \to 0$ and by $m_{23} \mapsto 0$.
    Then, each equation $f_i$, $i=1, \dots, 5$, vanishes on the sub-$\P^2$-bundle.
    The sub-$\P^2$-bundle is contained in $V$, hence $V$ is reducible, which is a contradiction.
\ex{13}
    This case is the trivial bundle, hence the Picard group $\Pic V$ is of rank $5$, which is an excluded case.

\newpar \Label{subsec:df5}
    From the above argument (\ref{subsec:exc5}), we have the following list of possible values of $a_i$ and $k_i$ as in Table \ref{tbl:df5}.
\begin{table}[htbp]
\caption{}\Label{tbl:df5}
$$\begin{array}{ccccccrrrrrlc}
{\rm Nos.} & a_1 & a_2 & a_3 & a_4 & a_5 & k_1 & k_2 & k_3 & k_4 & k_5 & -K_V & (-K_V)^3 \\
21 & 1 & 1 & 2 & 2 & 2 & -3 & -3 & -2 & -2 & -2 & H_V & 22 \\
37 & 0 & 0 & 0 & 0 & 1 &  0 &  0 &  0 &  0 &  0 & H_V+F_V & 20 \\
24 & 1 & 1 & 1 & 1 & 2 & -2 & -2 & -2 & -1 & -1 & H_V & 18 \\
 6 & 1 & 1 & 1 & 1 & 1 & -2 & -1 & -1 & -1 & -1 & H_V & 16 \\
29 & 0 & 1 & 1 & 1 & 1 & -1 & -1 & -1 & -1 &  0 & H_V & 14 \\
33 & 0 & 0 & 1 & 1 & 1 & -1 & -1 &  0 &  0 &  0 & H_V & 12 \\
 2 & 0 & 0 & 0 & 1 & 1 &  0 &  0 &  0 &  0 &  0 & H_V & 10 \\
38 & 0 & 0 & 0 & 0 & 1 &  0 &  0 &  0 &  1 &  1 & H_V &  8 \\
16 & 0 & 0 & 0 & 0 & 0 &  0 &  1 &  1 &  1 &  1 & H_V &  6 \\
28 & 1 & 1 & 1 & 1 & 1 & -1 & -1 & -1 & -1 &  0 & H_V-F_V & 4 \\
31 & 0 & 1 & 1 & 1 & 1 & -1 & -1 &  0 &  0 &  0 & H_V-F_V & 2 \\
\end{array}$$
\end{table}
The column Nos. in Table \ref{tbl:df5} denotes the numbers in Tables \ref{tbl:t1}--\ref{tbl:t3}.

    Now we construct the del Pezzo fibrations with each conditions in the rest cases on Table \ref{tbl:df5}, and show Theorem(\ref{thm:df5}) and Supplement(\ref{thm:df5s}).
\cons{21}
    In this case, we have $X=\P[0,1^2,2^3]$ and $V$ is defined by $M$ with
$$w(M)=\WM(-3,-3,-2,-2,-2).$$
    Regard $X$ as an intersection of three members of $|H_Y-F_Y|$ of the $\P^9$-bundle $Y= \P[-w_{12},-w_{13},\dots,-w_{45}] =\P[0,1^6,2^3]$ over $\P^1$.
    Denote by $G$ the subvariety of $Y$ defined by the Pfaffian of the $4\times4$ diagonal minors of a skew-symmetric matrix $M_Y$ with $m_{ij} \in H^0(Y, \O(H_Y+w_{ij}H_Y))$.
    The variety $G$ is the Grassmannian bundle whose fiber is $Gr(5,2)=Gr(\P^1\subset\P^4)$.
    For general $M_Y$, the variety $G$ is smooth.
    By Bertini's theorem, general $V$ is smooth outside the minimal section $s_0=\P[0] \subset Y$, since $V$ is the intersection of $G$ and three members of $|H_Y-F_Y|$ with $Bs|H_Y-F_Y|=s_0$.
    Taking the local coordinates along $s_0$ of $X$ and making a bit of calculation, we can see that general $V$ is smooth also along $s_0$.

    We now consider the dual variety of $G \subset Y$ over $\P^1$ in order to show that $\rho(V)=2$.
    Denote by $\dual{Y}$ the dual projective space bundle $\P[0,(-1)^6,(-2)^3]$ of $Y$ over $\P^1$, parameterizing hyperplanes in each projective space fiber.
    Denote by $\dual{G} \subset \dual{Y}$ the subvariety parameterizing hyperplanes tangent to $G$ in each fiber.
    By self-duality of the Grassmannian $Gr(5,2)$, the variety $\dual{G}$ is also defined by the Pfaffian of the $4\times4$ diagonal minors of a skew-symmetric matrix $M_{\dual{Y}}$ with $\dual{m_{ij}} \in H^0(\dual{Y}, \O(H_{\dual{Y}}-w_{ij}F_{\dual{Y}}))$.
    The variety $\dual{X} \subset \dual{Y}$ parameterizing hyperplanes containing $X$ is the sub-$\P^3$-bundle $\P[(-1)^4]$.
    The intersection $C_V = \dual{G} \cap \dual{X}$ is a curve quintuple covering over $\P^1$.
    Kanev's approach \cite{kanev} shows $\rho(V)=2$ if the curve is irreducible.
    (We explain this fact by more elementary argument in (\ref{subsec:kanev}).)

    Thus, we have only to study the curve $C_V$ to see $\rho(V)=2$.
    The curve $C_V$ is defined by the skew-symmetric matrix $M_C$ with $\dual{m_{ij}} \in H^0(\dual{X}, \O(\dual{H}-w_{ij}\dual{F}))$, where $\dual{H}$ and $\dual{F}$ are the tautological line bundle and a general fiber of $\dual{X}\to\P^1$, respectively.

    Hence $\dual{m_{12}}=0$ and $\dual{m_{1j}}, \dual{m_{2j}} \in H^0(\dual{X}, \O(\dual{H}+\dual{F})) \cong H^0(\P^1, \O^{\op4}) \cong \C^{\op4}$ for $j=3,4,5$.
    Choosing the appropriate homogeneous coordinates $[x_0,x_1,x_2,x_3]\times[t_0,t_1]$ on $\dual{X}\cong\P^3\times\P^1$ and making a linear transformation of $M_C$, we have
$$M_C=
\left(\begin{array}{ccccc}
   0 &   0 & x_0 & x_2 &  a \\
   0 &   0 & x_1 & x_3 &  b \\
-x_0 &-x_1 &   0 &   0 &  g \\
-x_2 &-x_3 &   0 &   0 &  h \\
  -a &  -b &  -g &  -h &  0
\end{array}\right),$$
where $a$ is a linear form in $x_1$ and $x_3$, $b$ is a linear form in $x_0$, $x_1$, $x_2$ and $x_3$, and $g$ and $h$ are bilinear forms in $x_i$'s and $t_i$'s.
    Hence we have
\begin{eqnarray*}
\dual{f_5} &=& x_0x_3 - x_1x_2, \\
\dual{f_4} &=& x_0 b - x_1 a, \\
\dual{f_3} &=& x_2 b - x_3 a, \\
\dual{f_2} &=& x_0 h - x_2 g, \\
\dual{f_1} &=& x_1 h - x_3 g.
\end{eqnarray*}
    Now consider the projection $p : \P^3\times\P^1 \cdots\to S=\P^1\times\P^1$ eliminating $x_2$ and $x_3$.
    Denote by $D$ the strict transform of $C_V$ by $p$.
    We will see that $D$ is irreducible, that $p : C_V \to D$ is generically bijective, and that no component of $C_V$ contracts to a point by $p$.
    Since $\dual{f_5}$ and $\dual{f_4}$ are linear in $x_2$ and $x_3$, there is a unique solution of the system of equations $\dual{f_5} = \dual{f_4} = 0$ for $x_2$ and $x_3$; there are expressions $x_2=\xi x_0/\eta$ and $x_3=\xi x_1/\eta$, by $x_0$ and $x_1$, where $\xi$ and $\eta$ are quadratic forms in $x_0$ and $x_1$.
    Substituting the above expressions for $x_2$ and $x_3$ in $\dual{f_2}$ or $\dual{f_1}$, we have the homogeneous polynomial
$$f = (h_0\eta^2+(h_2-g_0)\eta\xi-g_2\xi^2)x_0 + (h_1\eta^2+(h_3-g_1)\eta\xi-g_3\xi^2)x_1$$
of bidegree $(5,1)$ in $[x_0:x_1]$ and $[t_0:t_1]$.
    Here, $g=\sum_{i=0}^3 g_i$, $h=\sum_{i=0}^3 h_i$ for linear forms $g_i$, $h_i$ in $t_i$'s.
    The polynominal $f$ defines $D$ in $S=\P^1\times\P^1$.
    By regarding $f\in\C[\eta,\xi,x_0,x_1,t_0,t_1]$ as a polynomial in the six variables $\eta$, $\xi$, $x_0$, $x_1$, $t_0$ and $t_1$, the coefficient $h_i\eta^2+(h_{i+2}-g_i)\eta\xi-g_{i+2}\xi^2$ of $x_i$ is irreducible in $\C[\eta,\xi,t_0,t_1]$ for general $g_i$, $h_i$, hence $f$ is irreducible in $\C[\eta,\xi,x_0,x_1,t_0,t_1]$.
    The natural surjection $\C[\eta,\xi,x_0,x_1,t_0,t_1] \to \C[x_0,x_1,t_0,t_1]$, by substituting the quadratic forms $\eta(x_i)$ and $\xi(x_i)$ in the variables $x_0$, $x_1$ for the variables $\eta$ and $\xi$, induces an inclusion $\P(1^2,2^2)\times\P^1 \supset \P^1\times\P^1 = S$, and $D$ is the intersection of $S$ and $V(f) \subset \P(1^2,2^2)\times\P^1$.
    For general $a$ and $b$, $D$ is irreducible by Bertini's theorem.
    Since the expressions $x_2=\xi x_0/\eta$ and $x_3=\xi x_1/\eta$ are unique, $p : C_V \to D$ is bijective, and no component of $C$ contracts to a point by $p$.
    Thus we can see that $\rho(V)=2$.

    Since $-K_V\sim H_V$, the minimal section $s_0=\P[0]\subset X$ is the unique flopping curve of $V$.
    The $(-F_V)$-flop $V'$ of $V$ is also a weak Fano $3$-fold with one extremal ray $R'$.
    We will fix the ray $R'$ and its contraction morphism.
    For any curve $C\equiv al+bs_0$ in $V$, $n$ denotes $\#\{C\cap s_0\}$.
    Since each member of $|H_V-2F_V|$ has singularities along $s_0$ and since $Bs|H_V-2F_V|=s_0$, it follows that $2n \leq C\cdot (H_V-2F_V) = a-2b$, and that $\mu(C') = a/(b+n) \geq 2$ by (\ref{eqn:slope}).
    Cutting $V$ by two general members of $|H_V-2F_V|$, we obtain the curve $D\subset \P[0,1^2,2]$ defined by $M_D$ with
$$w(M_D)=\WM(-3,-3,-2,-2,-2).$$
    The bihomogeneous coordinate ring of $\P[0,1^2,2]$ is $R=\C[x,y_1,y_2,z,t_0,t_1]$ with $\deg x=(1,0)$, $\deg y_i=(1,-1)$, $\deg z=(1,-2)$, and $\deg t_j=(0,1)$.
    Choosing the appropriate generators of $R$ and making a linear transform of $M_D$, we have
$$
 M_D=\pmatrix{
  0 &  x & a & b &y_1\cr
 -x &  0 & c & d &y_2\cr
 -a & -c & 0 & z & 0 \cr
 -b & -d &-z & 0 & 0 \cr
-y_1&-y_2& 0 & 0 & 0},
$$
where $a$, $b$, $c$ and $d$ are general linear forms in $y_1$, $y_2$, $zt_0$ and $zt_1$.
    Thus the curve defined by the system of equations below:
\begin{eqnarray*}
f_5 &=& xz-ad+bc, \\
f_4 &=& -ay_2+cy_1, \\
f_3 &=& -by_2+dy_1, \\
f_2 &=& y_1z, \\
f_1 &=& y_2z.
\end{eqnarray*}
    It follows that $y_1=y_2=0$ from $f_4=f_3=f_2=f_1=0$.
    Therefore, the curve $D$ is decomposed into $D \equiv s_0+D_1$, where $s_0 = \{ y_1=y_2=z=0 \}$ and $D_1 = \{ x=y_1=y_2=0 \} \equiv s_0+2l$.
    The strict transform $D_1'$ of $D_1$ is of slope $\mu(D_1')=2$, because $D_1$ does not meet with the flopping curve $s_0$.
    Consequently, the extremal ray $R'$ is of slope $\mu(R')=2$, and its contraction morphism is defined by $|H_V-2F_V|$.
    We can see that the morphism gives a $\P^1$-bundle structure, and the ray is of type $C_2$.

    The anti-canonical model $\ol{V}$ is an intersection of quadrics in $\P^{11}$ with one ODP. This is (\ref{thm:df5}.1).
\cons{37}
    In this case, $X=\P[0^5,1]$, and $V$ is defined by $M$ with $m_{ij} \in H^0(X,\O(H))$.
    Consider $X$ as the intersection of a member of $|H_Y+F_Y|$ and three members of $|H_Y|$ in $Y=\P[0^{10}]=\P^9\times\P^1$.
    The subvariety $V$ is the intersection $X\cap G$, where $G$ is defined by $M_Y$ with $m_{ij} \in H^0(Y, \O(H_Y))$, i.e., $G\cong Gr(5,2)\times\P^1 \subset \P^9\times\P^1$.
    Thus, $V$ for general $M$ is smooth by $Bs|H_Y| = Bs|H_Y+F_Y| = \emptyset$.
    A similar argument as (\ref{subsec:df5}.1) shows $\rho(V)=2$ if $C_V$ is irreducible, where $C_V$ is defined by $M_C$ with $\dual{m_{ij}} \in H^0(\dual{X}, \O(\dual{H}))$ in $\dual{X}=\P[0^3,1]$.
    The curve $C_V$ is the intersection $\dual{X}\cap {}_1C_5\times\P^1$ in $\P^4\times\P^1 = \P[0^5]$ for the quintic elliptic curve ${}_1C_5\subset\P^4$, the general hyperplane section of a del Pezzo surface of degree $5$ in $\P^5$.
    Here, $C_V$ is irreducible because $\dual{X}$ is ample in $\P[0^5]$, and $\rho(V)=2$.

    Since $-K_V \sim H_V+F_V$ is ample, $V$ is a Fano $3$-fold and has two extremal rays.
    The one corresponds to the fibration.
    The other $R'$ is generated by minimal sections $s_0=\P[0] \subset X$.
    Its contraction morphism is defined by $|H_V|$ with $\dim|H_V|=6$ and $(H_V)^3=5$.
    The image is a Fano $3$-fold $B_5\subset\P^6$ with index $2$ and degree $5$.
    The morphism contracts $V\cap\P[0^3]$ to a quintic elliptic curve ${}_1C_5$ in $B_5$.
    This is (\ref{thm:df5}.2).
\cons{24}
    In this case, we have $V$ is defined by $M$ with
$$w(M)=\WM(-2,-2,-2,-1,-1)$$
in $X=\P[0,1^4,2]$.
    The $\P^5$-bundle $X$ is the intersection $H_1\cap H_2\cap H_3\cap H_4$ in $Y=\P[0^3,1^6,2]$, $H_1, H_2\in |H_Y|$ and $H_3, H_4\in |H_Y-F_Y|$.
    Let $G$ be the $Gr(5,2)$-bundle defined by $M_Y$ with $w(M_Y)=w(M)$.
    Then, $V$ is the intersection $G\cap X$, and $V$ is smooth outside $s_0=\P[0] \subset X$.
    Calculation with local coordinates along $s_0$ shows that general $V$ has no singularity along $s_0$.
    We now consider the duals $\dual{Y}=\P[0^3,(-1)^6,-2]$, $\dual{X}=\P[0^2,(-1)^2]$ and $\dual{G}$ defined by $\dual{M}$ with
$$w(\dual{M})=\WM(2,2,2,1,1).$$
    Denote by $\dual{Z}$ the sub-$\P^4$-bundle $\P[0^3,(-1)^2]$ containing $\dual{X}$.
    We have $\dual{Z} = \dual{L_1}\cap\dual{L_2}\cap\dual{L_3}\cap\dual{L_4}\cap\dual{L_5}$, where $\dual{L_1} \in |H_{\dual{X}}+2F_{\dual{X}}|$ and $\dual{L_2}, \dual{L_3}, \dual{L_4}, \dual{L_5}\in |H_{\dual{X}}+F_{\dual{X}}|$.
    On $\dual{L_1}$, the linear system $|H_{\dual{X}}+F_{\dual{X}}|$ is free from base points.
    By Bertini's theorem, $\dual{G}\cap\dual{Z}$ is irreducible.
    We have $\dual{m_{12}}=\dual{m_{13}}=\dual{m_{23}}=0$ on $\P[(-1)^2]$, and $\dual{m_{i4}}, \dual{m_{i5}} \in |H_{\dual{X}}|$, hence $\dual{f_5}=\dual{f_4}=0$ and
\begin{eqnarray*}
 \dual{f_3} = -\dual{m_{14}}\dual{m_{25}}+\dual{m_{15}}\dual{m_{24}}, & & \\
 \dual{f_2} = -\dual{m_{14}}\dual{m_{35}}+\dual{m_{15}}\dual{m_{34}}, & & {\rm and} \\
 \dual{f_1} = -\dual{m_{24}}\dual{m_{35}}+\dual{m_{25}}\dual{m_{34}}. & &
\end{eqnarray*}
    For general $\dual{m_{ij}}$, the system of equations $\dual{f_3}=\dual{f_2}=\dual{f_1}=0$ has no common zero on $\P[(-1)^2]$.
    Hence, $|H_{\dual{X}}|$ has no base point on $\dual{G} \cap \dual{Z}$ because of $Bs|H_{\dual{X}}\big|_{\dual{Z}}| = \P[(-1)^2]$.
    Therefore, $C_V = \dual{G}\cap\dual{X} = (\dual{G}\cap\dual{Z})\cap H_{\dual{Z}}$ is irreducible.

    Since $-K_V\sim H_V$, the minimal section $s_0=\P[0]$ is the unique flopping curve on $V$.
    We obtain another weak Fano $3$-fold $V'$ with one extremal ray $R'$.
    Now we will fix the ray $R'$ and its contraction morphism.
    For any curve $C\equiv a l+b s_0$ on $V$, the inequality $n=\#\{C\cap s_0\} \leq C\cdot(H-F)=a-b$ shows that the slop $\mu(C')$ of the strict transform of $C$ is $a/(b+n)\geq1$ from (\ref{eqn:slope}).
    Denote by $L$ the sub-$\P^4$-bundle $\P[0,1^4] \subset X$ which is the unique member of $|H_X-2F_X|$.
    Since $m_{45}=0$ on $L$, the zero locus of $\{f_3,f_2,f_1\}$ is ${}_0C_3\times\P^1 \subset \P^3\times\P^1\cong\P[1^4] \subset L$.
    Each $\alpha\in{}_0C_3$ determines a curve $D_\alpha = V\cap S_\alpha$ where $S_\alpha=\P[0,1] \subset L$.
    The curve $D_\alpha$ is decomposed into $s_0+C_\alpha$, where $C_\alpha\equiv s_0+l$, i.e., $\mu(C_\alpha)=1$, and $C_\alpha\cap s_0=\emptyset$ for general $\alpha$.
    Hence $R'$ is generated by $C_\alpha'$ and $\mu(R')=1$.
    Consequently, the contraction morphism of $R'$ is defined by $|H_{V'}-F_{V'}|$ with $\dim|H_{V'}-F_{V'}|=5$ and $(H_{V'}-F_{V'})^3=4$;
    the image of the morphism is a Fano $3$-fold $B_4\subset\P^5$ with index $2$ and degree $4$;
    and the unique member of $|H_{V'}-2F_{V'}|$ is contracted to a twisted cubic curve ${}_0C_3\subset B_4\cap\P^3$.

    The anti-canonical model $\ol{V}$ is an intersection of quadrics in $\P^{11}$, and has one ODP.
    This is (\ref{thm:df5}.3).
\cons{6}
    In this case, $X=\P[0,1^5]$ and $V$ is defined by $M$ with
$$w(M)=\WM(-2,-1,-1,-1,-1).$$
    Consider $X$ as the intersection $H_1\cap H_2\cap H_3\cap H_4$, $H_1, H_2, H_3\in |H_Y|$ and $H_4\in |H_Y-F_Y|$ in $Y=\P[0^4,1^6]$,.
    Denote by $G$ the $Gr(5,2)$-bundle in $Y$ defined by $M_Y$ with $w(M_Y)=w(M)$.
    Then, $V$ is the intersection $G\cap X$ and is smooth outside $s_0=\P[0]\subset X$.
    Calculation using the local coordinates on $X$ along $s_0$ shows that $V$ is smooth also along $s_0$.
    As similar to (\ref{subsec:df5}.3), consider the duals $\dual{Y}=\P[0^4,(-1)^6]$, $\dual{X}=\P[0^3,-1]$, and $\dual{G}$ defined by $\dual{M}$ with
$$w(\dual{M})=\WM(2,1,1,1,1).$$
    Then $\dual{X} = \dual{L_1}\cap\cdots\cap\dual{L_6}$, where $\dual{L_1},\dots,\dual{L_5}\in |H_{\dual{Y}}+F_{\dual{Y}}|$ and $\dual{H_6}\in |H_{\dual{Y}}|$.
    Denote by $\dual{L}=\P[0^4,-1]$ the intersection of general five members $\dual{L_1},\dots,\dual{L_5}\in |H_{\dual{Y}}+F_{\dual{Y}}|$.
    The curve $C_V=\dual{G}\cap\dual{X}=\dual{G}\cap\dual{L}\cap H_{\dual{L}}$ is irreducible because $Bs|H_{\dual{L}}|=Bs|H_{\dual{Y}}|\cap\dual{L}=\P[-1]$ and $\dual{G}\cap\dual{L}\cap\P[-1]=\emptyset$.
    Hence $\rho(V)=2$.

    Since $-K_V\sim H_V$, the minimal section $s_0\subset V$ is the unique flopping curve.
    We have another weak Fano $3$-fold $V'$ with one extremal ray $R'$.
    The equation (\ref{eqn:slope}) implies $\mu(R')\geq1$ as above (\ref{subsec:df5}.3).
    One can see that $f_1=Q_0(x)y_0+Q_1(x)y_1$ and $f_2=Q_2(x)$ by choosing homogeneous coordinates $[x_0:\cdots:x_4]\times[y_0:y_1]$ on $\P^4\times\P^1\cong\P[1^5]\subset X$, where $Q_i(x)=l_{i0}x_0+l_{i1}x_1+l_{i2}x_2$ for $l_{ij}$ linear forms in $x_0,\dots,x_4$.
    Denote by $C$ the curve $\{x | Q_0(x)=Q_1(x)=Q_2(x)=0\}\subset\P^4$.
    The curve $C$ lies on a del Pezzo surface of degree $4$ in $\P^4$ and the line $l=\{x | x_0=x_1=x_2=0\}\subset\P^4$ is an irreducible component of $C$.
    Hence, $C$ is decomposed into $l+C_0$ and $C_0$ is an irreducible curve of degree $7$ and genus $3$ for general $l_{ij}$.
    For any $\alpha\in C_0$, each of the three functions $f_5$, $f_4$, and $f_3$ defines the same curve $D_\alpha$ on the ruled surface $S_\alpha=\P[0,1] \subset X$.
    The curve $D_\alpha$ is decomposed into $s_0+D_{\alpha,1}$, where $D_{\alpha,1}$ is irreducilbe with $\mu(D_{\alpha,1})=1$ for general $\alpha\in C_0$.
    The strict transform $D_{\alpha,1}'$ of $D_{\alpha,1}$ is of slope $1$, hence the ray $R'$ is generated by $D_{\alpha,1}'$ and $\mu(R')=1$.
    The contraction morphism of $R'$ is defined by $|H_{V'}-F_{V'}|$ with $\dim|H_{V'}-F_{V'}|=4$ and $(H_{V'}-F_{V'})^3=2$.
    Consequently, the morphism is the reverse of the blow-up of $\Q^3$ along ${}_3C_7$.

    The anti-canonical model $\ol{V}$ is an intersection of quadrics in $\P^{10}$ with one ODP.
    This is (\ref{thm:df5}.4).
\cons{29}\Label{subsubsec:no29}
    In this case, $V$ is defined by $M$ with
$$w(M)=\WM(-1,-1,-1,-1,0)$$
in $X=\P[0^2,1^4]\subset Y=\P[0^6,1^4]$.
    Then, $X=H_1\cap H_2\cap H_3\cap H_4$ for $H_i\in|H_Y|$, and $V$ is smooth because $Bs|H_Y|=\emptyset$.
    Consider the duals $\dual{Y}=\P[0^6,(-1)^4]$, $\dual{X}=\P[0^4]$ and $\dual{G}$ defined by $\dual{M}$ with
$$w(M)=\WM(1,1,1,1,0)$$.
    Denote by $\dual{L}$ the intersection $\P[0^6]=\dual{L_1}\cap\cdots\cap\dual{L_4}$ for $\dual{L_i}\in|H_{\dual{Y}}+F_{\dual{Y}}|$ with $Bs|H_{\dual{Y}}+F_{\dual{Y}}|=\emptyset$.
    Then, $\dual{X}$ is the intersection of two members of $|H_{\dual{L}}|$ with $Bs|H_{\dual{L}}|=\emptyset$.
    Hence the curve $C_V=\dual{G}\cap\dual{X}=(\dual{G}\cap\dual{L})\cap\dual{X}$ is irreducible, and $\rho(V)=2$.

    Since $-K_V\sim H_V$, the minimal sections $\P[0]\subset X$ are flopping curves.
    Restricting $V$ to $\P[0^2]\subset X$, we can see that there are two minimal sections $s_0$ and $s_1$ in $V$.
    The weak Fano $3$-fold $V'$ flopped $V$ has the extremal ray $R'$ with $\mu(R')\geq a/(b+n)\geq1$ by (\ref{eqn:slope}), because $Bs|H_V-F_V|=s_0\cup s_1$.
    We will fix the ray $R'$ and the contraction morphism.
    Let $R$ be the bihomogeneous coordinate ring $\C[x_0,x_1,y_2,y_3,y_4,y_5,t_0,t_1]$ of $X$ with $\deg x_i=(1,0)$, $\deg y_i=(1,-1)$ and $\deg t_j=(0,1)$.
    Choosing the generator of $R$ and making the linear transform of $M$, we may write
$$M=\pmatrix{
  0  & x_0 & x_1 &  a  & y_2 \cr
-x_0 &  0  &  b  &x_1+c& y_3 \cr
-x_1 & -b  &  0  & x_0 & y_4 \cr
 -a &-x_1-c&-x_0 &  0  & y_5 \cr
-y_2 &-y_3 &-y_4 &-y_5 &  0},$$
where $a$, $b$ and $c$ are bilinear forms in $[y_i]$ and $[t_j]$.
    Consider the intersection $D_\lambda=V\cap L_\lambda$ for $\lambda=[\lambda_2:\lambda_3:\lambda_4:\lambda_5] \in \P^3\cong \P(H^0(X,\O(H-F)))$, defined by
\begin{eqnarray*}
f_5 &=& x_0^2 - x_1(x_1+cy) + aby^2, \\
f_4 &=& (\lambda_4x_0 - \lambda_3x_1 +b\lambda_2y)y, \\
f_3 &=& (\lambda_5x_0 - a\lambda_3y + \lambda_2(x_1+cy))y, \\
f_2 &=& (\lambda_5x_1 - a\lambda_4y + \lambda_2x_0)y, \qquad {\rm and} \\
f_1 &=& (\lambda_5b-\lambda_4(x_0+cy)+\lambda_3x_0)y,
\end{eqnarray*}
where $a$, $b$ and $c$ are bilinear forms in $\lambda=[\lambda_i]$ and $[t_j]$.
    For general $\lambda\in\P^3$, $D_\lambda$ is the union $s_0\cup s_1$ of the minimal sections.
    Let $C = \{ \lambda\in\P^3 | a_j(\lambda_3^2-\lambda_4^2) + b_j(\lambda_5^2-\lambda_2^2) - c_j(\lambda_2\lambda_3+\lambda_4\lambda_5) = 0, j=0,1 \}$ for $a=a_0t_0+a_1t_1$, $b=b_0t_0+b_1t_1$ and $c=c_0t_0+c_1t_1$.
    The curve $C\subset\P^3$ is decomposed into three irreducible components $Z_1+Z_2+C_0$, where $Z_1=\{ [\mu:\nu:\nu:-\mu]\in\P^3 \}$, $Z_2=\{ [\mu:\nu:-\nu:\mu]\in\P^3 \}$ are lines, and meet $C_0$ at $4$ points. The curve $C_0$ is of degree $7$ and genus $4$, becouse of $\deg C = 9$ and $p_a(C) = 10$.
    For $\lambda\in C_0$, the intersection $D_\lambda$ has the third component $C_\lambda\equiv l+s_0$ disjoint to the minimal sections $s_0$ and $s_1$.
    The strict transform $C_\lambda'$ is of slope $\mu(C_\lambda')=1$, and hence $\mu(R')=1$.
    As a result, the contraction morphism of $R'$ is defined by $|H_{V'}-F_{V'}|$, the reverse of the blowing-up of $\P^3$ along the curve $C_0$ with $\deg C_0 = 7$ and $g(C_0) = 4$.
    The anti-canonical model $\ol{V}$ is an intersection of quadrics in $\P^9$ with two ODP's. This is (\ref{thm:df5}.5).
\cons{33}\Label{subsubsec:no33}
    In this case, $V$ is defined by $M$ with
$$w(M)=\WM(-1,-1,0,0,0)$$
in $X=\P[0^3,1^3]$.
    By regarding $X$ as a subbundle of $Y=\P[-1,0^6,1^3]$, it follows that $X=H_1\cap H_2\cap H_3\cap H_4$ for $H_1\in|H_Y+F_Y|$, $H_i\in|H_Y|$, $i=2,3,4$.
    Since $Bs|H_Y+F_Y|=\emptyset$ and $Bs|H_Y\big|_{H_1}|=\emptyset$, the intersection $V=G\cap H_1\cap H_2\cap H_3\cap H_4$ is smooth for genenral $H_i$.
    Consider the duals $\dual{Y}=\P[1,0^6,(-1)^3]$, $\dual{X}=\P[1,0^3]$ and $\dual{G}$.
    Since $\dual{X}$ is the intersection $\bigcap_{i=1}^6 \dual{L_i}$ for $\dual{L_1}, \dual{L_2}, \dual{L_3}\in |H_{\dual{Y}}+F_{\dual{Y}}|$ and $\dual{L_4}, \dual{L_5}, \dual{L_6}\in |H_{\dual{Y}}|$, and since $Bs|H_{\dual{Y}}+F_{\dual{Y}}|=\emptyset$ and $Bs|H_{\dual{Y}}\big|_{\dual{L}}|=\emptyset$ for $\dual{L}=\P[0,1^6]=\bigcap_{i=1}^3 \dual{L_i}$, the intersection curve $C_M = \dual{G}\cap\dual{X} = \dual{G}\cap (\dual{L}\cap \dual{L_4}\cap\dual{L_5}\cap\dual{L_6})$ is irreducible, and hence $\rho(V)=2$.

    There are three $K$-trivial curves on $V$, because $-K_V\sim H_V$ and the skew symmetric matrix $M$ defining $V$ is
$$M=\pmatrix{
  0  &  0  & x_0 & x_1 & x_2 \cr
  0  &  0  &  a  &  b  &  c  \cr
-x_0 & -a  &  0  &  0  &  0  \cr
-x_1 & -b  &  0  &  0  &  0  \cr
-x_2 & -c  &  0  &  0  &  0}$$
on $\P[0^3]\subset X$ for linear forms $a$, $b$, $c$ in $[x_0,x_1,x_2]$, where $[x_0,x_1,x_2]$ is the appropriate homogeneous coordinates on a fiber $\P^2$ of $\P[0^3]\to\P^1$.
    Therefore, we obtain the $(-F_V)$-flop $V'$ of $V$, which is also a weak Fano $3$-fold with one extremal ray $R'$.
    A similar argument to (\ref{subsubsec:no29}) shows that the ray $R'$ is of slope $1$, that its contraction morphism is defined by $|H_{V'}-F_{V'}|$, and that the morphism is a conic bundle over $\P^2$ with discrimant locus $\Delta$ of degree $5$.
    The anti-canonical model $\ol{V}$ is an intersection of quadrics in $\P^8$ with three ODP's. This is (\ref{thm:df5}.6).
\cons{2}
    In this case, $V\subset X=\P[0^4,1^2]$ is defined by $M$ with
$$w(M)=\WM(0,0,0,0,0).$$
    This case is very similar to the above case (\ref{subsubsec:no33}) and easer than that.
    Considering $X$ as the subbundle of $Y=\P[0^{10}]\cong\P^9\times\P^1$ obtained by the intersection $H_1\cap H_2\cap H_3\cap H_4$ for $H_1, H_2\in |H_Y+F_Y|$ and $H_3, H_4\in |H_Y|$, we can see that $V$ is smooth for general $H_i$ because $Bs|H_Y+F_Y|=\emptyset$ and $Bs|H_Y|=\emptyset$.
    For the duals $\dual{Y}=\P[0^{10}]$, $\dual{X}=\P[1^2,0^2]$ and $\dual{G}$, we have that the intersection curve $C = \dual{G}\cap\dual{X} = \dual{G}\cap (\dual{L_1}\cap\dual{L_2}\cap\dual{L_3}\cap\dual{L_4})$ is irreducible for general $\dual{L_1}, \dual{L_2}\in |H_{\dual{Y}}+F_{\dual{Y}}|$ and general $\dual{L_3}, \dual{L_4}\in |H_{\dual{Y}}|$, and that $\rho(V)=2$.
    Since $-K_V\sim H_V$, and since the restriction $V\cap \P[0^4]$ of $V$ to $\P[0^4]\subset X$ is the direct product $\{ 5~{\rm points} \}\times\P^1$, the each $\P^1$ of the product is a $K$-trivial curve in $V$.
    Hence, there are only five $K$-trivial curves on $V$, and there exists the $(-F_V)$-flop $V'$ of $V$.
    The $(-F_V)$-flop $V'$ is also a weak Fano $3$-fold with one extremal ray $R'$ of slope $\mu(R')=1$, because $Bs|H_V-F_V|$ is the union $\{ 5~{\rm points} \}\times\P^1$ of the flopping curves, and because the strict transform of the curve $V\cap\P[0^3,1]$ is of slope $1$.
    Consequently, the contraction morphism of $R'$ is defined by $|H_{V'}-F_{V'}|$ and the morphism is again the del Pezzo fibration of degree $5$ over $\P^1$.
    The anti-canonical model $\ol{V}$ is the intersection $Gr(5,2)\cap H_1\cap H_2\cap Q\subset \P^9$, where $H_1$ and $H_2$ are general hyperplanes, and $Q$ is a quadric hypersurface of rank $4$.
    The singularities of $\ol{V}$ come from the singular locus $\Sigma\cong\P^5$ of $Q$, and are the intersection $\Sigma\cap Gr(5,2)\cap H_1\cap H_2$, five ODP's. This is (\ref{thm:df5}.7).
\cons{38}
    This is the case that $V$ is defined by $M$ with
$$w(M)=\WM(0,0,0,1,1)$$
in $X=\P[0^5,1]$.
    This case is also similar to the case (\ref{subsubsec:no33}).
    Regarding $X$ as the subbundle of $Y=\P[(-1)^3,0^6,1]$, and considering the duals $\dual{Y}=\P[1^3,0^6,-1]$, $\dual{X}=\P[1^3,0]$ and $\dual{G}$, we can see that $V$ is smooth and that $\rho(V)=2$.
    Choosing the appropriate generators of the bihomogeneous coordinate ring $R=\C[x_0,x_1,x_2,x_3,x_4,y,t_0,t_1]$ of $X$, where $\deg x_i=(1,0)$, $\deg y=(1,-1)$, $\deg t_j=(0,1)$, and making a linear transform on $M$, we can assume that the defining matirx of $V$ is
$$M = \pmatrix{
  0 & g_1& g_2& x_0& x_1\cr
-g_1&  0 & g_3& x_2& x_3\cr
-g_2&-g_3&  0 & x_3& x_4\cr
-x_0&-x_2&-x_3&  0 &  y \cr
-x_1&-x_3&-x_4& -y &  0}$$
for $g_p\in R$ of bidegree $(1,1)$, i.e., $g_p=g_{p0}t_0+g_{p1}t_1+yg_{py}$ with linear forms $g_{p0}, g_{p1}$ in $[x_0,\dots,x_4]$ and a quadratic form $g_{py}$ in $[t_0,t_1]$.
   Then, $V\subset X$ is defined by
\begin{eqnarray*}
f_5 &=& g_1x_3 - g_2x_2 + g_3x_0, \\
f_4 &=& g_1x_4 - g_2x_3 + g_3x_1, \\
f_3 &=& g_1 y - x_0x_3 + x_1x_2, \\
f_2 &=& g_2 y - x_0x_4 + x_1x_3, \qquad {\rm and} \\
f_1 &=& g_3 y - x_2x_4 + x_3^2.
\end{eqnarray*}
     Denote by $L$ the subbundle $\P[0^5]\subset X$ defined by $y$.
    From $-K_V\sim H_V$, it follows that the $K$-trivial curves of $V$ are in $V\cap L$ and that they are minimal sections of $L\to\P^1$.
    On $L$, the zero locus $\{ f_i\big|_L=0 \mbox{ for } i=1,2,3 \}$ is the trivial fibration $\Sigma_1\times\P^1$ of the rational ruled surface $\Sigma_1\cong\P[0,1]$ in a fiber $\P^4$ of $L\to\P^1$.
    Moreover, on $L$, we can write $f_i\big|_L = f_{i0}t_0 + f_{i1}t_1$, $i=4,5$, by using the quadratic forms $f_{ij} = g_{1j}x_{8-i} - g_{2j}x_{7-i} + g_{3j}x_{5-i}$, $i=4,5$, in $[x_0,\dots,x_4]$.
    The desired minimal sections are obtained by $\{ \alpha \}\times\P^1$ for $\alpha\in A=\{ \alpha\in\Sigma_1 | f_{ij}(\alpha)=0 \mbox{ for } i=4,5, j=0,1 \}=\{ \alpha\in\P^4 | f_i(\alpha,0)=0 \mbox{ for } i=1,2,3 \mbox{ and } f_{ij}(\alpha)=0 \mbox{ for } i=4,5, j=0,1 \}$.
    Denote by $\Lambda$ the linear system $\{ (h_1x_2 - h_2x_1 + h_3x_0)_0 | h_i \mbox{ are linear forms in } [x_0,\dots,x_4] \}$ on $\Sigma_1$.
    The curve $\{ \alpha\in\Sigma_1 | f_{5j}(\alpha)=0 \}$ is in $\Lambda$, and is decomposed into two curves $f+C_j$ for each $j=0,1$, because $Bs \Lambda=f$, where $f=\{ x_0=x_1=x_2=0 \}\subset\Sigma_1$ is a fiber of $\Sigma_1\to\P^1$.
    The curve $C_j=\{ \alpha\in\Sigma_1 | f_{4j}(\alpha)=f_{5j}(\alpha)=0 \}$, $j=0,1$, is numerically equivalent to $2e + 3f$, where $e$ is the minimal sectoin of $\Sigma_1\to\P^1$.
    The curves $C_j$ for general $g_p$ are smooth by $Bs (\Lambda-f)=\emptyset$.
    Since $A=C_0\cap C_1$, the number of the minimal sections, $K$-trivial curves, is $(C_0\cdot C_1)_{\Sigma_1} = (2e + 3f)^2 = -4+12=8$.
    Consequently, there exists the $(-F_V)$-flop $V'$ of $V$, and $V'$ is a weak Fano $3$-fold with one extremal ray $R'$.

    Now, we will study the ray $R'$ and its contraction morphism.
    To do this, we consider the linear system $|2H_V-F_V|$ on $V$.
    The linear system $|2H_V-F_V|$ is not $|2H-F|\big|_V$, but contains the other member $D_V$ coming from an irreducible member $D\in|2H|$ such that $D\big|_V=D_V+F_V$.
    Therefore, $Bs|2H_V-F_V|$ is the intersection of $Bs|2H-F|\cap V$ and $\bigcap\{ \mbox{ such } D_V \}$.
    Such the irreducible member $D\in|2H|$ is $(f_{50}\mu_0+f_{51}\mu_1+yf_{5y}(\mu))_0$ for $f_{5y}(\mu) = g_{1y}(\mu)x_3-g_{2y}(\mu)x_2+g_{3y}(\mu)x_0$, where $g_{py}(\mu)$ is obtained by substituting $\mu_j$ for $t_j$ in $g_{py}$.
    Thus we have that $Bs|2H_V-F_V| = (\{ y=0 \}\cap V)\cap \{ f_{50}=f_{51}=0 \}$ is the union of the eight flopping curves.
    Similarly as others, since $n=\#\{ C\cap\{ 8 \mbox{ flopping curves\} } \}\leq (C\cdot 2H_V-F_V) = 2a-b$ for any curve $C\equiv al+bs_0 \subset V$, (\ref{eqn:slope}) shows that $\mu(C')=a/(b+n)\geq1/2$ and hence $\mu(R')\geq1/2$.
    Let $P_\lambda$ be the sub-$\P^1$-bundle $\P[0^2]\subset X$ associated to $[x_0:x_1:x_2:x_3:x_4:y] = [\lambda_0z_0:\lambda_1z_0:\lambda_0^2z_1:\lambda_0\lambda_1z_1:\lambda_1^2z_1:0]$, and $C_\lambda$ the intersection curve $V\cap P_\lambda$.
    Then, $C_\lambda$ is defined by $g_1(\lambda,z)\lambda_1z_1-g_2(\lambda,z)\lambda_0z_1+g_3(\lambda,z)z_0$ on $P_\lambda$, and numerically equivalent to $2s_0+l$, i.e., $\mu(C_\lambda)=1/2$.
    For gneral $\lambda=[\lambda_0:\lambda_1]\in\P^1$, $C_\lambda$ is irreducible and does not meet the eight minimal sections, hence $\mu(R')=\mu(C_\lambda')=1/2$.
    Therefore, the contraction morphism of $R'$ is defined by the linear system $|2H_{V'}-F_{V'}|$, contracting the strict transform of $V\cap L$ to a line.
    Since $\dim|2H_{V'}-F_{V'}|=\dim|2H_V-F_V|=8$ and $(2H_{V'}-F_{V'})^3=(2H_V-F_V)^3+8=12$, the image of the contracion morphism is a Fano $3$-fold $V_{12}\subset\P^8$ of index $1$ with genus $7$.
    Consequently, $V'$ is the blowing-up of $V_{12}$ along a line.

    The anti-canonical model $\ol{V}$ is the complete intersection $Q_1\cap Q_2\cap Q_3\subset\P^6$ of three smooth quadrics, defined by \begin{eqnarray*}
Q_1 &=& \{ g_{10}x_5+g_{11}x_6+g_{1y}(x_5,x_6) - x_0x_3 + x_1x_2=0 \}, \\
Q_2 &=& \{ g_{20}x_5+g_{21}x_6+g_{2y}(x_5,x_6) - x_0x_4 + x_1x_3=0 \}, \qquad {\rm and}\\
Q_3 &=& \{ g_{30}x_5+g_{31}x_6+g_{3y}(x_5,x_6) - x_2x_4 + x_3^2=0 \},
\end{eqnarray*}
for the homogeneous coordinates $[x_0:\cdots:x_6]$ on $\P^6$.
    Since the linear forms $g_{p0}$, $g_{p1}$, and the quadratic forms $g_{py}$ are general, the singularities of $\ol{V}$ lie on $\P^4 = \{ x_5=x_6=0 \} \subset \P^6$.
    The calculation of the Jacobian shows that the singularities are eight ODP's $\{ \alpha\in\P^6 | f_{ij}(\alpha)=0 \mbox{ for } i=4,5, j=0,1, \mbox{ and } \alpha_5=\alpha_6=0 \}$, which correspond to the minimal sections in $V\subset X$.
    Thus, we obtain (\ref{thm:df5}.8).
\cons{16}
    In this case, $V$ is defined by $M$ with
$$w(M)=\WM(0,1,1,1,1)$$
in $X=\P[0^6]$.
    As simlir to others, the smoothness of $V$ comes from the identification $X$ with the subbundle of $Y=\P[(-1)^4,0^6]$, and $\rho(V)=2$ from the consideration of the duals $\dual{Y}$, $\dual{X}$ and $\dual{G}$.
    Without loss of generality, we can choose the homogenious coordinates on $X\cong\P^5\times\P^1$ and write down the defining equations of $V$ as follows:
\begin{eqnarray*}
f_5 &=& g_1x_3-g_2x_1+g_3x_0, \\
f_4 &=& g_1x_4-g_2x_2+g_4x_0, \\
f_3 &=& g_1x_5-g_3x_2+g_4x_1, \\
f_2 &=& g_2x_5-g_3x_4+g_4x_3, \qquad {\rm and} \\
f_1 &=& x_0x_5-x_1x_4+x_2x_3,
\end{eqnarray*}
where $g_p$ are bilinear forms in $[x_0,\dots,x_5]\times[t_0,t_1]$.
    Each $f_i$, $i=2,\dots,5$, is decomposed into $f_{i0}t_0+f_{i1}t_1$ by quadric forms $f_{ij}$ in $[x_0,\dots,x_5]$.
    Since $-K_V\sim H_V$, the $K$-trivial curves in $V$ are the minimal sections of $V\subset\P^5\to\P^1$, and these sections correspond to the points of the set $A = \{ \alpha\in\P^5 | f_1(\alpha)=f_{ij}(\alpha)=0 \mbox{ for } i=2,\dots,5 \mbox{ and } j=0,1 \}$.
    Since the equation $f_1=0$ defines the quadric $4$-fold $\Q^4$, the set $A$ is identified the intersection of two surfaces $S_j$ ($j=0,1$) in $\Q^4$, where $S_j$ is defined by $f_{2j}=\cdots=f_{5j}=0$.
    The intersection of general hyperplane sections of $\Q^4$ decomposes into two irreducible $2$-cycles $D_1+D_2$, and the intersection numbers of these cycles are $(D_i)^2=1$ and $(D_1 \cdot D_2)=0$.
    We can see that each $S_j$ is equivalent to $2D_1+3D_2$ as $2$-cycles and that the intersection number is $(S_0 \cdot S_1)=4+9=13$.
    For general $g_i$, each point of $A$ are simple, and hence the number of $K$-trivial curves in $V$ equals to $\#A=13$.
    Thus, the $(-F_V)$-flop $V'$ of $V$ exists, and has one extremal ray $R'$. 

    To study the ray $R'$, we consider the linear system $|2H_V-F_V|$.
    We can calculate that $\dim|2H_V-F_V|=3$ and $(2H_V-F_V)^3=8\times6-12\times5=-12$.
    No member of $|2H_V-F_V|$ comes from the restriction of the member of $|2H-F|$ on $X$, because $|2H-F|=\emptyset$.
    Each member $D_V\in|2H_V-F_V|$ is the irreducible component of the reducible member $D\big|_V=D_V+F_V\in|2H_V|$, where $D$ is a member of the linear system $|2H|$ on $X$.
    Let $D_i$ be the restriction of $(f_{i0})_0 \in |2H|$ to $V$ for $i=2,\dots,5$, then $D_i$ is reducible and has the irreducible component $D_{iV}=D_i-F_V\in|2H_V-F_V|$.
    Therefore, $Bs|2H_V-F_V|=D_{2V}\cap\cdots\cap D_{5V}$ is the union of the thirteen flopping curves.
    Let $C'\subset V'$ be any curve, and $C\equiv al+bs_0$ its birational transform.
    The inequality $\mu(C')=a/(b+n)\geq1/2$ comes from (\ref{eqn:slope}), because $n=\#\{ C\cap \{ 13 \mbox{ flopping curves\} }\}\leq (C\cdot 2H_V-F_V)=2a-b$.
    Thus, the contraction morphism of the extremal ray $R'$ of the flop $V'$ is defined by the linear system $|2H_{V'}-F_{V'}|$ with $\dim |2H_{V'}-F_{V'}|=3$ and $(2H_{V'}-F_{V'})^3 = (2H_V-F_V)^3+13 = 1$.
    We can see that the image of the morphism is $\P^3$ and the morphism is the (reverse of) blowing-up of $\P^3$ along a curve $_{8}C_{9}$ of degree $9$ and of genus $8$.

    The anti-canonical model $\ol{V}$ is the compete intersection of the quadric $Q$ defined by $x_0x_5-x_1x_4+x_2x_3$ and the cubic $4$-fold defined by $g_{[3,4]}x_0-g_{[2,4]}x_1+g_{[2,3]}x_2+g_{[1,4]}x_3-g_{[1,3]}x_4+g_{[1,2]}x_5$, where $g_{[i,j]}=g_{i0}g_{j1}-g_{i1}g_{j0}$.
    The singularities of $\ol{V}$ are $f_1=f_{ij}=0$, which are $13$ ODP's.
\cons{28}
    In this case, $V$ is defined by $M$ with
$$w(M)=\WM(-1,-1,-1,-1,0)$$
in $X=\P[0,1^5]$.
    As simlir to others, the smoothness of $V$ comes from the identification $X$ with the subbundle of $Y=\P[0^6,1^4]$, and $\rho(V)=2$ from the consideration of the duals $\dual{Y}$, $\dual{X}$ and $\dual{G}$.
    By a suitable choice of the bihomogenious coordinates on $X$, the defining equations of $V$ are
\begin{eqnarray*}
f_5 &=& x_0(x_0+y_5t_0)-g_1g_4+g_2g_3, \\
f_4 &=& x_0y_3-g_1y_2+g_3y_1, \\
f_3 &=& x_0y_4-g_2y_2+g_4y_1, \\
f_2 &=& g_1y_4-g_2y_3+(x_0+y_5t_0)y_1, \qquad {\rm and} \\
f_1 &=& g_3y_4-g_4y_3+(x_0+y_5t_0)y_2,
\end{eqnarray*}
where $g_p$ are bilinear forms in $[y_1,\dots,y_5]\times[t_0,t_1]$.
    Let $S_\alpha=\P[0,1]\subset X$ be a ruled surface correspoinding to $\alpha=[\alpha_i]\in\P^4$ by $y_i=\alpha_iy$.
    The intersection $C_\alpha=V\cap S_\alpha$ is defined by
\begin{eqnarray*}
f_5(\alpha) &=& x_0(x_0+\alpha_5t_0y)-(g_1g_4-g_2g_3)y^2, \\
f_4(\alpha) &=& x_0\alpha_3-(g_1\alpha_2-g_3\alpha_1)y, \\
f_3(\alpha) &=& x_0\alpha_4-(g_2\alpha_2-g_4\alpha_1)y, \\
f_2(\alpha) &=& x_0\alpha_1+(g_1\alpha_4-g_2\alpha_3+\alpha_5\alpha_1t_0)y, \qquad {\rm and} \\
f_1(\alpha) &=& x_0\alpha_2+(g_3\alpha_4-g_4\alpha_3+\alpha_5\alpha_2t_0)y.
\end{eqnarray*}
    Since $-K_V\sim H_V-F_V$, the $K$-trivial curves are contained in $S_\alpha$ for $\alpha$ when $f_i(\alpha)=0$ defines the same curve on $S_\alpha=\P[0,1]$.
    These $\alpha$'s build up the set
$$A=\{ \alpha\in\P^5 |\: {\rm rank}\left|
\begin{array}{cccc}
\alpha_3 & -g_{10}\alpha_2+g_{30}\alpha_1 & -g_{11}\alpha_2+g_{31}\alpha_1 \\
\alpha_4 & -g_{20}\alpha_2+g_{40}\alpha_1 & -g_{21}\alpha_2+g_{41}\alpha_1 \\
\alpha_1 & g_{10}\alpha_4-g_{20}\alpha_3+\alpha_5\alpha_1 & g_{11}\alpha_4-g_{21}\alpha_3 \\
\alpha_2 & g_{30}\alpha_4-g_{40}\alpha_3+\alpha_5\alpha_2 & g_{31}\alpha_4-g_{41}\alpha_3
\end{array}\right|\leq1 \},$$
where $g_p=g_{p0}t_0+g_{p1}t_1$.
   The set $A$ containds $\alpha_0=[0:0:0:0:1]$, and $C_{\alpha_0}=\{ x_0(x_0+t_0y)-(g_1g_4-g_2g_3)y^2 = 0\}$ is a bisection of $S_\alpha\to\P^1$, hence of $V\to\P^1$.
    With a slight effort, we can calculates the number $\#|A\setminus\{\alpha_0\}|$, which is $21$.
    Thus there exist $22$ $K$-trivial curves on $V$; one of them is a bisection, the others are sections.
    The anti-canonical model $\ol{V}$ is a singular quartic $3$-fold in $\P^4$, which has exactly $22$ ODP's.

    We now study the flop $V'$ and its extremal ray $R'$.
    The members of the linear system $|3H_V-4F_V|$ does not come from the restriction of the member of $|3H-4F|$, because $|3H-4F|=\emptyset$.
    The members of $|3H-3F|$ whose restriction to $V$ decompose as $D+F_V$ are the zero loci of the following functions:
\begin{eqnarray*}
h_0&=&y_3(-g_{20}y_2+g_{40}y_1)-y_4(-g_{10}y_2+g_{30}y_1), \\
h_1&=&y_3(g_{10}y_4-g_{20}y_3+y_1y_5)-y_1(-g_{10}y_2+g_{30}y_1), \\
h_2&=&y_3(g_{30}y_4-g_{40}y_3+y_2y_5)-y_2(-g_{10}y_2+g_{30}y_1), \\
h_3&=&y_4(g_{10}y_4-g_{20}y_3+y_1y_5)-y_1(-g_{20}y_2+g_{40}y_1), \quad {\rm and} \\
h_4&=&y_4(g_{30}y_4-g_{40}y_3+y_2y_5)-y_2(-g_{20}y_2+g_{40}y_1).
\end{eqnarray*}
    The irreducible components $D_i$ of $(h_i)_0\cap V$ generate $|3H_V-4F_V|$.
    Since $Bs|3H_V-4F_V|$ is the union of the $K$-trivial curves found above, for any curve $C$ on $V$, we have $\mu(C')\geq 4/3$, and hence $\mu(R')\geq 4/3$.
    We can find the curve of slope $4/3$, and the contraction morphism corresponding to the extremal ray $R'$ is defined by $|3H_{V'}-4F_{V'}|$ with $\dim |3H_{V'}-4F_{V'}|=4$ and $(3H_{V'}-4F_{V'})^3 = (3H_V-4F_V)^3 + 21 + 8 = 2$.
    The slightly hard calculation shows the morphism is the (reverse of the) blowing-up of $\Q^3$ along a curve of degree $11$ and with genus $9$.
\cons{31}
    The case is very complicated.
    We will discuss in the additional paper.

\newpar\Label{subsec:kanev}
    We here give a briaf explanation the relation between the irreducibility of the curve $C_V$ and $\rho(V)=2$.

    Let $F$ be a $5$-dimensional vector space, and $Gr(F,2) \subset \P(\bigwedge^2F)$ the Grassmannian manifold parametrizing $2$-dimensional quotient vector spaces of $F$, i.e., parametrizing lines in $\P^4=\P(F)$.
    The Grassmannian manifold $Gr(F,2)$ has the self-dual property: it means that hyperplanes tangent to $Gr(F,2) \subset \P(\bigwedge^2F)$ form a Grassmannian manifold $Gr(\dual{F},2)$ in the dual projective space $\P(\bigwedge^2\dual{F})$ parametrizing hyperplanes in $\P(\bigwedge^2F)$.
    Each point $[H] \in Gr(\dual{F},2)$ determines the hyperplane $H$ tangent to $Gr(F,2)$ in $\P(\bigwedge^2F)$, and the intersection $H \cap Gr(F,2)$ has a singular locus $\Sigma$ isomorphic to $\P^2$.
    For each point $[l] \in \Sigma \subset Gr(F,2)$, the corresponding line $l \subset \P^4=\P(F)$ lies on the same plane $P_H=\P^2 \subset \P^4$ determined by $[H]$.
    Moreover each point in $H \cap Gr(F,2)$ corresponds to a line in $\P^4$ meeting with $P_M$.

    Del Pezzo surfaces of degree $5$ are hyperspace sections of Grassmannian manifold parametrizing lines in $\P^4$; i.e., for any del Pezzo surface $S$ of degree $5$, there exist four hyperplanes $H_i \subset \P(\bigwedge^2F)\cong\P^9$, $i=1,\dots,4$, such that $S=Gr(F,2) \cap H_1 \cap \cdots \cap H_4$, for a $5$-dimensional vector space $F$.
    A $3$-dimensional projective subspace spanned by the duals of $H_i$ in $\P(\bigwedge^2\dual{F})$ meets at $5$ points with $Gr(\dual{F},2)$, say $\{L_1,\dots,L_5\}$.
    Note that $L_1 \cap \cdots \cap L_5 = H_1 \cap \cdots \cap H_4$, i.e., $S=Gr(F,2) \cap L_1 \cap \cdots \cap L_5$.
    Distinct pairs of the $5$ points naturally correspond to lines in $S$ as follows.

    Each $H_i$ defines a plane $P_{L_i} \subset \P^4=\P(F)$, and $P_{L_i}$ meets $P_{L_j}$ at one point $p_{ij}$ for any $i\ne j$.
    The lines meeting $P_{L_3}$ and $P_{L_4}$ through the point $p_{12}$ in $\P^4=\P(F)$ meets $P_{L_5}$.
    These line form a $1$-dimensional family parametrized by a line $l_{12}$ in $Gr(F,2)$.
    The line $l_{12}$ lies on each $L_i$, hence $l_{12}$ is a line in $S$.
    For each pair $(i,j)$ of $i,j=1,\dots,5$, we obtain a line $l_{ij}$ on the del Pezzo surface $S$.
    Thus the ten pairs $(i,j)$ correspond one to one to the ten lines on $S$.

    In the case of del Pezzo fibration,  replace a vector space $F$ by a vector bundle $\F$, and the parallel argument runs.
    Thus the curve $C_V$ parametrizes lines on each fiber of the del Pezzo fibration $\varphi : V \to \P^1$.
    To be $\rho(V)>2$, the horizontal divisor $H_V$ must decompose into at least two horizontal divisors, because every fiber is reducible in our case.
    But $H_V$ may be decomposable only if $C_V$ is reducible.

\section{Del Pezzo fibrations of degree 1} \label{sec:df1}

\newpar
    Finaly, we treat del Pezzo fibrations of degree $1$ over $\P^1$, and derive Theorem(\ref{thm:df1}) and Supplement(\ref{thm:df1s}).
    Let $\S = \bigoplus_{d \geq 0} \S_d$ be a graded $\O_{\p^1}$-algebra which satisfies the following conditions:

(1) ~ $\S_1 = \E = \bigoplus_{i=0}^n \O(a_i)$, ~ $a_0 \leq a_1 \leq \cdots \leq a_n$,

(2) ~ $\S_2 = S^2(\E) \op \F$, ~ $\F = \O(b)$, ~ $b \in \Z$,

(3) ~ $\S_3 = S^3(\E) \op \E \otimes \F \op \G$, ~ $\G = \O(c)$, ~ $c \in \Z$, ~ and

(4) ~ $\S$ is generated by $\S_1$, $\F$, and $\G$ as an $\O_{\p^1}$-algebra.

\noindent
    Namely, $\S_d$ is described as $\bigoplus_{i+2j+3k=d} S^i(\E)\otimes S^j(\O(b))\otimes S^k(\O(c))$ for $\E=\bigoplus_{i=0}^n \O(a_i)$, $a_0\leq a_1\leq \dots\leq a_n$.
    Denote by $\P[a_0, a_1, \dots a_n; b; c]$ the weighted projective space bundle $\proj\S$ over $\P^1$.

    Let $X = \P[0,a;b;c]$, and let $H_X$ and $F_X$ be the tautological \Q-divisor and a general fiber of $\pi : X \to \P^1$, respectively.
    Note that $6H_X$ is a Cartier divisor.
    Each smooth hypersurface $V$ linearly equivalent to $6H_X + 6kF_X$ is a del Pezzo fibration of degree 1 over $\P^1$ for $k \in \frac16\Z$.
    Note that any del Pezzo fibration $V$ of degree 1 over $\P^1$ can be constructed in this way (see (\ref{subsec:setup})):
$$\begin{array}{ccl}
  V & \Longarrow{\psi}{} & X = \P[0,a;b;c] \vspace{5pt}\\
  \Downarrow{\varphi} & \Swarrow{\pi} & \qquad V \in |6H + 6kF| \vspace{5pt}\\
  \P^1
\end{array}$$
    Assume that $V$ is a weak Fano $3$-fold with only finite $K$-trivial curves and with $\rho(V)=2$.

\newpar
    To prove Theorem(\ref{thm:df1}), it is suffice to fix the quadruple $(a, b, c, k)$ determining the weighted projective space bundle $X = \P[0, a; b; c]$ and the linear equivalence class of $V \sim 6H_X + 6kF_X$.

    Let $s_1$, $s_2$, and $s_3$ be sections of $\pi : X \to \P^1$ associated to surjections $\S \to S(\O)$, $\S \to S(\F)$, and $\S \to S(\G)$, respectively.
    Then, along the sections $s_2$ and $s_3$, $X$ has quotient singularities.
    These singularities are of type $\frac12(1,1,1)$ along $s_2$, and of type $\frac13(1,1,2)$ along $s_3$.
    Two rational ruled surfaces $L = \P[0,a] \to \P^1$ and $M = \proj{\cal M} \to \P^1$ are naturally considered as subbundles of $X$, where ${\cal M}$ is a graded $\O_{\p^1}$-subalgebra of $\S$ generated by $\F=\O(b)$ and $\G=\O(c)$.
    Then, $L \cong \Sigma_a$ contains $s_1$ as the minimal section, and $M \cong \Sigma_{|3b-2c|}$ contains $s_2$ and $s_3$ as disjoint sections.
    Let $f_L$ and $f_M$ be a general fiber of $L \to \P^1$ and of $M \to \P^1$, respectively.
    The restrictions $V|_L$ and $V|_M$ can be discribed by
$$
  V|_L \sim 6(s_1 + (a+k)f_L), \qquad
  V|_M \sim s_2 + 2(c+3k)f_M \sim s_3 + 3(b+2k)f_M.
$$
    From the effectivity of $V|L$, it follows that
\Eq(ak){
  a + k \geq 0.
}\\
    The smoothness of $V$ implies that $V$ is disjoint from singular loci $s_2 \cup s_3$ of $X$, hence
\Eq(bc){
  b = -2k, \qquad c = -3k,
}\\
and $M \cong \P^1 \times \P^1 \to \P^1$.
    Moreover, $V$ meets $M$ along a section $s = \{\mbox{a point}\} \times \P^1$ of $M \to \P^1$, and $s_2$ and $s_3$ are algebraically equivalent to $s$.
    Since $b$ and $c$ are integers, $k$ is also an integer by (\ref{eqn:bc}).

    Let $f \in H^0(X, \O_X(6H_X + 6kF_X))$ be a global section defining $V$ as its zero locus.
    There is a natural identification
\begin{eqnarray*}
  \lefteqn{H^0(X,\O_X(6H_X+6kF_X))} \\
  &\cong& \hspace{-5pt}
  H^0(\P^1, S^6(\E)(6k)) \op H^0(\P^1, S^4(\E)(4k)) \op H^0(\P^1, S^3(\E)(3k)) \\
  & & {} \op H^0(\P^1, S^2(\E)(2k)) \op H^0(\P^1, \E(k)) \op H^0(\P^1, \O) \op H^0(\P^1, \O),
\end{eqnarray*}
by (\ref{eqn:bc}).
    The smoothness of $V$ implies $\dim H^0(\O(a+6k)) > 0$ or $\dim H^0(\O(4k)) > 0$ or $\dim H^0(\O(3k)) > 0$, i.e.,
\Eq(neg){
  a + 6k \geq 0 \qquad {\rm or} \qquad k \geq 0.
}

\newpar
    Let $\rho : Y \to X$ be a minimal resolution of singularities of $X$, $E = \rho^{-1}(s_2)$, $P \cup G = \rho^{-1}(s_3)$, where $E$ and $P$ are $\P^2$-bundles over $\P^1$ and $G$ is a $\Sigma_2$-bundle over $\P^1$.
    Denote by $F$ a genenral fiber of $\pi\circ\rho : Y \to \P^1$.
    Then, the Picard group $\Pic Y$ is isomorphic to $\Z^{\op 5}$ and generated by $E$, $P$, $G$, $F$, and $H_0$, where $H_0$ is a divisor satisfing $\rho^{\ast} 6H_X \sim 6H_0 + 3E + 4P + 2G$.
    Since $V$ is disjoint from $s_2$ and $s_3$, $V$ can be regarded as a divisor of $Y$ and $V$ is linearly equivalent to $6H_0 + 3E + 4P + 2G + 6kF$.
    Denote by $H_V$ and $F_V$ the restrictions of $H_0$ and $F$ to $V$, respectively.
    The anti-canonical divisor $-K_V$ of $V$ is linearly equivalent to $H_V - (a+k-2)F_V$, because of $K_Y \sim - 7H_0 - 3E - 4P - 2G - (2-a-b-c)F$ and (\ref{eqn:bc}).

    Now consider the intersection numbers of a section $s$ of $M \cong \P^1 \times \P^1$ in $Y$.
    The section $s$ does not meet $E$, $P$ and $G$, hence $(s \cdot E)_Y = (s \cdot P)_Y = (s \cdot G)_Y = 0$.
    Since $s$ is regarded as a section of $\pi\circ\rho : Y \to \P^1$ by the natural embedding $M \subset Y$, it follows that $(s \cdot F)_Y = 1$.
    Since $V$ meets $M$ only along $s$ and since the minimal sections of $M \to \P^1$ are disjoint and algebraically equivalent to each other, we have $(V \cdot s)_Y = 0$, and hence $(H_0 \cdot s)_Y = k$.

    If $V$ is a weak Fano 3-fold, the anti-canonical divisor $-K_V$ is nef, and
$$
  (-K_V \cdot s) = (H_0 - (a+k-2)F \cdot s)_Y \geq 0.
$$
    This inequality implies that
$$
  2(1-k) \geq a.
$$
    Then, in the case $a + 6k \geq 0$ in (\ref{eqn:neg}), we have also $k \geq 0$ because $k$ is an integer.
    Nonnegativity of $a$ shows $k \leq 1$, hence $k = 0$ or $k = 1$.
    Thus the possible values of a pair $(k, a)$ are enumerated as follows :
$$
  (0, 0), \quad (0, 1), \quad (0, 2), \quad {\rm and} \quad (1, 0).
$$

\newpar
    Thus we have the four possibilities of quadruple $(a,b,c,k)$ as in Table \ref{tbl:df1}.
\begin{table}[htbp]
\caption{}\label{tbl:df1}
$$\begin{array}
{c@{\qquad}cccc@{\qquad}ll}
{\rm Nos.} & a & b & c & k & \quad V & \quad -K_V \\
 1 & 0 & 0 & 0 & 0 & 6H_X      & H_V+2F_V \\
 2 & 1 & 0 & 0 & 0 & 6H_X      & H_V+F_V \\
 3 & 2 & 0 & 0 & 0 & 6H_X      & H_V \\
 4 & 0 &-2 &-3 & 1 & 6H_X+6F_X & H_V+F_V \\
\end{array}$$
\end{table}
    In the rest of this section, we consider the realization of the each possibilities.

\newpar
    We now exclude two possibilities Nos.1 and 3 on Table \ref{tbl:df1}.
\ex{1}
    This case becomes a trivial fibration $\varphi : V = D_1 \times \P^1 \to \P^1$, where $D_1$ is a del Pezzo surface of degree $1$.
    The Picard number of $V$ is $10$, hence this case is excluded.
\ex{3}
    The case gives a weak Fano 3-fold $V$ in $X = \proj \S$, where $\S$ is generated by $\S_1 = \O \op \O(2)$, $\F = \O \subset \S_2$, and $\G = \O \subset \S_3$.
    The Mori cone $NE(V)$ has two edges, one of them is the extremal ray of type $D_1$ corresponding to our del Pezzo fibration of degree $1$.
    We will fix the generator of the other edge of $NE(V)$.
    Let ${\cal T}$ be a graded $\O_{\p^1}$-subalgebra generated by $\O (\subset \S_1)$, $\F$, and $\G$, and denote $\proj{\cal T}$ by $\Lambda$.
    Then, $\Lambda$ can be regarded as a divisor of $X$ through the natural surjection $\S \to {\cal T}$.
    The intersection $V$ and $\Lambda$ is a surface $S = V \cap \Lambda$ with the trivial fibration $S = C \times \P^1 \to \P^1$, where $C$ is an elliptic curve in a del Pezzo surface $V \cap F$ of degree $1$.
    Consider each section $s = \{\mbox{a point}\} \times \P^1$ of $S \to \P^1$ as a curve in $V$.
    Its normal bundle is ${\cal N}_{s/V} \cong \O \op \O(-2)$, i.e., the each section is non-isolated $K$-trivial curve.
    These curves generates the edge of $NE(V)$.
    Thus the case is excluded.

\newpar
    Next, we construct the del Pezzo fibrations with each conditions in the rest cases on Table \ref{tbl:df1}, and show Theorem(\ref{thm:df1}) and Supplement(\ref{thm:df1s}).
\cons{2}
    Considering $X = \P[0,1;0;0]$ as an ample divisor of $Y = \P[0,0,0;0;0]$, we can see that $V$ is an ample divisor of $W \cong B_1 \times \P^1$.
    Here $W$ is a general member of $|6H_Y|$ and $B_1$ is a Fano 3-fold of index $2$ with $\bigl(\frac{-K_{B_1}}{2}\bigr)^3 = 1$.
    Hence $V$ is a Fano 3-fold with $\rho(V) = 2$ by (\ref{thm:ampleLeff}), and has two extremal rays.
    One of them is the ray of type $D_1$ associated to the fibration considered here.
    The other ray $R$ is of type $E_1$ with $\psi = {\rm cont}_R : V \to B_1$ corresponding to the linear system $|6H_V|$.
    The exceptional divisor for $\psi$ is the unique member of $|H_V-F_V|$, which contracts to a curve $C$ of degree $\bigl(\frac{-K_{B_1}}{2} \cdot C\bigr) = 1$ with genus $g(C) = 1$.
    This is (\ref{thm:df1}.1).
\cons{4}
    The case gives a weak Fano 3-fold $V$ in $X = \proj \S$, where $\S$ is generated by $\S_1 = \O \op \O$, $\F = \O(-2) \subset \S_2$, and $\G = \O(-3) \subset \S_3$.
    As in the above argument, consider the rational ruled surface $M = \proj{\cal M} \cong \P^1 \times \P^1 \to \P^1$ corresponding to a graded $O_{\p^1}$-submodule ${\cal M}$ generated by $\F$ and $\G$.
    Then, $V$ meets $M$ along $s = \{\mbox{a point}\} \times \P^1$.
    The curve $s$ generates an edge of the Mori cone $NE(V)$.
    To see this, consider $V$ and $X$ as subvarieties in the $\P^{22}$-bundle $\P(\S'_6) = \P[0^2,1^2,2^3,3^4,4^5,6^7]$ over $\P^1$ through the standard isomorphisms $X = \proj \S \cong \proj \S' \cong \proj \S'^{(6)} \subset \P(\S'_6)$.
    Here $\S' = \bigoplus_{d \geq 0} \S'_d$ is a graded $\O_{\p^1}$-algebra defined by $\S'_d = \S_d \otimes \O(d)$, and $\S'^{(6)} = \bigoplus_{d \geq 0} \S'_{6d}$.
    The ruled surface $M = \proj{\cal M} \to \P^1$ is regarded as in $\P(\S'_6)$ by a surjection $\S_6 \otimes \O(6) \to {\cal M}_6 \otimes \O(6) \to 0$, i.e., $\S'_6 \to \O\op\O \to 0$.
    The tautological line bundle $H'$ of $\P(\S'_6)$ gives a horizontal ruling of $M \cong \P^1\times\P^1 \to \P^1$.
    The anti-canonical divisor $-K_V \sim H_V + F_V$ is the restriction of $\frac16H'$.
    The equality $(C \cdot -K_V)_V = (C \cdot \frac16H')_{\p(\S'_6)} = 0$ for a curve $C\subset V\subset\P(\S'_6)$ implies that $C$ is a section of the ruled surface $M \to \P^1$, because $H'$ gives a horizontal ruling of $M$.
    Thus the curve $s$ is the unique generator of the edge of $NE(V)$.
    Moreover, the normal bundle ${\cal N}_{s/V}$ is isomorphic to $\O(-1) \op \O(-1)$, i.e., $s$ is the isolated $K$-trivial curve in $V$; hence there is the $(-F_V)$-flop $\chi : V \dots\to V'$ along $s$.
    The new 3-fold $V'$ is also a weak Fano 3-fold having an extremal ray of type $D_1$, and its contraction morphism $\varphi' : V' \to \P^1$ is a del Pezzo fiberation of degree 1.
    The composition map $\varphi'\circ\chi : V \dots\to \P^1$ is defined by a linear system $|H_V|$.
    We can see that $\rho(V) = 2$ by (\ref{thm:ampleLeff}), considering $X \cong \proj \S' = \P[1,1;0;0]$ as an ample effective divisor linearly equivalent to $H_Y+2F_Y$ of $Y = \P[0,0,0;0;0]$, and $V \subset X$ as an ample effective divisor of $B_1 \times \P^1 \subset Y$, where $B_1$ is a Fano 3-fold of index 1 with $\bigl(\frac{-K_{B_1}}{2}\bigr)^3 = 1$.
    The (pluri-)anti-canonical model $\ol{V}$, the image of the morphism defined by $|-6K_V|$, is the complete intersection of $B_1^4$ and $Q^\ast$ in the weighted projective space $\ol{X}=\P(1^4,2,3)$.
    Here $B_1^4$ is a weighted hypersurface of degree 6 in $\ol{X}$, i.e., the smooth del Pezzo 4-fold of degree 1, and $Q^\ast$ is the weighted hypersurface of degree 2 defined by a quadric in the coordinate functions of weight 1 on $\ol{X}$, $Q^\ast$ has the singularities along the line $L$ through the two singularities of $\ol{X}$.
    The model $\ol{V}$ has only one ODP, which is the intersection of $B_0^4$ with the singular locus $L$.
    Thus we derive (\ref{thm:df1}.2).


~\\
Kiyohiko Takeuchi, \\
School of Education, Sugiyama Jogakuen University, \\
Hosigaoka-motomoachi, Chikusa-ku, Nakoya, Aichi, 464-8662 Japan \\
e-mail: takeuchi@sugiyama-u.ac.jp

\end{document}